\documentclass[a4paper,10pt]{article}
\pdfoutput=1
\usepackage[utf8]{inputenc}
\usepackage[T1]{fontenc}
\usepackage[francais]{babel}
\usepackage[all]{xy}
\usepackage{graphicx}
\usepackage{psfrag}
\usepackage{eufrak}
\usepackage{pgf}
\usepackage{pinlabel}
\usepackage{placeins}
\usepackage{amsmath}
\usepackage{amsthm}

\usepackage[pdftex]{hyperref}

\usepackage[all]{hypcap}

\addtocontents{toc}{
  \protect\setlength{\baselineskip}{10pt}}

\hypersetup{linkcolor=blue,colorlinks=true,citecolor=magenta}

\newcommand{\beqa}{\begin{eqnarray}}
\newcommand{\eeqa}{\end{eqnarray}}
\newcommand{\beqan}{\begin{eqnarray*}}
\newcommand{\eeqan}{\end{eqnarray*}}  
\newcommand{\bmarr}{\beq\begin{array}}
\newcommand{\emarr}{\end{array}\eeq} 
\newcommand{\ssi}{si et seulement si }

\newcommand{\pr}{{\it Preuve. }}

\newcommand{\esp}{~~~~~} 
\newcommand{\cad}{c'est-\`a-dire } 
\newcommand{\cadab}{c.-à-d.\ }

\def\re{\mathop{\rm R\acute{e}}\nolimits}
\def\im{\mathop{\rm Im}\nolimits}
\def\aut{\mathop{\rm Aut}\nolimits}

\def\rang{\mathop{\rm rang}\nolimits}
\def\tr{\mathop{\rm Tr}\nolimits}
\def\endo{\mathop{\rm End}\nolimits}
\def\hom{\mathop{\rm Hom}\nolimits}

\def\pgcd{\mathop{\rm pgcd}\nolimits}
\def\pf{\mathop{\rm Pf}\nolimits}
\def\sym{\mathop{\rm Sym}\nolimits}
\def\id{\mathop{\rm Id}\nolimits}

\newcommand{\C}{{\bf C}}
\newcommand{\D}{{\bf D}}

\newcommand{\G}{{\bf G}}
\newcommand{\GL}{{\bf GL}}

\newcommand{\HH}{{\bf H}}
\newcommand{\K}{{\bf K}}

\newcommand{\N}{{\bf N}}
\newcommand{\OO}{{\bf O}}

\newcommand{\Q}{{\bf Q}}
\newcommand{\R}{{\bf R}}
\newcommand{\SL}{{\bf SL}}
\newcommand{\PSL}{{\bf PSL}}
\newcommand{\SO}{{\bf SO}}

\newcommand{\SU}{{\bf SU}}
\newcommand{\PSU}{{\bf PSU}}
\newcommand{\Sp}{{\bf Sp}}

\newcommand{\SSS}{{\bf S}}

\newcommand{\U}{{\bf U}}
\newcommand{\Z}{{\bf Z}}

\newcommand{\al}{\alpha}
\newcommand{\be}{\beta}
\newcommand{\ga}{\gamma} \newcommand{\Ga}{\Gamma}
\newcommand{\de}{\delta} \newcommand{\De}{\Delta}
 \newcommand{\eps}{\epsilon}
 \newcommand{\la}{\lambda}
\newcommand{\La}{\Lambda} \newcommand{\ro}{\rho}   
\newcommand{\sig}{\sigma}  
\newcommand{\tht}{\theta}
 \newcommand{\ph}{\varphi}
\newcommand{\om}{\omega} \newcommand{\Om}{\Omega}

 \newcommand{\mcB}{{\mathcal B}}
\newcommand{\mcC}{{\mathcal C}} \newcommand{\mcD}{{\mathcal D}}
\newcommand{\mcE}{{\mathcal E}} \newcommand{\mcF}{{\mathcal F}}
\newcommand{\mcG}{{\mathcal G}} \newcommand{\mcH}{{\mathcal H}}
\newcommand{\mcI}{{\mathcal I}} 
 
 \newcommand{\mcO}{{\mathcal O}}
\newcommand{\mcP}{{\mathcal P}} 
 \newcommand{\mcS}{{\mathcal S}}
\newcommand{\mcT}{{\mathcal T}} 
\newcommand{\mcV}{{\mathcal V}} \newcommand{\mcW}{{\mathcal W}}

\newcommand{\mfA}{{\mathfrak A}} \newcommand{\mfa}{{\mathfrak a}}
\newcommand{\mfB}{{\mathfrak B}} 
\newcommand{\mfF}{{\mathfrak F}} \newcommand{\mff}{{\mathfrak f}}
\newcommand{\mfH}{{\mathfrak H}}
\newcommand{\mfP}{{\mathfrak P}}
\newcommand{\mfQ}{{\mathfrak Q}}

\newcommand{\mfS}{{\mathfrak S}}

\newcommand{\smallmat}[4]{\renewcommand{\arraystretch}{0.5}
\renewcommand\arraycolsep{0.2ex}
({\scriptsize \begin{array}{cc}  #1 & #2 \\ #3 & #4
\end{array}})\renewcommand{\arraystretch}{1}}
\newcommand{\smallvect}[2]{\renewcommand{\arraystretch}{0.5}
({\scriptsize \begin{array}{cc}  #1  \\ #2
\end{array}})  \renewcommand{\arraystretch}{1}}

\newcommand{\subsetsim}{\subset_{\rm eq}}

\def\thebibliography#1{\section*{Bibliographie}\list
 {[\arabic{enumi}]}{%
 \settowidth\labelwidth{[#1]}\leftmargin\labelwidth
 \advance\leftmargin\labelsep
 \usecounter{enumi}}
 \def\newblock{\hskip .11em plus .33em minus .07em}
 \sloppy\clubpenalty4000\widowpenalty4000
 \sfcode`\.=1000\relax}

\theoremstyle{plain}
\newtheorem{theo}{Théorème}[section]
\newtheorem*{theo*}{Théorème}
\newtheorem{prop}[theo]{Proposition}
\newtheorem{coro}[theo]{Corollaire}
\newtheorem{lemm}[theo]{Lemme}
\newtheorem{affi}[theo]{Affirmation}
\newtheorem*{affi*}{Affirmation}

\theoremstyle{definition}
\newtheorem{defi}[theo]{Définition} 
\newtheorem*{defi*}{Définition} 

\theoremstyle{remark}
\newtheorem{rema}[theo]{Remarque}
\newtheorem{remas}[theo]{Remarques}
\newtheorem{exem}[theo]{Exemple}

\newtheorem*{exem*}{Exemple}
\newtheorem*{exems*}{Exemples}

\newtheorem{tabl}{Table}

\pgfdeclareimage[interpolate=true,height=4cm]{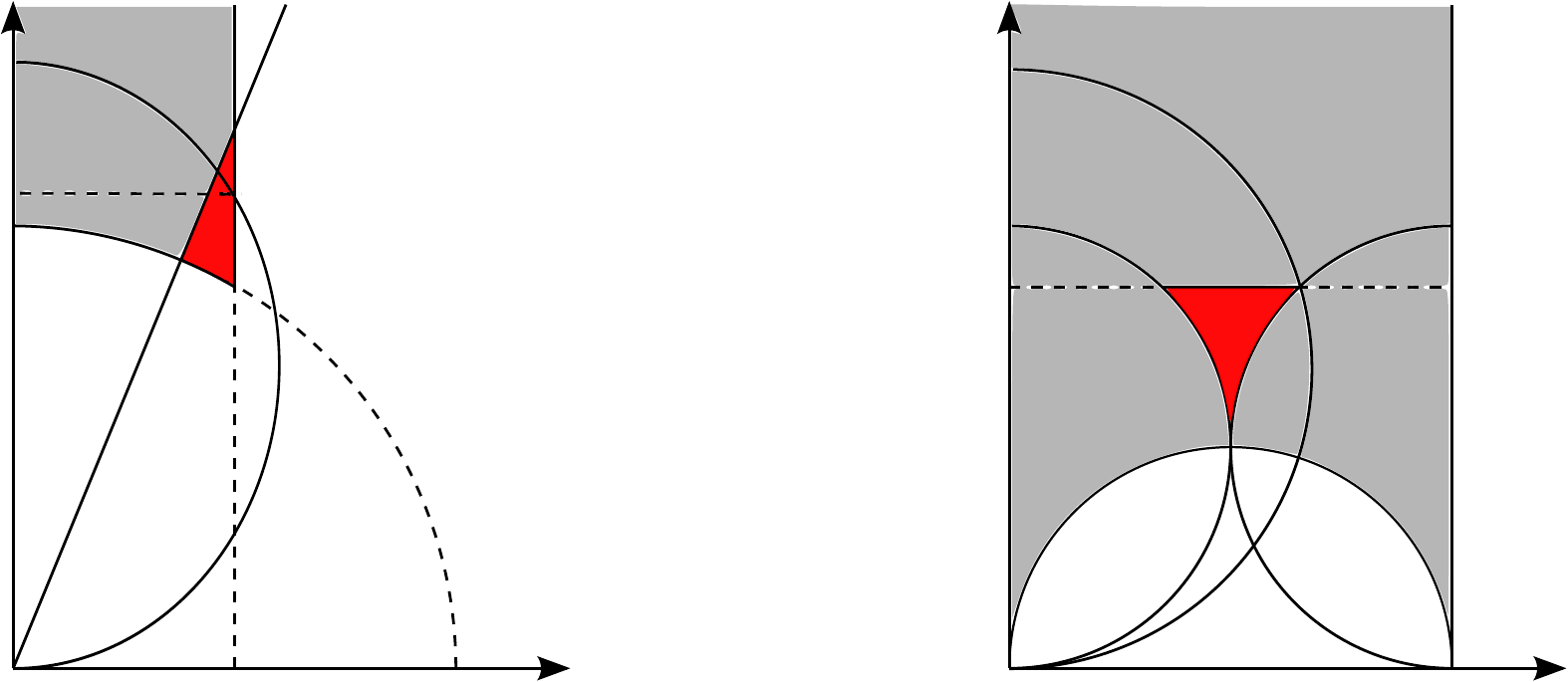}{I1F4_art}
\pgfdeclareimage[interpolate=true,height=4cm]{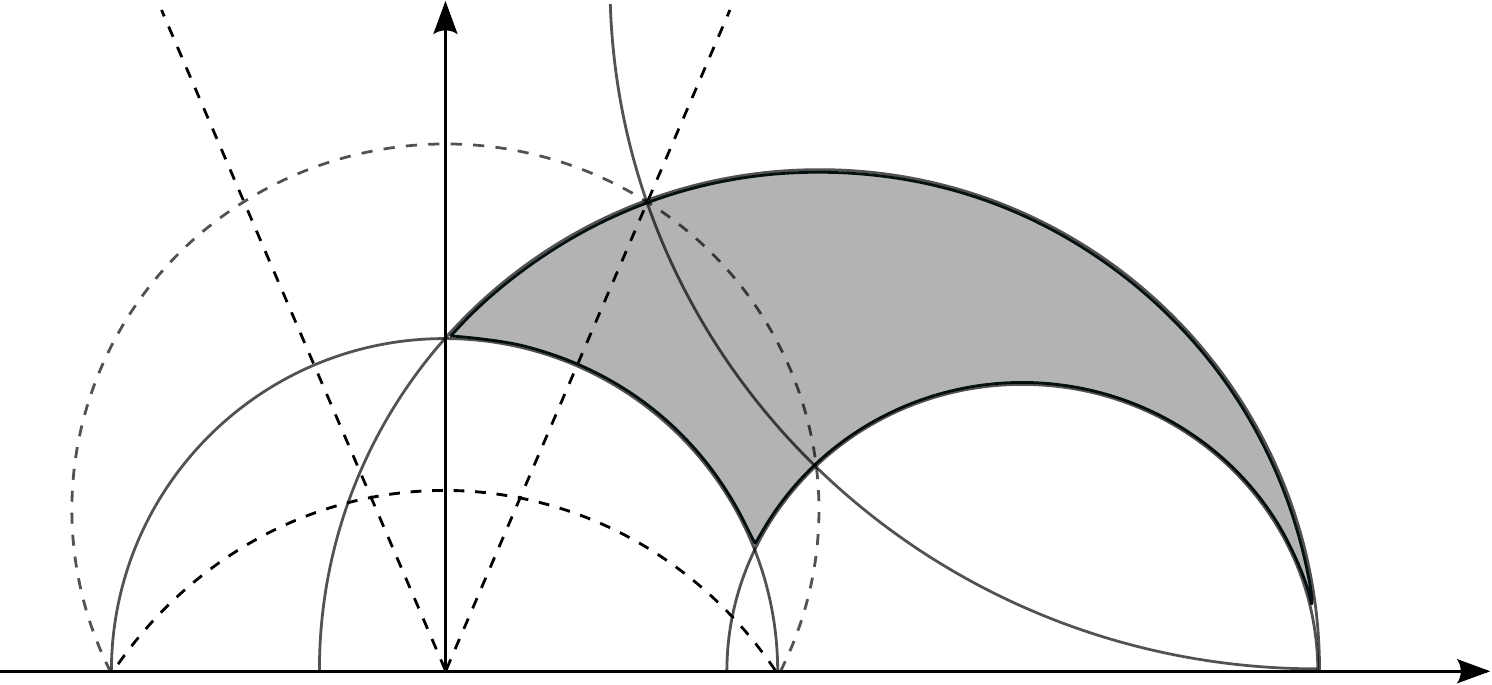}{I2mG3_art}
\pgfdeclareimage[interpolate=true,height=4cm]{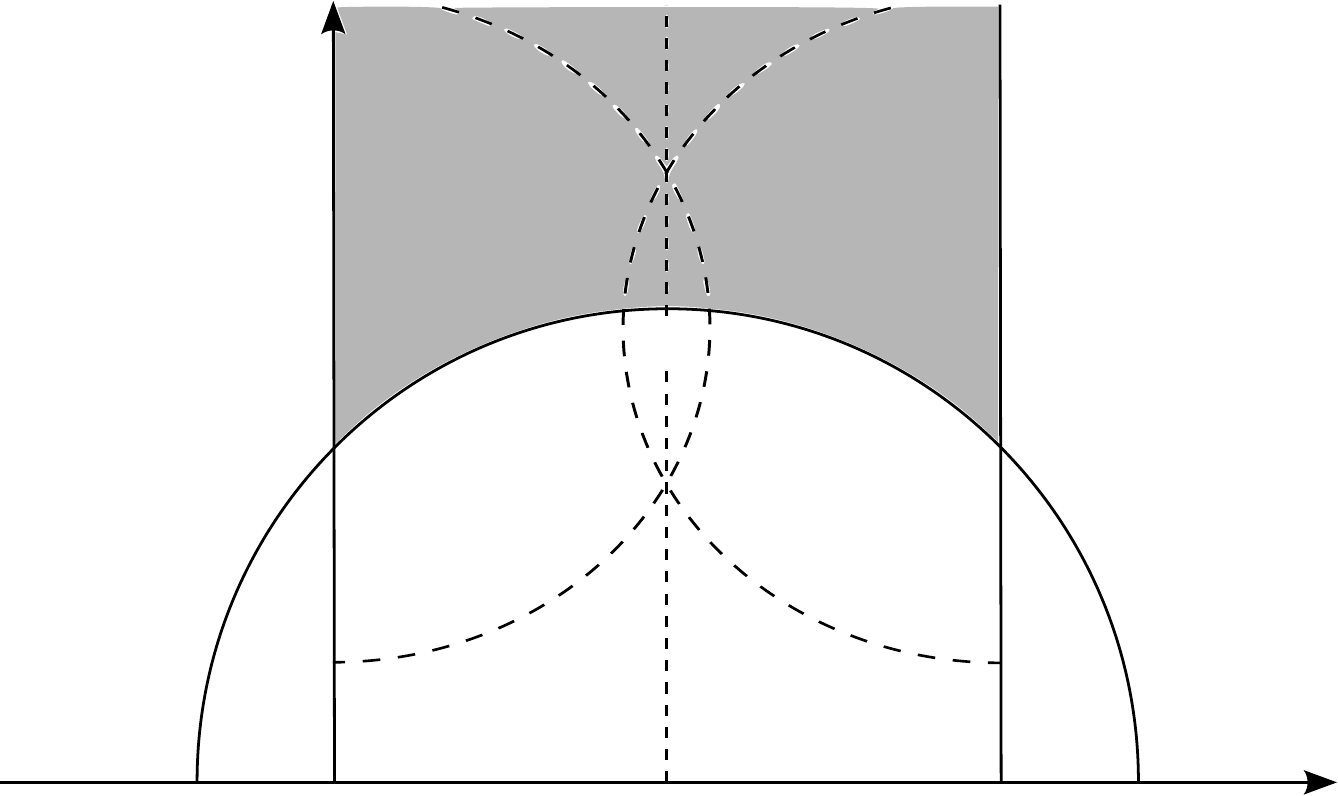}{G3H3_art}
\pgfdeclareimage[interpolate=true,height=4cm]{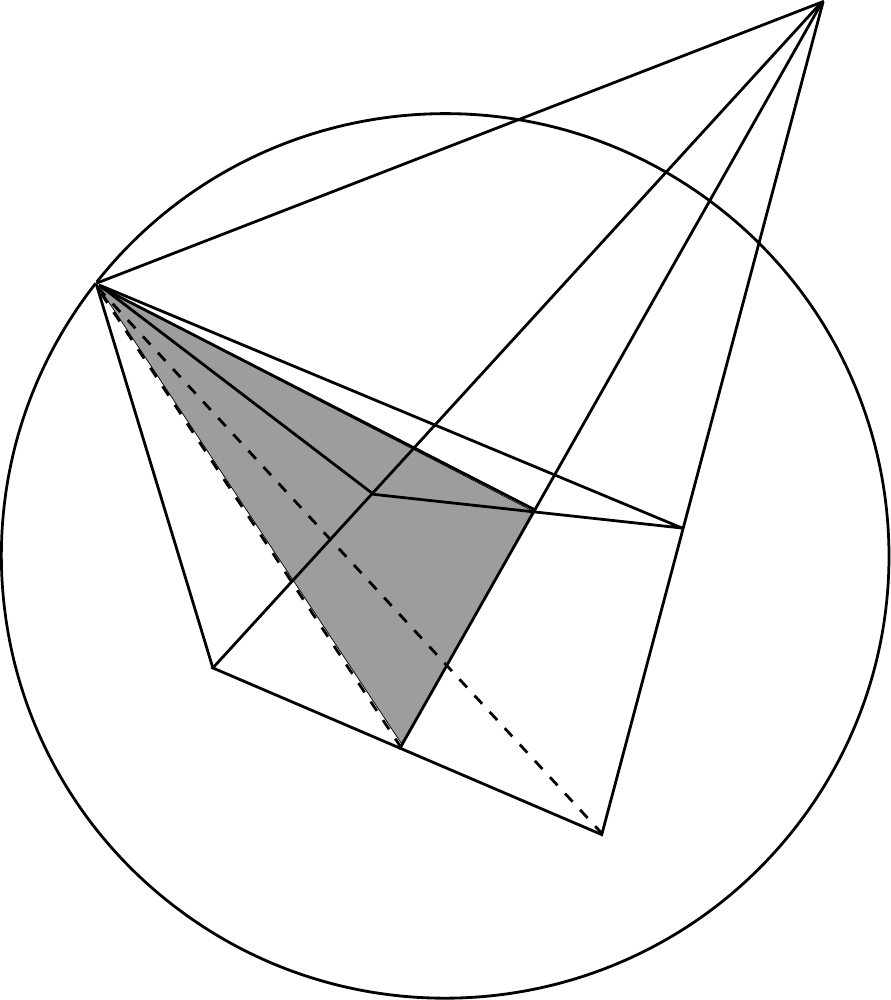}{I3mG3_art}

\begin{document}

\def\contentsname{Contenu de l'article}
\renewcommand{\thesection}{\arabic{section}}
\renewcommand{\abstractname}{}

\title{Isodualité des réseaux euclidiens\\
en petite dimension}
\author{Christophe \sc{BAVARD}}
\maketitle

\begin{abstract}
{\sc Résumé.}  Nous proposons une classification algébrique et géométrique
des réseaux euclidiens isoduaux de rang fixé. Premièrement, nous prouvons
que ces réseaux se répartissent suivant un nombre fini de types
algébriques. Deuxièmement, nous montrons qu'ils sont paramétrés par des
espaces symétriques associés aux groupes classiques $\SO_0(p,q)$,
$\Sp(2g,\R)$ et $\SU(p,q)$, en nombre fini.  Nous obtenons une description
complète des types algébriques et des matrices de Gram des réseaux isoduaux
jusqu'en rang 7.  Nous discutons également la question de la densité
maximale.
\end{abstract}

\begin{abstract}
{\sc Abstract.}  We propose an algebraic and a geometric classification of
euclidean isodual lattices of fixed rank. First, we prove that these
lattices are distribued according to a finite number of algebraic types.
Second, we show that they are parametrized by a finite number of symmetric
spaces associated to the classical groups $\SO_0(p,q)$, $\Sp(2g,\R)$ and
$\SU(p,q)$. We obtain a complete discription of algebraic types and Gram
matrices of isodual lattices up to rank 7. The maximal density problem is
also discussed.
\end{abstract}

\section{Introduction et résultats principaux}
Un réseau $\La$ d'un espace euclidien est {\em isodual} s'il est
isométrique à son dual~$\La^*$ (voir \cite{Conw-Sloa1994a}).  Tous les
réseaux euclidiens de rang~1 ou 2 sont isoduaux, à homothétie près.  En
dimension 3, Conway et Sloane \cite{Conw-Sloa1994a} décrivent les réseaux
isoduaux et déterminent leur densité maximale ; la description de
\cite{Conw-Sloa1994a} s'appuie sur un système de paramètres liés à la
cellule de Voronoï, particulièrement bien adapté à la dimension 3 (voir
\cite{Conw-Sloa1992a}).  

\par Parmi les réseaux isoduaux, on trouve les réseaux $\La$ engendrés par
les bases symplectiques d'un espace hermitien positif ($\La$ et $\La^*$
sont alors isométriques par la structure complexe). Il s'agit des réseaux
{\em symplectiques} (\cite{Buse-Sarn1994a,Bavard2001a}), famille importante
comprenant les réseaux de périodes des surfaces de Riemann, ou jacobiennes.
Une autre classe intéressante de réseaux isoduaux est introduite dans
\cite{Berg-Mart1995a} : les réseaux {\em orthogonaux} (voir ci-dessous).
Bergé et {\mbox Martinet} \cite{Berg-Mart1995a} caractérisent à la Voronoï
les maxima locaux de densité parmi les réseaux symplectiques ou orthogonaux
; en dimension 4 ils étudient le lien entre l'isodualité et les
configurations de vecteurs minimaux et classent les réseaux isoduaux de
rang~4 ayant au moins 7 vecteurs minimaux.

\par 
Un premier objectif est ici de proposer une double classification 
des réseaux isoduaux : algébrique et géométrique. 
Pour cela, nous appellerons   {\it réseau isodual}  
un couple $(\La,\sig)$ où $\La$ est un réseau d'un espace euclidien $E$
et $\sig$ est une isométrie de  $E$ telle que $\sig(\La) = \La^*$. 
La donnée supplémentaire de l'isométrie $\sig$ permet de définir, 
comme dans \cite{Berg-Mart1995a}, une forme bilinéaire sur $E$  par 
\begin{equation} \label{equa:intro_alpha}
\al(u,v)= \langle u, \sig(v) \rangle \esp (u, v \in E),
\end{equation}
où $\langle\cdot,\cdot\rangle$ désigne le produit scalaire de $E$.
Nous appellerons  {\em type algébrique} de $(\La,\sig)$  la classe 
d'isométrie du $\Z$-module bilinéaire régulier $(\La,\al_{|\La \times \La})$
(définition~\ref{defi:ta_tg}). La forme $\al$ n'a pas de symétrie particulière
{\it a priori}. Le cas symétrique correspond aux réseaux orthogonaux, 
le cas antisymétrique aux réseaux symplectiques. 
Cependant $\al$ n'est pas arbitraire ; nous donnons une caractérisation
algébrique des formes $\al$ provenant des réseaux isoduaux 
(théorème~\ref{theo:intro_finitude}, assertion~(1)). 
\par
Pour aborder l'aspect géométrique,  considérons 
l'ensemble $P_n$ des matrices  de Gram des réseaux euclidiens de rang~$n$ 
et de covolume~1.
Rappelons que $P_n$ est  l'espace symétrique (riemannien) associé au groupe
semi-simple $\SL_n(\R)$, voir \S\ref{subs:types_reels_aspects_diff}. 
Nous montrons que {\em les réseaux isoduaux de rang $n$
partageant un type isodual fixé sont paramétrés  par une 
sous-variété totalement géodésique de $P_n$}. Plus précisément,
si  $(\La,\sig)$  est muni d'une base, 
nous dirons que la matrice  $F\in \GL_n(\Z)$ de la forme $\al$ 
{\em représente le type isodual de $(\La,\sig)$}.
Une vérification élémentaire montre alors que la matrice de Gram de 
$\La$ est localisée dans un  sous-ensemble $V_F$ de $P_n$, à savoir 
\begin{equation} \label{equa:intro_V_F}
V_F = \{A \in P_n; A F^\vee \! A = F\},
\end{equation}
où $F^\vee = F'\null^{-1}$ ($F'$ désigne la transposée de~$F$).
Il se trouve que  $V_F$ est une  sous-variété totalement  
géodésique de $P_n$  (proposition~\ref{prop:VF_symetrique}).
Nous définissons le {\em type géométrique} de  $(\La,\sig)$ comme 
l'orbite de $V_F$  dans $P_n$ sous l'action entière du groupe $\GL_n(\Z)$
(définition~\ref{defi:ta_tg}). 
Le point de départ de notre classification est le résultat suivant.
\begin{theo}
\label{theo:intro_finitude}
{\em (caractérisation algébrique et finitude des types).}
\begin{enumerate}
\renewcommand{\itemsep}{0mm}
\item Un élément $F\in GL_n(\Z)$ représente un type isodual
\ssi $FF^\vee$ est d'ordre fini.
\item L'ensemble des types algébriques isoduaux
de rang fixé est fini.
\item L'ensemble des types géométriques isoduaux
de rang fixé est fini.
\end{enumerate}
\end{theo}
L'assertion (1) donne  une caractérisation purement algébrique
des types isoduaux ; la condition nécessaire est évidente car 
$FF^\vee$ représente $\sig^2$ qui est un automorphisme de $\La$.
  L'assertion (2) contient la finitude des réseaux
unimodulaires entiers de rang fixé, définis ou indéfinis (cas symétrique). 
On notera qu'à partir du rang~2,  les classes d'isométries de $\Z$-modules 
bilinéaires réguliers ne sont pas en nombre fini, voir 
remarques~\ref{rema:finitude}. L'assertion (3) signifie qu'il existe 
une famille finie $(V_i)$ de sous-variétés de $P_n$ telle que tout 
réseau isodual de rang~$n$ admet une matrice de Gram appartenant 
à l'une des $V_i$. 
\par 
Le deuxième objectif de l'article est d'obtenir des informations 
sur la densité maximale des réseaux isoduaux (voir plus bas), ce qui 
nécessite d'avoir une description suffisamment précise de ces 
réseaux. 
Des paramétrages explicites se déduisent du fait que les 
sous-variétés $V_F$ définies par l'équation \eqref{equa:intro_V_F}
sont des espaces symétriques d'un type bien particulier. 
\begin{theo}[structure et paramétrage des types géométriques]
\label{theo:intro_structVF}
Tout type géométrique isodual est (représenté par) un produit riemannien
d'espaces symétriques associés aux groupes $\SO(p,q)$, 
$\Sp(2g,\R)$ ou $\SU(p,q)$. 
\end{theo}
La structure métrique des sous-variétés $V_F$  (voir plus précisément 
\S\ref{subs:composantes}) ne dépend que du {\em type algébrique  réel}, 
notion introduite  au \S\ref{subs:types_reels_aspects_diff} à partir de 
la propriété  caractéristique (1) du théorème~\ref{theo:intro_finitude}.
Les espaces symétriques des groupes $\SO(p,q)$, 
$\Sp(2g,\R)$ ou $\SU(p,q)$ admettent comme  modèles  des ouverts
d'espaces vectoriels, par exemple le demi-espace de Siegel 
pour $\Sp(2g,\R)$. En explicitant le type algébrique réel de $F$,
 on en déduit facilement  un paramétrage de $V_F$
par un ouvert  de $\R^N$, $N$ étant la dimension de $V_F$ (voir
\S\ref{subs:param} pour les détails). Ainsi, {\em il est possible
d'écrire explicitement les matrices de Gram de n'importe quelle
 famille  de réseaux isoduaux dont le type algébrique est donné}.
\par 
Précisons quelques aspects algébriques. 
Nous définissons  {\em l'ordre} d'un 
réseau isodual $(\La,\sig)$ (ou de son type algébrique) comme 
 l'ordre de $\sig^2$.  Selon une observation de Bergé et Martinet
\cite[p. 315]{Berg-Mart1995a}, le couple $(\La,\sig^m)$ est encore 
isodual pour tout $m$ impair. Si l'on s'intéresse seulement au réseau 
$\La$ (par exemple à sa densité), on peut donc se restreindre aux
types algébriques dont l'ordre est une puissance de 2, que nous
appellerons {\em types principaux}. L'étape préliminaire à  notre 
classification %
consiste à faire l'inventaire des éléments d'ordre fini de $\GL_n(\Z)$. 
Pour les cas traités ici ($n$ petit, types principaux), il s'agit
d'un simple exercice à partir de la théorie des représentations 
entières des groupes cycliques, voir \S\ref{subs:of_GLnZ}.
Ensuite, étant donné $R\in \GL_n(\Z)$ d'ordre fini, la détermination 
des $F \in \GL_n(\Z)$ tels que $FF^\vee = R$ -- s'il en existe -- 
s'appuie sur la connaissance des modules quadratiques (resp. 
hermitiens)  sur $\Z$ (resp.  sur  certains anneaux d'entiers), 
voir \S\ref{subs:method_class} et 
\S\S\ref{subs:typalg1-4}-\ref{subs:typalg7} ; en particulier,
la classification des formes quadratiques  entières de petit déterminant 
est cruciale  pour déterminer  la  \og composante symétrique \fg\ des 
types algébriques  (proposition~\ref{prop:dec_can}).
 Il serait intéressant de caractériser  les éléments $R$ d'ordre fini 
de la  forme $R=FF^\vee$ avec $F \in \GL_n(\Z)$. Ce problème est discuté 
au \S\ref{susb:existence} :  nous y donnons des conditions nécessaires,
 mais la caractérisation reste ouverte.
\par 
Revenons maintenant à la question de la densité maximale. Soulevée
dans \cite{Buse-Sarn1994a} pour les réseaux symplectiques, elle se 
pose également  pour  les réseaux isoduaux en général.  Soit
$\ga^\mathrm{isod}_n$  (resp. $\mu_F$) le maximum de l'invariant d'Hermite 
sur les réseaux isoduaux (resp. sur  $V_F$). Il s'agit bien de maxima car 
les $V_F$ vérifient un 
critère de compacité de Mahler (proposition~\ref{prop:VF_Mahler})
et $\ga^\mathrm{isod}_n$ est le maximum des $\mu_F$, en nombre fini
(théorème~\ref{theo:intro_finitude}, (3)).
Pour $n \in \{1,2,4,8,24\}$, on a $\ga^\mathrm{isod}_n=\ga_n$ (constante 
d'Hermite usuelle). Grâce à Conway et Sloane \cite{Conw-Sloa1994a}
on sait aussi que $\ga^\mathrm{isod}_3 = 1/2+1/\sqrt{2}$. Ce maximum
est réalisé par un réseau orthogonal  lorentzien. 
La densité maximale du type orthogonal lorentzien est connue
jusqu'en rang~12 d'après \cite{Bavard2007a}.
Dans \cite{Buse-Sarn1994a}, Buser et Sarnak montrent
que la densité maximale des réseaux symplectiques croît comme le rang,
tandis que celle des jacobiennes croît comme le logarithme
du rang.
\par
Les dimensions $n=5$, $6$ et $7$ présentent un intérêt particulier 
car  $\ga^\mathrm{isod}_n < \ga_n$, les réseaux les plus denses n'étant
pas isoduaux. Nous avons donc choisi de classer les types isoduaux jusqu'en
dimension 7. La classification complète des types de rang $n$ fixé est utile
même si $\ga^\mathrm{isod}_n$ est connue, et ce pour au moins
deux raisons. D'abord, ces types 
interviennent dans les types décomposés des dimensions supérieures. 
Ensuite, la connaissance du type algébrique $F$ donne une estimation
plus précise de la densité {\it via} $\mu_F$. Pour ne pas allonger 
l'article, nous nous limitons aux types principaux à partir de la 
dimension~5. Les classifications algébrique et géométrique des types 
non principaux pour \mbox{$5 \leq n \leq 7$} 
ainsi que la classification des types
en dimension 8 peuvent être abordées par la même méthode.
\par
La détermination de $\ga^\mathrm{isod}_n$ nécessite une étude spécifique
de chaque type géométrique (ou au moins des types maximaux, 
voir ci-dessous).
Pour tous les types $V_F$ classés ici, nous estimons $\mu_F$ en
exhibant {\em un maximum local de densité} ; nous déterminons
autant que possible la valeur de $\mu_F$, ce qui est relativement 
facile dans certains cas ($V_F$  de dimension $\leq 3$, type décomposé), 
voir  tables~\ref{tabl:types_geom_1_3}-\ref{tabl:types_geom_6} et
\ref{tabl:types_geom_7}. Nous obtenons ainsi la valeur de $\mu_F$
pour tous les types de rang $\leq 4$ et pour tous les types 
principaux de rang $5$ à l'exception du type orthogonal de signature
$(3,2)$, avec comme conséquence que $\ga^\mathrm{isod}_5 = \mu_{3,2}$
(corollaire~\ref{coro:ga5isod}).
Les types géométriques sont naturellement ordonnés par l'inclusion.
Par exemple toute sous-variété  $V_F$ est incluse dans une sous-variété
$V_G$ où  $G$ est un type algébrique principal. 
Nous donnons un critère d'inclusion 
(proposition~\ref{prop:inclusion}) permettant d'ordonner
les types  géométriques (tables~\ref{tabl:types_geom_1_3},
\ref{tabl:types_geom_4}, \ref{tabl:types_geom_5},
\ref{tabl:incl_tg6} et \ref{tabl:incl_tg7}). 
En particulier, nous pouvons décrire  explicitement 
 tous les types géométriques maximaux jusqu'en en dimension 7.
\begin{theo}[matrices de Gram des réseaux isoduaux de rang $\leq 7$]
\label{theo:intro_param1_7}
Tout réseau euclidien isomorphe à son dual et de rang $n\leq 7$ 
admet une matrice de Gram appartenant à l'une 
des sous-variétés de la table suivante. 
\begin{table}[h]
\centering
\renewcommand{\arraystretch}{1.2}
\begin{tabular}{lllll}
\hline\vspace{-4pt}
n & orthogonaux ou &  autres & décomposés  & \#  \\
  & symplectique & indécomposables & & \\
\hline
1 & $\{I_1\}$ & & & 1  \\
2 & $\mfS_1$  & & & 1 \\ 
3 & $V_{2,1}$  & & $\{I_1\} \mfS_1$ & 2 \\ 
4 & $V_{3,1}, V_{2,2}, V_{2,2}^{\rm II}, \mfS_2$ & & & 4\\
5 & $V_{4,1}, V_{3,2}$ & $V_{I_1F_4}$ & $I_{2,1} \mfS_1$,
$\{I_1\} \mfS_2$ & 5\\
6 & $V_{5,1}, V_{4,2},V_{3,3}, V_{3,3}^{\rm II},$ & 
$V_{I_2 F_4}, V_{I_{1,1}F_4}, V_{U_2F_4}$ & 
$V_{3,1} \mfS_1, V_{2,2} \mfS_1$, &  \\
  & $ \mfS_3$ & & $V_{2,2}^{\rm II} \mfS_1$ & 11 \\
7 & $V_{6,1}, V_{5,2}, V_{4,3}$ & 
$V_{I_3 F_4}, V_{I_{2,1} F_4}, V_{I_1 J_2 F_4}$,  & $V_{4,1}\mfS_1, V_{3,2}\mfS_1,$ & \\
  & & $V_{I_1 G_3 G_3^-}, V_{K_4 G_3}$ & $V_{2,1}\mfS_2, \{I_1\}\mfS_3$ & 12 \\
\hline
\end{tabular}
\caption{Types géométriques maximaux en dimension $n\leq 7$}
\label{tabl:intro_types_geom}
\end{table}
\end{theo}
Les types orthogonaux ou symplectique sont toujours maximaux (sauf 
si $\dim V_F \leq 1$, voir proposition~\ref{prop:orth_symp_max}, (2)).
Dans la table~\ref{tabl:intro_types_geom}, les sous-variétés
$\mfS_g$, $V_{p,q}$, $V_{2,2}^{\rm II}$ et $V_{3,3}^{\rm II}$
sont données  par les équations 
\eqref{equa:Sig_g}, \eqref{equa:Vpq}, \eqref{equa:V22_pair} et 
\eqref{equa:V33_pair}  respectivement. et les autres types 
indécomposables (deuxième colonne) sont explicités
au \ref{subs:ann_tg_non_scin} de l'annexe. 
Enfin, pour les types décomposés (définition~\ref{defi:tg_scinde}) 
 on pose $VW=V\oplus W$, sous-variété associée aux  sommes  
orthogonales des  réseaux des deux familles
(voir proposition~\ref{prop:som_prod_se_sym} pour sa 
structure métrique). 
\par
Il résulte du théorème~\ref{theo:intro_param1_7} que tout réseau euclidien
de rang $4$ isomorphe à son dual  admet une structure 
orthogonale ou symplectique.  En rang 3, nous retrouvons les deux
familles de \cite{Conw-Sloa1994a}.

\tableofcontents

\section{Types associés aux réseaux isoduaux}
\label{sect:types}
\subsection{Type algébrique, type géométrique}
\label{subs:def_types}
Soit $E$ un espace euclidien. Nous appellerons {\it réseau isodual} un 
couple $(\La,\sig)$ où $\La$ est un réseau de $E$ et $\sig$ est une 
isométrie de $E$ telle que $\sig(\La) = \La^*$. On a nécessairement 
$\det \La = 1$ et  $\sig(\La^*) = \La$.
Deux réseaux isoduaux $(\La_i,\sig_i)$ de $E_i$ ($i=1,2$) sont
isométriques  s'il existe une isométrie $\ph$ de $E_1$ sur $E_2$ 
telle que $\ph(\La_1)= \La_2$ et $\ph \sig_1 \ph^{-1} = \sig_2$. 
\par
Soit $(\La,\sig)$ un réseau isodual de $E$. Alors la forme bilinéaire
\begin{equation} \label{equa:alpha}
\al(u,v)= \langle u, \sig(v) \rangle \esp (u, v \in E),
\end{equation}
est par définition entière en restriction à $\La \times \La$.
\'Etant donnée une base $\mcB$ de~$\La$, notons $F \in \GL_n(\Z)$ 
la matrice de $\al$ dans $\mcB$ ($n=\dim E$) et $B$ la matrice de
Gram de $\mcB$ ; il est facile de voir que $F$ et $B$ sont liées par
la relation $B F^\vee \! B = F$.
De plus, les classes d'équivalence entière de $F$ et de $B$ sont
indépendantes du choix de $\mcB$. Rappelons que le groupe $\GL_n(\Z)$
agit par équivalence sur les matrices carrées $M$ par
$T \cdot M = T M T'  
(T \in  \GL_n(\Z)),$
 et qu'il agit en particulier sur l'espace $P_n$ 
des matrices de Gram de déterminant~1. 
Posons, comme dans l'introduction, $F^\vee = {F'}^{-1}$ et 
$$%
V_F = \{A \in P_n; A F^\vee \! A = F\}.
$$%
\begin{defi}\label{defi:ta_tg} 
Soit $(\La,\sig)$ un réseau isodual de rang~$n$ et soit
$F \in \GL_n(\Z)$ comme ci-dessus. Nous appellerons\\
\indent (1) {\em type algébrique} de $(\La,\sig)$ la classe 
d'équivalence entière de $F$, notée $[F]$, \cad la  classe
d'isométrie du $\Z$-module bilinéaire $(\La,\al_{|\La \times \La})$,\\
\indent (2) {\em type géométrique} de $(\La,\sig)$ la classe 
d'équivalence entière de $V_F$ dans l'espace $P_n$ des matrices de 
Gram.
\end{defi} 
Le type algébrique est clairement invariant par isométrie de
réseaux isoduaux. Il en est de même du type géométrique puisque
$V_{T\cdot F} = T \cdot V_F$ pour  $T\in \GL_n(\Z)$.

L'ensemble $P_n$ des matrices de Gram a une structure géométrique
très riche : il s'agit d'un espace symétrique riemannien (voir
les détails au \S \ref{subs:types_reels_aspects_diff}).
Les sous-variétés connexes, complètes et totalement géodésique
de $P_n$ sont automatiquement simplement connexes et symétriques ; pour
abréger, nous les appellerons {\em \og sous-espaces symétriques\fg}. 
Dans ce contexte,  on sait (proposition~\ref{prop:VF_symetrique})
 que les sous-ensembles   $V_F$  sont des   sous-espaces
 symétriques  de $P_n$. L'action par équivalence de $\GL_n(\Z)$ sur 
$P_n$ étant isométrique (voir~\S \ref{subs:types_reels_aspects_diff}), 
la  structure métrique -- riemannienne -- de $V_F$ est un invariant
du type algébrique et du type géométrique.

\subsection{Caractérisation algébrique des types isoduaux}
\label{subs:carac_types}
On peut se demander à quelle condition une matrice donnée
$F \in \GL_n(\Z)$ représente un type algébrique de réseau
isodual, ou {\em \og type algébrique isodual\fg}.  Un premier
critère évident est la non vacuité de $V_F$~: %
si $A=PP' \in V_F$ ($P \in \GL_n(\R))$, alors dans
l'espace euclidien usuel le réseau $(P'\Z^n,P^{-1}FP^\vee)$
est isodual de type $[F]$ (voir \S \ref{subs:def_types}).
Voici maintenant  un critère purement algébrique.
\begin{prop} 
\label{prop:carac_types}
{\rm (caractérisation algébrique des types algébriques isoduaux)}
Un élé\-ment $F$ du groupe $\GL_n(\Z)$ représente un type algébrique
isodual \ssi $FF^\vee$ est d'ordre fini.
\end{prop}
\pr Supposons que $F$ représente un type algébrique isodual.
Si $A\in V_F$, on a alors $FF^\vee \! A (FF^\vee)' = A$ (traduction de 
$\sig^2 \in \aut \La$ pour $(\La,\sig)$ isodual) et la matrice $FF^\vee$
est d'ordre fini par  compacité du stabilisateur  $\GL_n(\R)_A$.
Réciproquement, posons $\pi_F = \smallmat{0}{F}{F^\vee}{0}$,
$C_n = \cup_{\la >0}\la P_n$ et considérons 
\renewcommand{\arraystretch}{0.6} 
\begin{equation}\label{equa:Wn}
W_n = \{\tau_A; A \in C_n\},\esp \mbox{o\`u} \esp
\tau_A = \left( \begin{array}{cc}
A & 0 \\
0 & A^{-1}
\end{array} \right).
\end{equation}
L'action de $\pi_F$ sur $P_{2n}$ laisse stable $W_n$ ; de plus
 $A \in V_F$ \ssi $\tau_A$ est fixe par $\pi_F$ (noter 
que la relation $AF^\vee A= F$ entraîne $\det A = 1$). 
On vérifie (voir \S 
\ref{subs:finitude_types}, preuve du théorème \ref{theo:finitude_types})
que $W_n$ est un sous-espace symétrique de $P_{2n}$.
Si $FF^\vee$ est d'ordre fini, $\pi_F$ l'est aussi et 
on conclut grâce au théorème de point fixe d'\'Elie Cartan 
que $V_F$ est non vide. 

\subsection{Finitude des types isoduaux}
\label{subs:finitude_types}
\begin{theo}\label{theo:finitude_types}
L'ensemble des types algébriques (resp. géométriques) isoduaux
de rang fixé est fini.
\end{theo}

\pr Nous utilisons les propriétés topologiques de certains sous-groupes
de $\GL_{2n}(\R)$ combinées avec celles des réseaux symplectiques
semi-eutactiques établies dans \cite{Bavard2005a}. Soit $\mfS_n$
le sous-espace symétrique de $P_{2n}$ formé par les matrices
de Gram des réseaux symplectiques. Noter que $\mfS_n$ est 
stable par l'action de  $\Om= \smallmat{I_n}{0}{0}{-I_n}$, par
conséquent  le lieu des points fixes $W_n = \mfS_n^\Om$
est un sous-espace symétrique de $\mfS_n$. Posons
$H(\R)=\{\de_P;P\in \GL_n(\R)$\} où $\de_P=\smallmat{P}{0}{0}{P^\vee}$,
$G(\R)=H(\R) \cup \{\pi_F;F \in \GL_n(\R)\}$, et considérons le groupe
\begin{equation}\label{equa:gamma}
G(\Z) = G(\R) \cap \SL_{2n}(\Z) =
\{\de_P;P \in \GL_n(\Z)\}  \cup \{\pi_F;F \in \GL_n(\Z)\}.
\end{equation}
L'équivalence entière entre deux matrices $F$ et $G$ de  $\GL_n(\Z)$
équivaut clairement à l'existence d'un élément $\de_P$
($P \in \GL_n(\Z)$) qui conjugue $\pi_F$ et $\pi_G$. La finitude
des types isoduaux algébriques (et par suite celle des types
géométriques) résulte donc de la propriété suivante.

\begin{affi}\label{affi:sous-groupes_finis}
Les sous-groupes finis de $G(\Z)$  sont en nombre fini modulo
con\-ju\-gaison par $H(\Z)=H(\R) \cap \SL_{2n}(\Z)$.
\end{affi}

Le groupe $G(\Z)$ agit sur $W_n$. Si $\Pi$ est un sous-groupe fini
de $G(\Z)$, le lieu des points fixes $W_n^\Pi$ est non vide 
par le théorème du point fixe d'\'E. Cartan (loc. cit.). Nous
allons montrer l'existence d'un compact $K$ de $W_n$ tel que pour tout
sous-groupe fini  $\Pi$ de $G(\Z)$, l'un des translatés 
$\ga \cdot W_n^\Pi$ ($\ga \in H(\Z)$) coupe $K$. Le conjugué
$\ga \Pi \ga^{-1}$ sera alors  inclus dans
$\{\al \in \GL_{2n}(\Z);\al(K) \cap K \neq \emptyset\}$, ensemble
fini puisque l'action de $\GL_{2n}(\Z)$ sur $P_{2n}$ est
propre et discontinue.\\
Observons d'abord que $H(\R)$ est \og pseudo-algébrique \fg\  au sens
de \cite[définition~2.1]{Bavard2005a}. En effet, soit $M_{2n}(\C)$
l'espace des matrices carrées complexes d'ordre $2n$, soit 
$\ph : \SL_{2n}(\C) \to \GL(M_{2n}(\C)^2)$ définie par
$$\ph(M)\cdot (X,Y)= (MXM',MYM^{-1}),$$
 et soit $v= (\smallmat{0}{I_n}{-I_n}{0},\Om)$ ; alors $\ph(\SL_{2n}(\Z))$
laisse stable  le réseau $M_{2n}(\Z)^2$ et 
 $H(\R) = \{M \in \SL_{2n}(\R) ;\ph(M).v=v\}$. En outre $H(\R)$
contient toutes les transvections du sous-espace symétrique $W_n$.
Par suite, si $G^\Pi(\Z)$ désigne le commutant de $\Pi$ dans
$G(\Z) \cap \SL_{2n}(\Z)$, l'application
\begin{equation}\label{equa:MalherWnPi}
G^\Pi(\Z) \backslash W_n^\Pi \to  \SL_{2n}(\Z) \backslash P_{2n}
\end{equation}
est propre d'après \cite[proposition~2.7, 3)]{Bavard2005a}. Soit $\mu$
l'invariant d'Hermite défini sur $P_{2n}$. Par le théorème
de compacité de Mahler, la restriction de $\mu$ à $W_n^\Pi$
admet donc un maximum en $A \in W_n^\Pi$. Un tel point $A$
est nécessairement semi-eutactique relativement à $W_n^\Pi$,
\cite[proposition~2.1 (2)]{Bavard1997a}, et aussi relativement à $\mfS_n$
(\cite[\S 2.5, lemme 2.8]{Bavard2005a}) puisque $W_n^\Pi$ est le
lieu des points de $\mfS_n$ fixes par $\Pi$ et $\Om$.
Par suite (\cite[\S 2.6, théorème 1, 1), cas symplectique]{Bavard2005a})
on a  $\mu(A) \geq 1$. Pour conclure il suffit de remarquer que
$H(\Z) \backslash W_n \to  \SL_{2n}(\Z) \backslash P_{2n}$ est
propre (application \eqref{equa:MalherWnPi} pour $\Pi$ trivial) et de
prendre $K$ compact tel que
$\{A \in W_n;\mu(A) \geq 1 \}=\cup _{\ga \in H(\Z)}\ga \cdot K$.

\begin{remas} 
\label{rema:finitude}
(1) Le théorème \ref{theo:finitude_types} contient en  particulier la
finitude des réseaux entiers unimodulaires de rang fixé, positifs
ou indéfinis.\\
\indent(2) L'affirmation \ref{affi:sous-groupes_finis}
contient la finitude des sous-groupes
finis de $\GL_n(\Z)$ modulo conjugaison.\\
\indent(3) Les classes d'équivalence entières ne sont pas en nombre fini.
Par exemple, pour $a \in \N$, les matrices $F_a=\smallmat{1}{a}{0}{1}$
sont deux à deux non équivalentes car les polynômes caractéristiques
 de $F_aF_a^\vee$ sont distincts (cf. \S \ref{subs:carac_types}).
\end{remas}

\subsection{Types principaux. Types géométriques maximaux}
\label{subs:ppaux_et_max}
Si $F \in \GL_n(\Z)$, la classe de conjugaison de $FF^\vee$ dans 
$\GL_n(\Z)$ est un invariant de $[F]$. Supposons maintenant que
$[F]$ soit un type isodual; l'ordre de $FF^\vee$ sera appelé {\em ordre}
de $[F]$. Pour tout $\ell \in \N$, la matrice $G=F(F^\vee F)^\ell$
représente un type isodual car $GG^\vee = (FF^\vee)^{2\ell +1}$. De
plus  $V_F$ est inclus dans $V_G$ (en termes de réseaux, si $(\La,\sig)$
est isodual, alors $(\La,\sig^{2\ell +1})$ l'est aussi).
\begin{defi} 
Un type  algébrique isodual est dit {\em principal} si son ordre
est une puissance de 2.
\end{defi}
Par extension nous appellerons aussi {\em principal} tout type 
géométrique $V_F$ associé à un type algébrique principal $[F]$.
Noter que l'on peut avoir $V_F=V_G$ avec $[G]$ non 
principal (exemple : $F=I_3$, $G=F_3$, cf. \ref{subs:geom_1_2_3}) ; 
en particulier l'ordre d'un type géométrique n'est pas bien défini.
\par
L'inclusion entre parties de $P_n$ induit un ordre partiel (encore
appelé inclusion) sur l'ensemble des types géométriques. 
D'après ce qui précède, tout type géométrique isodual est inclus
dans un type  géométrique principal. Les {\em types géométriques
maximaux} pour l'inclusion sont donc principaux.

\begin{exem}
Le type symplectique et les types orthogonaux indéfinis
(à l'exception de la signature $(1,1)$) sont toujours
maximaux (proposition~\ref{prop:orth_symp_max},~(2)).
\end{exem}

Pour la recherche de la densité maximale des réseaux isoduaux,
il suffira donc de ne considérer que les  types  géométriques
 principaux et on pourra évidemment se restreindre aux types maximaux.

\subsection{Méthode de classification}
\label{subs:method_class}
Nous décrivons ici, dans les grandes lignes,   les deux étapes 
(algébrique puis géométrique) de la classification.
La dimension $n$ étant fixée, nous déterminerons d'abord les
types isoduaux algébriques en partant d'un système fini $\mcF_n$
de représentants des classes de conjugaison dans $\GL_n(\Z)$
des éléments d'ordre fini de $\SL_n(\Z)$ (tables \ref{tabl:of1-4} et 
\ref{tabl:of5-7} de l'annexe). Tout type isodual est
représenté par un élément $F \in \GL_n(\Z)$  tel que  $F F^\vee \in
\mcF_n$. Pour chaque $R \in \mcF_n$, posons
$\mcT_R=\{F \in \GL_n(\Z);F F^\vee=R\}$, et observons que $\mcT_R$
est aussi l'ensemble des solutions $F$  entières et de
déterminant $\pm 1$ de l'équation {\em linéaire} F=RF'.
Fréquemment les solutions entières de cette équation ont un déterminant
divisible par un entier $\geq 2$, auquel cas  $\mcT_R$ est vide (voir
des  exemples au \S \ref{subs:typalg1-4}).

Il est clair que le commutant $C_R$ de $R$ dans $\GL_n(\Z)$ opère
dans $\mcT_R$. De plus deux matrices de $\mcT_R$ ne peuvent
être équivalentes que par un élément de $C_R$. Par conséquent
l'action de  $C_R$ sur $\mcT_R$ n'a qu'un nombre fini d'orbites
(théorème \ref{theo:finitude_types}) et les types associés sont
tous distincts. Pour $(P,F) \in C_R \times \mcT_R$, les coefficients 
de la matrice $P \cdot F$ sont des formes quadratiques entières 
en les coefficients de $P$. 
La description complète des orbites  de $C_R$ dans $\mcT_R$ repose sur 
la possibilité de réaliser de \og petites valeurs \fg\ de ces formes
par des éléments $P \in C_R$ ; cette étape est détaillée 
au \S \ref{subs:preliminaires}.
Pour chaque dimension étudiée, nous donnerons un ensemble fini de 
représentants ainsi qu'une liste finie 
d'invariants qui caractérisent les types
(voir \S\S \ref{subs:typalg1-4}-\ref{subs:typalg7}). 
\par 

Il existe une action du groupe $(\Z/2\Z)^2$ sur les types
algébriques isoduaux ; en effet si $F$ représente un type
algébrique,
alors $F'$ et  $-F$  aussi (noter que $[F^\vee]=[F]$ et 
$[F^{-1}]=[F']$). Concernant les types  géométriques associés, on a 
\begin{equation}\label{equa:op_types_geom}
V_{F'}=V_{-F}=V_F .
\end{equation}

Pour classer les types géométriques, nous commencerons par
déterminer les sous-variétés $V_F$ associées à un système de
représentants des types algébriques préalablement réduit 
par l'action de $(\Z/2\Z)^2$. Cette étape est facilitée 
par les propriétés de décomposition (\ref{prop:dec1} et 
\ref{prop:dec_can_R}), y compris par la notion de \og type réel \fg\ 
(définition~\ref{defi:type_reel}) car de nombreux types
algébriques sont indécomposables sur $\Z$ mais se scindent sur $\R$. 
Les sous-variétés $V_F$ sont ensuite comparées grâce au critère
d'inclusion (proposition~\ref{prop:inclusion}) : {\em $V_F$ est incluse dans
$V_G$ \ssi $FG^\vee$ fixe $V_F$ point par point}. Nous
pouvons ainsi décrire l'ensemble ordonné des types géométriques
et en déduire  les types maximaux.
\par
Rappelons que, pour les raisons invoquées dans l'introduction, nous
ne classons à partir du rang $5$ que les types principaux.

\section{Classification des types algébriques ($n \leq 7$)}
\label{sect:class_alg}

\subsection{Éléments d'ordre finis de $\GL_n(\Z)$}
\label{subs:of_GLnZ}
Les  éléments  d'ordre fini de $\GL_n(\Z)$ pour $n$ petit
peuvent s'obtenir à partir de la  théorie des représentations
entières des groupes  cycliques. Quand l'ordre est premier,
celles-ci sont décrites dans le chapitre 34 de \cite{Curt-Rien1981a}.
Je remercie Gabriele Nebe qui m'a signalé cette référence.

On s'intéresse ici aux  éléments $u \in \GL_n(\Z)$ d'ordre fini 
$d$ pour $n \leq 4$ et d'ordre $2^k$ pour $5 \leq n \leq 7$.
Si $d=p^\al$ avec $p$ premier et $\al \in \N^*$, alors
$u$ doit annuler le polynôme
cyclotomique $\Phi_{p^\al}$ et l'on a $p^{\al-1}(p-1) \leq n$.
De là (et avec des  arguments analogues), on déduit 
aisément que les valeurs possibles de $d$ sont $d\in \{1,2,3,4,6\}$
pour $n=2,3$ et  $d\in \{1,2,3,4,5,6,8,10,12\}$ pour $n=4$
(voir \cite{Vaidyanathaswamy1928a} pour le cas général).

Soit maintenant $p$ un entier premier, soit $\eps_p$ une racine
primitive $p$-ième de l'unité dans $\C$ et soit $\Pi_p=\Z/p\Z$.
Il y a  bijection entre les classes de conjugaison des éléments
d'ordre $p$ de $\GL_n(\Z)$ et les structures de $\Z[\Pi_p]$-module
non triviales  sur $\Z^n$ à isomorphisme près. 
On sait (cf. \cite[p. 729]{Curt-Rien1981a}) qu'un tel module
$M$  est isomorphe à une  somme d'indécomposables de trois types :
trivial ($\Z$), idéal fractionnaire
du corps cyclotomique $\Q[\eps_p]$ ou extension d'un tel idéal ;
en outre le nombre de composantes de chaque type et
le produit des idéaux dans le groupe des classes caractérisent
$M$ à isomorphisme près.  Supposons de plus que le groupe des classes
est trivial,
ce qui est le cas pour $p \in \{2,3,5\}$ (cf. \cite{Washington1997a}).
Les seuls indécomposables à isomorphisme près sont  alors $\Z$,
$\Z[\eps_p]$ et $\Z[\Pi_p]$ (extension de $\Z[\eps_p]$ {\it via}
le morphisme d'augmentation) et  $M$ est caractérisé par les 
multiplicités de ces trois types. 
Les éléments d'ordre $p \in \{2,3,5\}$ sont explicités dans les
tables \ref{tabl:of1-4} et \ref{tabl:of5-7} de l'annexe ($R_p$ et $V_{p-1}$
 correspondent  respectivement à $\Z[\Pi_p]$ et $\Z[\eps_p]$).
En particulier, il y a $[n^2/4]+n$ classes d'éléments d'ordre 2
dans $\GL_n(\Z)$.

Une autre situation agréable se présente quand  le polynôme 
minimal de $u \in \GL_n(\Z)$ d'ordre $d \geq 2$ est le polynôme
cyclotomique $\Phi_d$. C'est le cas pour
$(n,d) = (1,2),(2,3),(2,4),(2,6),(4,5),(4,8),(4,10)$
et pour $(n,d)=(4,12)$ si $u^6+ I_4 =0$.  Alors $\Z^n$ est un
module sur $Z[\eps_d]=\Z[X]/(\Phi_d)$, sans torsion (appliquer 
Bézout dans $\Q[X]$), donc isomorphe à une somme d'idéaux
fractionnaires de $\Q[\eps_d]$ (\cite[théorème 4.13 p. 85]{Curt-Rien1981a}). 
Pour les valeurs de $(n,d)$ ci-dessus, le groupe des classes
d'idéaux est trivial (\cite{Washington1997a}) et $\Z^n$ est
isomorphe à $Z[\eps_d]$ comme $Z[\eps_d]$-module.

Le cas $n=2$ résulte entièrement de ce qui précède.
Nous explicitons ensuite les $u \in \GL_n(\Z)$ d'ordre $d$ non
premier par des considérations élémentaires en procédant
suivant les valeurs croissantes de $n \geq 3$.
Observons d'abord que 
pour tout $P \in \Z[X]$, le noyau  $L$ de $P(u)$ est un facteur
direct de $\Z^n$ stable par $u$. Quand $L$ est non trivial,
$\Z^n$ apparaît comme une extension de $\Z[\Z/d\Z]$-modules
de rang sur $\Z$ plus petit que  $n$.  Les ordres $d$ considérés ici
sont tous pairs et nous prendrons  généralement $P(X)=X^{d/2}-1$
pour utiliser les  involutions. On peut supposer que 
$u=\smallmat{A}{X}{0}{B}$. Comme $u^{d/2}$ est une involution, le rang
$k$ de $L$ est égal à la multiplicité de la valeur propre 1 dans le
polynôme caractéristique de $u^{d/2}$. Par suite $B^{d/2}+I_{n-k}=0$.
Le cas éventuel où $u^{d/2}+I_n=0$ ($k=0$) se traite séparément ;
par exemple, si $d=4$ on a un $Z[\eps_4]$-module  (voir plus haut),
si $(n,d)=(3,6)$ alors $-u$ est d'ordre 3, \ldots   
Quand $L$ est non trivial, les matrices $A$ et $B$ sont supposées
connues  à conjugaison près puisque $A^{d/2}=I_k$ et $B^{d/2}=-I_{n-k}$
(éléments d'ordre fini et de rang au plus $n-1$). Dans la
suite on  fixe $A$ et $B$.
Soient $C_A$ et $C_B$ leurs commutants respectifs dans $\GL_k(\Z)$
et $\GL_{n-k}(\Z)$ et soit $T \in M_{k,n-k}(\Z)$. On a les relations 
\renewcommand{\arraycolsep}{0.3 em}
\begin{equation}\label{equa:conjT}
\renewcommand{\arraystretch}{1}
\left(
\begin{array}{cc}
I_k &  T\\
0 & I_{n-k}\\
\end{array}
\right)
\left(
\begin{array}{cc}
A &  X\\
0 & B\\
\end{array}
\right)\left(
\begin{array}{cc}
I_k &  T\\
0 & I_{n-k}\\
\end{array}
\right)^{-1}
=
\left(
\begin{array}{cc}
A &  X+TB-AT\\
0 & B\\
\end{array}
\right)
\end{equation}
et pour $(P,Q) \in C_A \times C_B$
\begin{equation}\label{equa:conjPQ}
\renewcommand{\arraystretch}{1}
\left(
\begin{array}{cc}
P &  T\\
0 & Q\\
\end{array}
\right)
\left(
\begin{array}{cc}
A &  X\\
0 & B\\
\end{array}
\right)\left(
\begin{array}{cc}
P &  T\\
0 & Q\\
\end{array}
\right)^{-1}
=
\left(
\begin{array}{cc}
A &  PXQ^{-1}\\
0 & B\\
\end{array}
\right).
\end{equation}
On voit que  l'ensemble $M_{A,B}=\{TB-AT;T \in M_{k,n-k}(\Z)\}$ est un 
sous-module d'indice fini du $\Z$-module $M_{k,n-k}(\Z)$. 
En effet, si $TB-AT=0$, on a $TB^j=A^jT$ pour tout $j \in \N$,
d'où  $T=0$ ($j=d/2$). Par \ref{equa:conjT} on peut donc ramener 
$X$ dans une liste finie de représentants modulo $ M_{A,B}$.
Cette liste peut éventuellement être encore réduite grâce à 
l'action des commutants de $A$ et de $B$, relation \ref{equa:conjPQ}.
Enfin il faut discerner les classes de conjugaison.
On distinguera toujours deux éléments $u_1$ et  $u_2$ 
(ayant le même polynôme caractéristique) en explicitant
les solutions entières $w$ de l'équation linéaire  $w u_1 - u_2 w =0$,
puis en vérifiant que le déterminant d'une telle solution
 est divisible par un entier distinct de $\pm 1$. 
Tous les résultats sont rassemblés dans les tables 
\ref{tabl:of1-4} et \ref{tabl:of5-7}  de l'annexe.

Détaillons un exemple pour illustrer la méthode : $(n,d)=(6,8)$.
Les notations sont définies dans l'appendice. Le polynôme minimal
de la matrice $B$ est ici $\Phi_8$, donc $B$ correspond à
un $\Z[\eps_8]$-module et nécessairement $\ker(u^4-I_6)$ est de rang
$k=2$. On peut donc prendre $B=X_4$ et $A \in \{I_2,-I_2,I_{1,1},
R_2,J_2\}$. Dans tous les cas, $M_{A,B}$ est d'indice 4 dans
$M_{2,4}(\Z)$. Après réduction supplémentaire par \ref{equa:conjPQ},
on se ramène à $u \in \{I_2X_4,I_1X_5\}$ si $A=I_2$,
 $u \in \{I_2^-X_4,I_1^-X_5^-\}$ si $A=-I_2$,
$u \in \{I_{1,1}X_4,I_1X_5^-,I_1^-X_5,Z_6\}$ si $A=I_{1,1}$,
$u \in \{R_2X_4,Z'_6,Z''_6\}$ si $A=R_2$ et
$u \in \{J_2X_4,X_6,Y_6\}$ si $A=J_2$. On vérifie enfin que cette
liste de 14 éléments ne contient pas de paire d'éléments conjugués en 
triant suivant le polynôme caractéristique, puis à l'aide 
d'un calcul modulo~2.

Certains cas se règlent plus directement grâce à des propriétés
spécifiques. Par exemple si $(n,d)=(4,12)$, le polynôme minimal
$\mu$ de $u$ vaut $\Phi_{12}$ (traité plus haut), $\Phi_3 \Phi_4$
ou $\Phi_4 \Phi_6$. Dans ces deux derniers cas on a $\mu=PQ$
où $(P,Q) \in \Z[X]^2$ vérifie une relation de Bézout dans
$\Z[X]$, à savoir $\Phi_3 \Phi_6 - X^2 \Phi_4 =1$. Donc 
$\Z^4=\ker P(u) \oplus \ker Q(u)$ se scinde comme 
$\Z[\Z/12\Z]$-module.

\smallskip
\begin{rema}
Les éléments d'ordre 4 peuvent également se déduire 
des 9 représentations entières indécomposables de $\Z/4\Z$
(cf. \cite{Roiter1960a}, \cite[ch. 34]{Curt-Rien1981a}).
\end{rema}

\subsection{Opérations sur les types. Décompositions}
\label{subs:oper_et_dec_types}
Soient $A = (a_{i,j}) \in M_m(\Z)$ et $B\in M_n(\Z)$  deux
matrices carrées.  On note $A \oplus B = \smallmat{A}{0}{0}{B}$ la
{\em somme directe}  et  $A \otimes B = (a_{i,j}B)$  
(matrice formée de $m^2$  blocs) le {\em produit tensoriel}.
Ces deux opérations sont 
associatives et commutent avec la multiplication matricielle
 et la transposition;
elles sont commutatives à équivalence près. Par conséquent, si
$F$ et $G$ représentent des types isoduaux, il en de même 
pour $F \oplus G$ et  $F \otimes G$. De plus, ces opérations
induisent deux structures de monoïdes  commutatifs sur l'ensemble
$\mcT(\Z)$ des types isoduaux.

\smallskip

\begin{rema}
On pourrait enrichir ces  structures comme dans la théorie des
formes bilinéaires symétriques, mais nous n'utiliserons pas ce
fait ici. En effet le groupe de Grothendieck de $(\mcT(\Z),\oplus)$
est muni d'une structure d'anneau unitaire commutatif grâce au
produit tensoriel.
Notons $\hat{T}(\Z)$ cet \og anneau des types isoduaux \fg. 
Noter que $\mcT(\Z)$ ne s'injecte pas dans $\hat{T}(\Z)$ 
(la forme paire $U_2=\smallmat{0}{1}{1}{0}$ est stablement
équivalente à la forme impaire $\smallmat{1}{0}{0}{-1}$). Par ailleurs,
on vérifie aisément  que  l'anneau de Grothendieck-Witt de $\Z$
s'injecte dans $\hat{T}(\Z)$ et de même que l'anneau de Witt de $\Z$
s'injecte dans $\hat{T}(\Z)/U_2 \hat{T}(\Z)$.
\end{rema}

Si $[F]$ est un type  isodual, la matrice $F$ définit une forme
 bilinéaire sur $\Z^{n*}$ et on a la notion d'orthogonalité
associée, en général non symétrique. 
Le résultat suivant est très utile pour les questions de 
décomposition.

\begin{lemm}%
\label{lemm:orthog_bilat}
{\rm (lemme d'orthogonalité bilatère)}
Soit $F$ un élément de $\GL_n(\Z)$ et 
soit $R=FF^\vee$. Soit $(P,Q) \in \Z[X]^2$.
On pose $P^*(X)= X^{\deg (P)} P(1/X)$ (polynôme ré\-ci\-pro\-que
de $P$) et on suppose que $P^*$ et $Q$ sont premiers entre
eux dans $\Q[X]$. Alors les sous-modules $\ker P(R^\vee)$ et  
$\ker Q(R^\vee)$ sont bilatéralement $F$-or\-tho\-go\-naux.
\end{lemm} 

\pr Noter   d'abord que  $FR^{\vee m}=R^mF$ pour tout $m \in \Z$.
Soit $x \in \ker P(R^\vee)$ et soit $y \in \ker Q(R^\vee)$.
Par hypothèse, il existe $(A,B) \in \Q[X]^2$ tel que
$AP^*+BQ=1$. On a donc $y=A(R^\vee)P^*(R^\vee)y$. Par suite,
on a 
$x'Fy=x' P^*(R)A(R) F y$. Mais $P^*(R')=
R'\null^{\deg (P)} P(R^\vee)$, d'où $x' P^*(R)=0$ et 
\mbox{$x'Fy=0$.}
Comme $Q^*$ et $P$ sont premiers entre eux, on a
aussi $y'Fx=0$ par le même raisonnement.

\begin{prop} 
\label{prop:dec1}
Soit $R = R_1 \oplus R_2 \in \GL_n(\Z)$ un élément d'ordre fini et 
soit $F\in \GL_n(\Z)$ tel que $FF^\vee=R$.
On suppose que les polynômes carac\-té\-ris\-tiques $\chi_{R_1}$
et $\chi_{R_2^\vee}$ ($=\chi_{R_2}$\footnote{$R_2$ est conjugué
sur $\R$ à un élément de $\OO(n)$}) sont premiers entre eux dans
$\Q[X]$. Alors\\
\indent {\rm a)} $F$ se scinde en $F=F_1 \oplus F_2$ avec 
$F_jF_j^\vee=R_j$ ($j=1,2$),\\
\indent {\rm b)} on a $V_F=\{A_1 \oplus A_2;(A_1,A_2) \in V_{F_1} \times
V_{F_2}\}$,\\
\indent {\rm c)} le groupe d'isométries de $F$ est isomorphe au
produit des groupes d'iso\-mé\-tries respectifs de $F_1$ et $F_2$.
\end{prop}

\pr Soient $\ro_1$ et $\ro_2$ les représentations complexes définies
par $R_1$ et $R_2$. L'hypothèse sur les polynômes caractéristiques
revient à dire  que   $\ro_1$ et $\ro_2^\vee$  n'ont pas de type
irréductible en commun (ce qui équivaut aussi à l'orthogonalité
de leurs caractères). Par conséquent la seule matrice complexe $M$
vérifiant $R_1 M= M R_2^\vee$ est $M=0$. De plus, comme  $R_2$ est
réelle, on a $\chi_{R_2^\vee}=\chi_{R_2}$ et la même
propriété est vraie si on remplace $R_2^\vee$ par $R_2$.
En décomposant $F$ en 4 blocs et en explicitant $F=RF'$,
on voit d'après ce qui précède les blocs non diagonaux sont nuls.
La matrice $F$ est alors diagonale par blocs, ce qui prouve~a). 
Cette assertion résulte aussi du lemme d'orthogonalité bilatère
\ref{lemm:orthog_bilat}. Soit maintenant $A \in V_F$. On a
$RAR'=A$ et l'on voit de même que $A$ est diagonale par blocs,
d'où l'assertion~b). Enfin toute matrice $P \in \GL_n(\Z)$ telle que
$P'FP=F$ (isométrie de $F$) doit commuter avec $R'$;
on trouve à nouveau que $P$ est diagonale par blocs, ce qui établit 
l'assertion~c).
 
\begin{prop} 
\label{prop:dec_can}
{\rm (décomposition canonique d'un type isodual)}\\
Soit \mbox{$F \in \GL_n(\Z)$} tel que $R= FF^\vee$ soit
 d'ordre fini~$d$.
Pour tout diviseur $k$ de~$d$ on pose $M_k=\ker \Phi_k(R^\vee)$
où $\Phi_k$ est le $k$-ième polynôme cyclotomique. Alors\\
\indent {\rm a)} la somme des $M_k$ est directe et  d'indice fini dans
$\Z^n$,\\
\indent {\rm b)} cette somme directe $\bigoplus_{k|d}M_k$ est
$F$-orthogonale {\em bilatère},\\
\indent {\rm c)} l'indice $[\Z^n:\bigoplus_{k|d}M_k]$
et les classes d'équivalence des $\Z$-modules bili\-né\-ai\-res
$(M_k,F)$ sont des invariants du type $[F]$.
\end{prop}

\pr On étend l'action de $R^\vee$ à $\Q^n$. 
Le polynôme minimal de $R^\vee$ est le produit de certains
$\Phi_k$ ($k|d$), tous  irréductibles sur $\Q$. L'espace vectoriel
$\Q^n$ est donc la somme directe des noyaux correspondants, d'où a).
L'assertion b)  résulte du lemme \ref{lemm:orthog_bilat} puisque
$\Phi_1 ^* =-\Phi_1$ et $\Phi_k^*=\Phi_k$ pour
$k \geq 2$. Soit maintenant $P \in \GL_n(\Z)$; on pose
$G=PFP'$ et $S=GG^\vee$. On a $S^\vee=P^\vee R^\vee P'$. Il est
clair que $P^\vee : (\Z^n,F) \to (\Z^n,G)$ est une isométrie
qui envoie  $\ker \Phi_k(R^\vee)$ sur $\ker \Phi_k(S^\vee)$ pour
tout diviseur $k$ de $d$, ce qui prouve l'assertion c).
\medskip

Le résultat suivant permet d'obtenir des types scindés à partir
de certains éléments $R \in \GL_n(\Z)$ décomposés ne vérifiant
pas l'hypothèse de la proposition~\ref{prop:dec1} (voir par exemple
au \ref{subs:typalg1-4} les cas  $I_1\oplus R_3$ ou $W_2\oplus W_2$).

\begin{lemm} 
\label{lemm:scind}
{\rm (lemme de scindement)}
Soit $R\oplus S \in \GL_{k+l}(\Z)$ et soit 
$F = \smallmat{A}{B}{C}{D}\in \GL_{k+l}(\Z)$ une solution  de
$FF^\vee=R\oplus S $. On suppose que  $A \in \GL_k(\Z)$. Alors $F$
est équivalente à une matrice décomposée $G=A \oplus M$
($M \in \GL_l(\Z)$) telle que $GG^\vee=R\oplus S$.
\end{lemm}

\pr L'hypothèse sur $F$ implique (entre autres) les relations
$A=RA'$, $B=RC'$ et $SC=C R^\vee $. Par suite, si on pose
$U=-CA^{-1}$, on a $SU = UR$ et la matrice triangulaire 
$T=\smallmat{I_k}{0}{U}{I_l}$ commute avec $R\oplus S$. Alors $G=TFT'$
convient.

\subsection{Sur l'existence de type associé à un élément de $\GL_n(\Z)$}
\label{susb:existence}
Considérons un élément $R \in \GL_n(\Z)$ d'ordre fini $d$
et soit $\Pi=\Z/d\Z$. Le  $\Z$-module $\Z^n$ sera muni
de sa base naturelle, son dual de la base duale et on
fera ainsi correspondre matrices et applications linéaires. 
Soit  $V$ le $\Pi$-module défini par $R$ et soit $V^\vee$ son
module contragrédient, défini par  $R^\vee$. L'équation
\begin{equation}
\label{equa:F=RFt}
F=RF'
\end{equation} 
entraîne $FR^\vee = RF$. En plus de la relation évidente $\det(R)=1$,
on a donc une deuxième condition nécessaire
à l'existence de type $F \in \GL_n(\Z)$ vérifiant (\ref{equa:F=RFt})~:
{\em il faut que $V$ et $V^\vee$ soient isomorphes comme $\Pi$-modules}.
Cette condition n'est pas toujours remplie (voir exemple
\ref{exem:contragrediente}) et elle
n'est pas suffisante (exemple \ref{exem:epsxi}).
Supposons qu'il existe un $\Pi$-isomorphisme $\psi$ de $V$ dans $V^\vee$ ;
sa matrice $\Psi$ vérifie  $\Psi R = R^\vee  \Psi$.  
On définit un anti-isomorphisme $\sig$ de l'anneau $\endo_\Pi V$ 
en posant $x^\sig=\psi^{-1}x^t \psi$  pour $x \in \endo_\Pi V$.
Soit $F$ une solution de (\ref{equa:F=RFt}) et soit 
$f \in \hom_\Pi(V^\vee,V)$ l'application  associée.
En notant  $G$ la matrice de $f\psi$, l'équation 
(\ref{equa:F=RFt}) se réécrit dans $\endo_\Pi V$  sous la forme
matricielle
\begin{equation}
\label{equa:G=RxiGsig}
G=R (\Psi^\vee \Psi) G^\sig .
\end{equation} 
Remarquer que $\psi^\vee \psi(\psi^\vee \psi)^\sig = id_V$ ; de plus
si $\psi^\vee \psi$ est  central dans $\endo_\Pi V$, alors $\sig$
est involutif. Ce sera presque toujours vérifié  pour les exemples
traités ici. Enfin, l'action par équivalence de $\ph \in \GL_\Pi(V)$
sur  $f$ se traduit dans $\endo_\Pi V$ par $(\ph f \ph^t) \psi =
\ph (f \psi) \ph^\sig$.

Examinons maintenant le cas -- qui se présente fréquemment  dans la
classification des types -- où le polynôme minimal de  $R$ sur $\Q$
est le polynôme cyclotomique $\Phi_d$. Soit $\eps$ une racine
primitive $d$-ième de l'unité et soit $\mcO=\Z[\eps]$.
Les $\Pi$-modules $V$ et $V^\vee$ sont aussi des  $\mcO$-modules de
type fini et sans torsion, donc décrits par des idéaux fractionnaires
du corps  $\Q[\eps]$ (\cite[théorème 4.13 p. 85]{Curt-Rien1981a}).
Soit $\mfa$ un tel idéal et soit $S$ la matrice de la multiplication
par $\eps$ dans une certaine base de $\mfa$ sur $\Z$. Il est
immédiat de constater que le $\Pi$-module défini par $S^{-1}$
est isomorphe à l'idéal conjugué $\overline{\mfa}$. Par ailleurs,
d'après \cite[p.~59]{Lang1994a}, $\hom_\Z(\mfa,\Z)$ s'identifie sur $\Z$ à
$\Phi_d'(\eps)^{-1}\mfa^{-1}$ {\it via} la forme trace. On en déduit que
le contragrédient $\mfa^\vee$ de $\mfa$ comme  $\Pi$-module est
isomorphe à l'idéal $\overline{\mfa}^{-1}$. En particulier, $\mcO$ et
$\mcO^\vee=\hom_\Z(\mcO,\Z)$ sont toujours $\Pi$-isomorphes.
Cependant, comme le montre l'exemple suivant, cette propriété n'est
pas vraie pour tout idéal de  $\Q[\eps]$.

\begin{exem}
\label{exem:contragrediente}
Soit  $C$ le groupe des classes de $\Q[\eps] $ et soit  
$C^+$ celui de  $\Q[\eps+\overline{\eps}]$.
On sait que $C^+$ s'injecte dans $C$ et que, pour  $d$ premier,
l'application norme de $C$ dans $C^+$ est surjective 
(\cite{Washington1997a}, théorèmes 4.14 et 10.1 respectivement).
Pour $d=257$, le groupe $C^+$ est non trivial
(\cite{Anke-Chow-Hass1965a}) ; il  existe donc un idéal 
$\mfa$ distinct de  $\overline{\mfa}^{-1}$ dans $C$, 
\cadab  $\mfa$ et $\mfa^\vee$ non isomorphes comme $\Pi$-modules. 
\end{exem}

À partir d'ici, nous supposerons que le nombre $h_d$ de classes du
corps $\Q[\eps]$ vaut 1, ce qui est le cas pour les entiers $d$ qui nous
intéressent. Dans ces conditions, on peut prendre $V=\mcO^m$ 
($m=n/\deg(\Phi_d)$) et  supposer que $R$
est l'homothétie de rapport~$\eps$. Considérons d'abord le cas $m=1$.
Soit $\psi_0$ un  $\mcO$-isomorphisme de $\mcO$ sur $\mcO^\vee$,
de matrice $\Psi_0$. La relation $R^\sig=\Psi_0^{-1}R'\Psi_0=R^{-1}$ 
montre que $\sig$ induit la conjugaison complexe de $\mcO$, 
indépendamment du choix de $\psi_0$. Par ailleurs,
$\psi_0^\vee \psi_0 \in \endo_\Pi(\mcO) \simeq \mcO$
s'identifie à une unité $\xi \in \mcO^\times$. De plus on a
$\xi \overline{\xi}=1$, donc $\xi$ est une racine de l'unité
puisque $\Q[\eps]$ et un corps CM.
Pour $m$ quelconque, on 
prend $\psi=\psi_0 \oplus \cdots \oplus \psi_0$ ($m$ facteurs).
Les éléments de $\endo_\Pi V$ se représentent par des matrices
$m \times m$ à coefficients dans $\mcO$ ($\mcO^m$ étant muni
de sa base canonique sur $\mcO$) ou par des matrices à coefficients 
dans $\Z$ formées de $m^2$ blocs. Ainsi on voit (compte tenu du cas
$m=1$) que  pour tout $X \in M_m(\mcO)$, $X^\sig$ est la transposée
de la conjuguée complexe de $X$, que nous noterons $X^*$.
Soit $\mfF \in M_m(\mcO)$ la matrice de $f\psi$. Dans 
$M_m(\mcO)$, l'équation (\ref{equa:G=RxiGsig}) prend la forme
\begin{equation}
\label{equa:frakF}
\mfF= \eps \xi \mfF^*.
\end{equation} 
En particulier, on a la relation $\det(\mfF)=(\eps \xi)^m
\overline{\det(\mfF)}$.
Si l'on change le choix  de $\psi_0$, l'unité $\xi$ est multipliée
par une unité de la forme $u/\overline{u}$ avec $u \in \mcO^\times$.
Nous avons donc prouvé le résultat suivant.
\begin{prop} 
\label{prop:cyclotomic}
Soit $F \in \GL_n(\Z)$ tel que $\Phi_d(FF^\vee)=0$.
On suppose que $h_d=1$. Soit $\eps$ une racine primitive $d$-ième de
l'unité et  soit $\psi_0$ un isomorphisme de
$\mcO=\Z[\eps]$ sur $\mcO^\vee$. On pose $m=n/\deg(\Phi_d)$ et
$\xi=\psi_0^\vee \psi_0 (1)$. Alors $(\eps \xi)^m$ est de la forme 
$u/\overline{u}$ avec $u \in \mcO^\times$, condition indépendante
des choix de $\eps$ et $\psi_0$.
\end{prop}

\begin{exem}
\label{exem:epsxi}
Cette condition est  toujours satisfaite  pour $m$ pair, comme 
on le voit en prenant $u=(\eps \xi)^{m/2}$.
Si $d=2^k$ ($k\geq 1$), on a $\xi=1$ et on sait qu'il n'existe pas 
d'unité $u \in \mcO^\times$ telle que $\eps  =u/\overline{u}$
(preuve du cor. 4.13, p. 40 de \cite{Washington1997a});
l'existence de type isodual $F$
tel que $\Phi_{2^k}(FF^\vee)=0$ n'est donc possible que si $m$ est pair,
\cadab $\rang F$ divisible par $2^k$. Un phénomène analogue se 
produit pour $d=3$ ($\xi=-\eps^2$) et pour  $d=5$ ($\xi=-\eps^2$). 
\end{exem}

\subsection{Commentaires et conventions préliminaires}
\label{subs:preliminaires}
Un des outils fondamentaux dans la détermination des types
est la classification des formes
quadratiques entières de petit déterminant (cf. par exemple
\cite[ch. 15]{Conw-Sloa1999a} ou \cite[ch. V]{Serre1970a} pour le
déterminant $\pm 1$) qui renseigne sur la composante symétrique
$(M_1,F)$ des types (proposition~\ref{prop:dec_can}). En particulier les types
d'ordre~1, \cadab symétriques,  de rang $\leq 7$ sont représentés
par $I_{p,q}$ (type impair) et par les sommes de $U_2$ (type pair).
Il est aussi bien connu qu'il n'existe  qu'un seul type
antisymétrique ($FF^\vee=-I_n$) représenté par
$J_{2m}=J_2 \otimes I_m$; plus généralement, d'après 
\cite[\S 5 n$^{\rm o}$ 1]{BourbakiAlg9-1959a}, la composante
antisymétrique $(M_2,F)$ d'un type d'ordre pair est isométrique à
$(\Z^{2m},a_1J_2 \oplus \cdots \oplus  a_mJ_2)$ pour une certaine
suite décroissante $(a_1), \cdots  ,(a_m)$ d'idéaux de $\Z$
($m \in \N$).

Nous procéderons  à la classification des types algébriques suivant
la méthode générale indiquée au \S \ref{subs:method_class}.
Les éléments d'ordre fini de $GL_n(\Z)$ sont donnés dans l'annexe
dont on reprend les notations. Comme les types symétriques
et antisymétriques sont connus, {\em nous pourrons toujours supposer
que $R=FF^\vee$ est distinct de $\pm I_n$}. 
Pour chaque matrice $R \in \SL_n(\Z)$ des  tables  \ref{tabl:of1-4}
et  \ref{tabl:of5-7} de l'annexe, il s'agira de classer les 
solutions $F \in  \GL_n(\Z)$ de l'équation linéaire $F=RF'$
modulo l'action par équivalence du commutant $C_R$.
Nous chercherons fréquemment des informations sur
la composante symétrique de $F$ -- quand elle non triviale.
Pour cela, nous noterons $\De$ la transposée de la matrice d'une
base du  sous-module $\bigoplus_{k|d}  \ker \Phi_k(R^\vee)$.
Si  $P\in C_R$, alors $Q= \De P \De^{-1} = \oplus_{k|d} Q_k$
est diagonale par blocs et chaque $Q_k$ s'identifie à un 
élément $\mfP_k$ de $\GL(W_k)$, où $W_k$ est un module sur
l'anneau $\mcO_k$ des entiers du corps des racines $k$-ièmes
de l'unité -- ici on aura toujours  $W_k=\mcO_k^{m_k}$ 
($m_k \in \N^*$). Comme la \og diagonalisation \fg\ n'est
que rationnelle il y a généralement des relations arithmétiques
entre les éléments $\mfP_k$: le commutant $C_R$ est isomorphe à
un sous-groupe, souvent propre, de  $\prod_{k|d}\GL(W_k)$.
Par ailleurs, on a évidemment
$Q \cdot (\De \cdot F)= \De \cdot (P \cdot F)$ et  l'action
de $C_R$ sur $F$ sera précisée {\it via} celle
du sous-groupe diagonal par blocs $C_R^\De=\De C_R\De^{-1}$. 
Nous utiliserons notamment les  matrices suivantes :
\begin{equation}
\label{equa:De3_De4}
{\small 
\De_3=\left(
\begin{array}{ccc}
2 & -1 & 1 \\
0 &  1 & 0 \\
0 & 0 &  1 
\end{array}
\right)
\esp {\rm et} \esp
\De_4=\left(
\begin{array}{cccc}
1 & 1 & 0 & 0 \\
0 & 0 & 1 & 1 \\
-1 & 1 & 0 & 0 \\
0 & 0 & 1 & -1 
\end{array}
\right),}
\end{equation}
vérifiant $\De_3 X_3 \De_3^{-1} = I_1 \oplus J_2$, 
$\De_4(R_2 \oplus R_2) \De_4^{-1} = I_{2,2}$ et $\det(\De_4)=4$.

Dans ce qui suit, nous noterons toujours $F$ une  solution
de $F=RF'$ et~$P$ une  matrice qui commute avec $R$.
Les coefficients de  $PFP'$ sont à $F$ fixée des 
formes quadratiques en $P$, certaines étant parfois
liées à la composante symétrique de~$F$. Comme  première
étape de la classification ($R$ étant donné), il pourra être
utile dans certains cas de ramener le premier coefficient
de $PFP'$, \cadab
\begin{equation}
\label{equa:FPFt_11}
q_F= (PFP')_{1,1},
\end{equation}
à des valeurs particulières en utilisant à nouveau la 
classification des formes entières combinée avec la description
du commutant $C_R$.

Les types de rang $n \leq 4$  et les types principaux de rang
5, 6 et 7  sont classés dans les tables
\ref{tabl:types_1_3} à \ref{tabl:types_7_b}
suivant leur ordre $d$;
les in\-dé\-com\-po\-sa\-bles sont signalés par un astérisque (*)
et explicités dans l'annexe, équations~(\ref{equa:indec_rg2})
à (\ref{equa:indec_rg7}).
Pour la notion de type réel, voir \S\S 
\ref{subs:types_reels_aspects_diff}-\ref{subs:composantes} et  
particulièrement la proposition~\ref{prop:types_alg_reels} 
page~\pageref{prop:types_alg_reels}. Le type réel $P_{2;k,l}$ est 
défini par l'équation (\ref{equa:type_Pkl}).
Dans toutes ces tables, les types avec $FF^\vee$ fixé 
(à conjugaison près) sont distingués par leur type réel et
à égalité des types réels par leur composante symétrique 
(voir proposition~\ref{prop:dec_can}), par exemple par sa parité.

\subsection{Types algébriques pour $n \leq4$}
\label{subs:typalg1-4}
Pour chaque matrice $R \in \SL_n(\Z)$ de la table  \ref{tabl:of1-4} 
(p.~\pageref{tabl:of1-4})  il est facile d'expliciter les
solutions entières $F \in  M_n(\Z)$ de l'équation linéaire
(\ref{equa:F=RFt}).
Dans certains cas on constate que $\det(F)$ est divisible par un
entier $\geq 2$, donc (\ref{equa:F=RFt}) n'a pas de solution dans
$\GL_n(\Z)$ (voir aussi \S \ref{susb:existence},
ex. \ref{exem:epsxi}).  Par exemple : si $R=J_2$,
alors $F=\smallmat{a}{a}{-a}{a}$ et 
$\det(F)=2a^2$ ; si $R=V_2$ alors $F=\smallmat{a}{-a}{2a}{a}$ et 
$\det(F)=3a^2$ ($a \in \Z$). On peut ainsi éliminer $V_2$, $J_2$ 
($n=2$),  $I_1^-R_2$, $I_1V_2$, $I_1J_2$ ($n=3$) et 
$I_{1,1}R_2$, $I_2V_2$, $I_2^-V_2$, $I_2J_2$, $I_2^-J_2$, $I_1^-X_3^-$, 
$V_4$, $I_1^-R_3^-$, $V_2W_2$, $T_4$, $X_4$, $V_2J_2$, $W_2J_2$ ($n=4$)
-- ici $AB=A \oplus B$.

Remarquons ensuite  que les types associés à $I_{1,2}$ et $I_1W_2$
sont tous scindés (proposition~\ref{prop:dec1}). Les seuls cas qui
restent à examiner pour $n \leq 3$ sont $R=W_2$, $R_3$ et $X_3$.
Une vérification immédiate montre que les seules solutions
$F \in \GL_n(\Z)$ de (\ref{equa:F=RFt}) sont alors $\pm F_2$
si $R=W_2$, $\pm F_3$ si $R=R_3$ et  $\pm G_3$, $\pm H_3$
si $R = X_3$ (équations (\ref{equa:indec_rg2}) et
(\ref{equa:indec_rg3}) de l'annexe). Dans la table
\ref{tabl:types_1_3} (types de rang $n \leq 3$) 
l'entier $t$ (resp. $t'$) désigne le nombre de types (resp. 
de types principaux)  de rang donné.

\begin{table}[h]
\centering
\renewcommand{\arraystretch}{1.2}
\renewcommand{\tabcolsep}{0.3 em}
{\small 
\begin{tabular}{llllll}
\hline
$n$ &$(t;t')$ & $d$  & $F F^\vee$ & $F$ &  type réel\\ \hline
1 &(2;2) &1   & $I_1$ 	& $\pm I_1^*$ &  $\pm I_1$\\
\hline
2  &(7;5)&1   & $I_2$ & $I_{p,2-p}$ &   $I_{p,2-p}$ \\
 & & &  & $U_2^*$ & $I_{1,1}$ \\
  & &2   & $I_2^-$  & $J_2^*$ & $J_2$\\
 & &6 & $W_2$  & $\pm F_2^*$ &   $\pm P_{2;6,1}$ \\
\hline
3  & (16;10)&1 & $I_3$ & $I_{p,3-p}$ & $I_{p,3-p}$ \\
   & &2 & $I_{1,2}$ & $\pm I_1 J_2$ & $\pm I_1 J_2$\\
  & &3 & $R_3$  & $\pm F_3^*$ & $\pm I_1P^-_{2;3,1}$\\
  & &4 & $X_3$ & $\pm G_3^*$& $\pm I_1P^-_{2;4,1}$\\
  & & & & $\pm H_3^*$ &  $\pm I_1P_{2;4,1}$\\
  & &6 &$I_1 W_2$  & $\pm I_{p,1-p} F_2$& $\pm I_{p,1-p} P_{2;6,1}$\\
\hline
\end{tabular}}
\caption{Types isoduaux de  rang $\leq 3$}
\label{tabl:types_1_3}
\end{table}

Pour $n=4$, certains éléments d'ordre fini $R$ ne donnent pas de types
(voir plus haut). Les types associés à $I_{2,2}$, $I_2W_2$ et 
 $I_2^-W_2$ sont tous décomposés (proposition~\ref{prop:dec1}),
voir table \ref{tabl:types_4}. Bien sûr $-I_4$ donne le type
antisymétrique $J_4$.  Il reste à examiner $R_2R_2$, $I_1R_3$
$V_2V_2$, $J_2J_2$, $I_1X_3$, $W_2W_2$, $W_4$ et $Y_4$.  Nous
prouverons que $I_1X_3$ donne aussi des types
scindés grâce au lemme \ref{lemm:scind}; les relations 
éventuelles entre les sommes (telle que (\ref{equa:I1^-G3})
p.~\pageref{equa:I1^-G3}) seront signalées au fur et à mesure.

\paragraph{Cas où $R=R_2 \oplus R_2$ (ordre $2$).}
Les solutions de
(\ref{equa:F=RFt}) sont de la forme
$$
{\small 
F=\left(
\begin{array}{cc}
a(I_2+R_2) &  bI_2+cR_2\\
cI_2+bR_2 &  d(I_2+R_2)\\
\end{array}
\right)}$$
avec $\det (F) =(b-c)^2[4(ad-bc)-(b-c)^2]$.
Si $F \in \GL_4(\Z)$, on a nécessairement $b-c=\pm 1$ (ainsi
$b$ et $c$ sont premiers entre eux), $ad=bc$ et  $\det(F)=-1$.
Il existe donc des entiers $u,v,w$ et $x$ tels que $a=uv$, $d=wx$,
$b=uw$ et $c=vx$ (en particulier $uw-vx = \pm 1$). La matrice
$P=\smallmat{w}{-v}{-x}{u}\otimes I_2$ commute avec $R=I_2 \otimes R_2$
et on vérifie que  $PFP'=\smallmat{0}{I_2}{U_2}{0}$,
qui définit  donc l'unique type associé à $R$, noté $F_4$.

\paragraph{Cas où $R=I_1 \oplus R_3$ (ordre $3$).} 
On a 
$${\small 
P=\left(
\begin{array}{cc}
u &  v \al \\
w \al' &  x I_3+ y R_3 +z R_3^2\\
\end{array}
\right)
\esp
{\rm et}
\esp
F=\left(
\begin{array}{cc}
a &  b \al\\
b \al'  &  c(I_3+R_3)+d R_3^2\\
\end{array}
\right),}$$
avec $\al=(1,1,1)$ -- noter que $\al R_3=\al$. Tout d'abord,
étant donné $(u,v)\in \Z^2$
vérifiant $\pgcd(u,3v)=1$, il existe $(w,x,y,z)\in \Z^4$ tel
que $P \in \GL_4(\Z)$. En effet on a $\det (P) =
(x^2+y^2+z^2 -xy-yz-zx)(u(x+y+z)-3vw)$; il existe $(t,w) \in  \Z^2$
tel que $tu-3vw=1$, où nécessairement $t=3k+\eps$, $\eps =\pm 1$ et
$k \in \Z$. Il suffit de poser $x=y=k$ et $z=k+\eps$ pour 
avoir $\det (P) =1$. Considérons maintenant la forme 
$q_F(u,v)= a u^2 +6buv+3(2c+d)v^2$ définie par l'équation
(\ref{equa:FPFt_11}),  de forme polaire entière
$A=\smallmat{a}{3b}{~3b}{~6c+3d}$ (qui correspond ici à la composante
symétrique de $F$, voir proposition~\ref{prop:dec_can}). Comme
$\det(F)=(c-d)^2 (a(2c+d)-3b^2)= \pm 1$, on voit que
$\det(A) = \pm 3$. Par conséquent, d'après la classification
des formes entières
(cf. \cite[pp. 360 et 362]{Conw-Sloa1999a}, la forme $q_F$ représente
$\pm 1$ si elle impaire et $\pm 2$ si elle paire. Il existe donc
$(u,v) \in \Z^2$ tel que $q_F(u,v)=\pm 1$ ou $\pm 2$, ce qui entraîne
$\pgcd(u,3v)=1$. En faisant agir le commutant $C_R$, on peut
donc supposer que $a= \pm 1$ ou $a= \pm 2$. Dans le premier cas,
$F$ est équivalente à un type décomposé (lemme de scindement
\ref{lemm:scind}) et on est ramené à des types de rang inférieur
(voir table \ref{tabl:types_4}). Quitte à changer $F$ en $-F$, il reste
à examiner le cas où $a=2$ (donc $b$ impair) avec $q_F$ paire
(donc $d$ pair). L'action de 
$\smallmat{1}{0}{n \al'}{I_3} \in C_R$ par équivalence sur $F$
conserve $a$ et change $b$ en $b+na$ ($n \in \Z$) : on peut 
de plus prendre $b=1$. Le couple $(c,d) \in \Z \times 2\Z$
est alors uniquement déterminé par $|\det(F)|=1$, il s'agit de
$c=1$ et $d=0$.   En résumé, nous avons trouvé 6 types, dont
2 indécomposables $\pm G_4$ (équation~(\ref{equa:indec_rg4}) 
p.~\pageref{equa:indec_rg4}) et 4 scindés, discernables
 par leurs composantes réelles (voir table~\ref{tabl:types_4}).

\paragraph{Cas où $R=I_1 \oplus X_3$ (ordre 4).}
\label{para:R=I_1X_3}
Considérons plus 
généralement $R=I_q \oplus X_3$ en dimension $n=q+3 \geq 4$.
Soit $\De=I_q \oplus \De_3$ (voir équation~(\ref{equa:De3_De4})),
de sorte que $\De R\De^{-1}=I_{q+1} \oplus J_2$ et $\det(\De)=2$.
Pour $F$ solution  inversible de (\ref{equa:F=RFt}), on pose 
$\De \cdot F = A_F \oplus B_F$. Les matrices $F$ et $A_F$ sont de
la forme
$$
{\small
F=\left(
\begin{array}{ccc}
A  &  v  & 0\\
v' &  f  & g \al \\
0  & g J_2 \al' & g (I_2 + J_2)\\
\end{array}
\right)
\esp {\rm et} \esp 
A_F=\left(
\begin{array}{cc}
A  &  2v  \\
2v' &  4f-2g  \\
\end{array}
\right),}
$$
avec $A \in M_q(\Z)$ symétrique, $v\in \Z^q$, $\al=(1,0)$ et $f$,
$g \in \Z$ ; en outre $B_F=g (I_2 + J_2)$.
On voit que $g^2$ divise $\det(F)=\pm 1$, d'où $g^2=1$, $\det(B_F)=\pm2$
et $\det (A_F)= \pm 2$. La forme $A_F$ appartient  donc à l'un des genres
${\rm I}_{r,s}(2)$ ou ${\rm II}_{r,s}(2)$
(cf. \cite[pp. 386-387]{Conw-Sloa1999a}). On sait de plus que les
genres indéfinis de rang $\geq 3$ et de petit déterminant ne
contiennent qu'une seule classe d'équivalence entière. Dans 
le cas indéfini, les genres  ci-dessus sont représentés
de façon unique par une forme du type
$B \oplus (\pm 2 I_1)$ avec $\det(B)=\pm 1$ ; c'est encore  vrai pour
${\rm I}_{1,1}(2)$ car toute forme binaire indéfinie représente $1$
(\cite[p. 362]{Conw-Sloa1999a}). Dans le cas défini, notons 
$\ga_k$ la constante d'Hermite en dimension $k$. Pour $2\leq k \leq 6$,
on~a \mbox{$2^{1/k} \ga_k < 2$}, donc toute forme définie de rang
$k$ et de déterminant $\pm 2$ est équivalente à la forme 
diagonale $\pm (I_{k-1}\oplus 2 I_1)$. Finalement, si 
$\rang (A_F) =q+1 \leq 6$, il existe dans tous les cas une
matrice $Q \in \GL_{q+1}(\Z)$ telle que $Q \cdot A_F =
B \oplus (\pm 2 I_1)$ avec  $\det(B)=\pm 1$.
En mettant~2 en facteur dans la dernière
ligne de $A_F$ puis en réduisant modulo~2, on voit que $A$ est 
inversible modulo~2. Il suit que nécessairement  $Q \in \Ga_0^{q,1}$
(voir l'énoncé du lemme \ref{lemm:Ga_0^kl},
p.~\pageref{lemm:Ga_0^kl}),
en particulier $Q[q,q]$ est impair. On vérifie ensuite que
$\De^{-1}(Q \oplus I_2)\De $ est entier, donc appartient au
commutant $C_R$. Par conséquent on peut supposer
que $A_F=B \oplus (\pm 2I_1)$ et tous les types associés à $R$
sont scindés (pourvu que $1 \leq q \leq 5)$. Pour $q=1$,
la liste des 8 types possibles à priori (voir table~\ref{tabl:types_1_3})
se réduit à 6 types (distincts) grâce à la relation 
\begin{equation}
\label{equa:I1^-G3}
[-I_1 \oplus G_3]=
[I_1 \oplus -H_3],
\end{equation}
où l'équivalence est  donnée  par $P=\smallmat{1}{\be}{\ga'}{I_3}$
avec $\be=(2,-1,1)$ et $\ga=(1,0,0)$.

\medskip
Il ne reste à examiner que des cas où le polynôme minimal de
$R$ est un polynôme cyclotomique $\Phi_d$ ($d=3,4,5,6,10,12$).
Posons $m=4/\deg(\Phi_d)$.
Nous procéderons comme indiqué en \ref{susb:existence} pour 
nous ramener au cadre des modules sur l'anneau 
$\mcO=\Z[\eps]$ des entiers cyclotomiques, $\eps$ étant
une racine primitive $d$-ième de l'unité. Le $\Z$-module $\mcO$
sera muni de la base $(1,\eps,\ldots)$ et $\mcO^m$ de sa base
canonique comme $\mcO$-module. Ainsi $\Z^4$ s'identifie à
$\mcO^m$ puisque $h_d=1$ et~$R$ est la
matrice (sur $\Z$) de la multiplication par $\eps$. Par un choix
convenable de $\Psi \in \GL_4(\Z)$ tel que $\Psi R = R^\vee \Psi$
on se ramène à l'équation (\ref{equa:frakF}) $\mfF=\eps \xi \mfF^*$.
Ses solutions $\mfF \in \GL_m(\mcO)$ sont à déterminer modulo
l'action de $\GL_m(\mcO)$ par $\mfP \cdot \mfF = \mfP \mfF \mfP^*$ 
($\mfP \in \GL_m(\mcO)$) qui correspond à l'action du commutant de
$R$ sur les solutions de $F=RF'$, voir \ref{susb:existence}.

\begin{lemm}
\label{lemm:antidiag}
Soit $\mcO$ l'anneau des entiers cyclotomiques comme ci-dessus et
soit $\tht \in \mcO$ tel que $\tht \bar{\tht}= 1$. On suppose
que $\{z \in \mcO; z=\tht \bar{z}\}=\{w+\tht \bar{w};w \in \mcO\}$.
Soit $\mfF= \smallmat{\al}{\be}{\ga}{\de} \in \GL_2(\mcO)$ vérifiant
$\mfF = \tht \mfF^*$ et $\al \de =0$. 
Alors il existe $\mfP \in \GL_2(\mcO)$ tel que
$\mfP \cdot \mfF = \smallmat{0}{1}{\tht}{0}$.
\end{lemm}

\pr On peut supposer que $\al=0$. On a $\be \in \mcO^\times$ et par
hypothèse il existe $w \in \mcO$ tel que $\de = w+\tht \bar{w}$. Alors
$\mfP=\smallmat{\be^{-1}}{~0}{-w\be^{-1}}{~1}$ convient.

\paragraph{Cas où $R=V_2 \oplus V_2$ (ordre $3$).} 
 On prend $\Psi_0=
\smallmat{1}{-1}{0}{1}$ et  $\Psi=\Psi_0 \oplus \Psi_0$, d'où
$\xi=-\eps^2$ (avec les notations de \ref{susb:existence}).
L'équation (\ref{equa:frakF}) s'écrit donc $\mfF+\mfF^*=0$
et on cherche  $\mfF=\smallmat{\al}{\be}{-\bar{\be}}{\de}\in\GL_2(\mcO)$
avec $\al$ et $\de$ éléments de $(\eps-\overline{\eps})\Z =
(1+2\eps)\Z$. En revenant sur $\Z$, on constate que $2q_F$
(équation~(\ref{equa:FPFt_11})) est une forme paire à 4 variables dont le
déterminant coïncide avec celui de $F$. Si $F \in \GL_4(\Z)$,
on voit donc que~$q_F$ représente 0 (et que $\det(F)=1$). 
Pour $\mfP=\smallmat{u}{v}{w}{x}\in\GL_2(\mcO)$ et $\mfF$ comme
plus haut on peut écrire
$$(\eps-\overline{\eps})q_F=\al |u|^2 -\bar{\be} v \bar{u}
+ \be u \bar{v} + \de |v|^2.$$
Il existe donc $(u,v) \in \mcO^2$, $(u,v)\neq (0,0)$ avec
$\pgcd(u,v)=1$ qui annule $q_F$. Comme $\mcO$ est principal,
$(u,v)$ est la première ligne d'un élément de $\GL_2(\mcO)$.
Par conséquent, on peut supposer que 
$\mfF=\smallmat{0}{\be}{-\bar{\be}}{\de}$.
Alors (lemme \ref{lemm:antidiag}) $\mfF$ est équivalente à $J_2$
et il existe un unique type associé à $R$. Le représentant $H_4$
(équation~(\ref{equa:indec_rg4}), p.~\pageref{equa:indec_rg4})) correspond à
$\mfF=\smallmat{0}{\eps}{~-\bar{\eps}}{~0}$.

\paragraph{Cas où $R=J_2 \oplus J_2$ (ordre $4$).}  On prend $\Psi_0=I_2$,
$\Psi=I_4$, d'où $\xi =1$. L'anneau $\mcO=\Z[i]$ est muni de la 
$\Z$-base $(1,i)$, de sorte que  $R$ correspond à la multiplication
par $-i$ dans $\mcO^2$ ; l'équation (\ref{equa:frakF}) s'écrit 
$\mfF+ i \mfF^*=0$ et on cherche $\mfF=
\smallmat{\al}{~\be}{-i\bar{\be}}{~\de} \in\GL_2(\mcO)$ avec
$\al$,  $\de \in (1-i)\Z$.
Noter que $\det(\mfF)=\pm i$ est un invariant du type.
Si $\al \de =0$, d'après le lemme \ref{lemm:antidiag},
$\mfF$ est équivalente dans $\GL_2(\mcO)$ à $\smallmat{0}{1}{-i}{0}$
qui correspond au type $K_4$ (voir  équation~(\ref{equa:indec_rg4})
p.~\pageref{equa:indec_rg4}). Nous déterminons les types restants
en nous inspirant de la méthode classique de réduction des 
formes binaires (cf. \cite[p. 14]{Buell1989a}). Pour 
$\mfP=\smallmat{0}{1}{1}{x}$, on a $\mfP \cdot \mfF =
\smallmat{\de~}{\de \bar{x}-i\bar{\be}}{$*$}{$*$}$, où l'on
peut choisir $x\in \mcO$ tel que $|\de \bar{x}-i\bar{\be}|
\leq |\de|/\sqrt{2}$. En supposant que les termes diagonaux sont
toujours non nuls, on est ainsi ramené au bout d'un nombre 
fini d'étapes à $\sqrt{2} |\be| \leq |\al| \leq |\de|$ ;
comme $\det(\mfF)=\pm i$, on a $|\al \de | \leq 1+ |\be|^2$
et forcément $|\al \de | \leq 2$. Ce qui conduit à 
$\al=\de = \pm (1-i)$ et $|\be|=1$. Enfin $\be$ peut être fixé
arbitrairement, par exemple $\be=-i$ qui donne $\pm L_4$
(équation~(\ref{equa:indec_rg4}) p.~\pageref{equa:indec_rg4}). Le cas
précédent $R=V_2 \oplus V_2$ se traite aussi par cette méthode. 

\begin{table}[t]
\centering
\renewcommand{\arraystretch}{1.2}
\renewcommand{\tabcolsep}{0.6 em}%
{\small
\begin{tabular}{llll}
\hline
$d$  & $F F^\vee$& $F$ & type réel \\ \hline
1   & $I_4$ & $I_{p,4-p}$  &$I_{p,4-p}$\\
    & &$U_2  U_2$ & $I_{2,2}$ \\
\hline
2   & $-I_4$ & $J_4$ & $J_2^2$\\
 & $I_{2,2}$ &  $\pm I_2  J_2$ & $\pm I_2  J_2$\\
 & &  $I_{1,1}  J_2$&   $I_{1,1}  J_2$\\
 & &  $U_2 J_2$&   $I_{1,1}  J_2$ \\
 &  $R_2  R_2$ & $F_4^*$ & $I_{1,1}  J_2$\\
\hline
3 & $I_1  R_3$  & $\pm I_{p,1-p}  F_3$ &  $\pm I_{p+1,1-p}  P_{2;3,1}^-$\\
  & & $ \pm G_4^*$ &  $\pm I_2  P_{2;3,1}$\\
  &$V_2  V_2$ & $H_4^*$ & $P_{2;3,1}P_{2;3,1}^-$\\
\hline
4 & $J_2  J_2$&  $K_4^*$&  $P_{2;4,1}P_{2;4,1}^-$\\
  & &$\pm L_4^*$  & $\pm P_{2;4,1}^2$\\
  & $I_1  X_3$ & $\pm I_{p,1-p} G_3$ &  $\pm I_{p+1,1-p}  P_{2;4,1}^-$\\
  & &  $\pm I_1   H_3$ & $\pm I_2  P_{2;4,1}$\\
\hline
6 &$I_2  W_2$  & $\pm I_{p,2-p} F_2$   &  $\pm I_{p,2-p} P_{2;6,1}$\\
  & &$\pm U_2  F_2$& $\pm I_{1,1}P_{2;6,1}$\\
  &$I_2^- W_2$  & $\pm J_2  F_2$ &  $\pm J_2  P_{2;6,1}$\\
  &$W_2  W_2$  & $\pm F_2  F_2$ & $\pm P_{2;6,1}^2$\\
  & &$F_2 F_2^-$& $P_{2;6,1} P_{2;6,1}^-$\\
\hline
10 &$W_4$  & $\pm M_4^*$&   $\pm P_{2;10,1}P_{2;10,3}$\\
   & & $\pm N_4^*$& $\pm P_{2;10,1}P_{2;10,3}^-$\\
\hline
12 & $Y_4$  & $\pm O_4^*$ & $\pm P_{2;12,1}P_{2;12,5}^-$\\
\hline
\end{tabular}}
\caption{Types isoduaux de rang $4$ (47 types, 21 principaux)}
\label{tabl:types_4}
\end{table}

\paragraph{Cas où $R=W_2 \oplus W_2$ (ordre $6$).}  Pour $\Psi_0=
\smallmat{1}{1}{0}{1}$ on trouve $\xi=1-\eps =\bar{\eps}$ et
$\mfF=\mfF^*$. Le lemme \ref{lemm:antidiag} s'applique à nouveau,
de  sorte que toute solution $\mfF=
\smallmat{\al}{\be}{\bar{\be}}{\de} \in \GL_2(\mcO)$ avec 
$\al \de = 0$ est équivalente à $U_2$.
Mais $U_2$ est équivalente par  $\mfP=
\smallmat{~\eps}{1}{-\bar{\eps}}{1} \in \GL_2(\mcO)$ à $I_{1,1}$,
qui correspond au type $F_2 \oplus (-F_2)$. En supposant maintenant
que les termes diagonaux sont toujours non nuls, on se ramène par 
un procédé de réduction comme dans le cas précédent à
 $\sqrt{3} |\be| \leq |\al| \leq |\de|$. Puisque $|\det(\mfF)|=1$,
on voit que $0 < |\al \de | \leq 3/2$, d'où $\al \in \mcO^\times$.
Grâce au  lemme de scindement \ref{lemm:scind} ($\al$ correspond
au premier bloc de $F$), ou par une vérification directe, 
on trouve  encore des types décomposés. 
Tous les types associés à $W_2 \oplus W_2$ sont donc scindés
(voir table~\ref{tabl:types_4}).

\paragraph{Cas où $R=W_4$ (ordre $10$).} Comme $m=1$, l'équation
(\ref{equa:frakF}) se réduit à une relation scalaire dans $\mcO$.
Soit $\Psi=(\psi_{i,j}) \in \GL_4(\Z)$  de coefficients
$\psi_{i,j}=1$ si $0 \leq j-i \leq 1$, $\psi_{i,j}=0$ sinon.
On a $\xi=\bar{\eps}$ ; il s'agit donc de trouver les unités
réelles $\mff \in \mcO^\times$ modulo multiplication  par
$u\bar{u}$ ($u \in \mcO^\times$). Toute unité de $\mcO$ s'écrit
$u=\tht \nu^k$ ($k \in\Z$) où $\tht^{10}=1$ et $\nu = \eps -\eps^4$
est une unité fondamentale réelle. Le sous-groupe des $u\bar{u}$
($u \in \mcO^\times$) est engendré par $\nu^2$ et on trouve 4
solutions $\mff=\pm 1$, $\pm \nu$ qui correspondent respectivement
aux types  $\pm M_4$ et $\pm N_4$ (équation~(\ref{equa:indec_rg4}) de
l'annexe).

\paragraph{Cas où $R=Y_4$ (ordre $12$).}  On pose $\Psi=(\psi_{i,j})
 \in \GL_4(\Z)$ avec $\psi_{i,j}=1$ si $0 \leq j-i =0$ ou $2$,
$\psi_{i,j}=0$ sinon. Alors $\xi=-\eps^4$ et (\ref{equa:frakF})
devient 
\begin{equation}
\label{equa:Y4}
\eps \mff = \bar{\mff}\esp  (\mff \in \mcO^\times).
\end{equation}
On a $\mcO^\times = \{\eps^l \nu^k ; (k,l) \in \Z^2\}$, où
$\nu=\eps^2-\eps^3$ est une unité fondamentale vé\-ri\-fi\-ant
(\ref{equa:Y4}). Il est alors immédiat de constater que toute
solution de (\ref{equa:Y4}) s'écrit $\pm \nu u \bar{u}$ 
($u\in \mcO^\times$). Les solutions  $\pm \nu$ correspondent aux 
types $\pm O_4$ (équation~(\ref{equa:indec_rg4}) de l'annexe).

\subsection{Types algébriques principaux de rang 5}
\label{subs:typalg5}
Dans $\GL_5(\Z)$ il y a  20 classes de conjugaison d'éléments
d'ordre $2^k$ et de déterminant 1, décrites  dans la table 
\ref{tabl:of5-7}. Pour $R=I_{2,1}R_2$, $I_3^-R_2$, $J_2X_3$,
$R_2X_3^-$, $I_{1,1}X_3^-$, $I_1Z_4$, $I_1Z'_4$ et $I_1Z''_4$,
une vérification directe montre que le déterminant d'une solution
entière $F$ de $F=RF'$ est divisible par 2. Pour la même raison,
ou en utilisant la proposition~\ref{prop:dec1} et l'exemple 
\ref{exem:epsxi}, on peut aussi  éliminer $I_3J_2$,  $I_{1,2}J_2$,
$I_1^-R_2J_2$ et $I_1X_4$. Ensuite, $I_5$ donne les types 
symétriques et $I_{3,2}$, $I_{1,4}$, $I_1J_2J_2$ et $I_2^-X_3$
ne donnent que des types scindés (proposition~ \ref{prop:dec1}), 
sans relations nouvelles (voir table~\ref{tabl:types_5}).  
Nous avons également prouvé (\S \ref{subs:typalg1-4}, cas $I_1X_3$
p.~\pageref{para:R=I_1X_3}) que tous les types associés à
$I_2X_3$ sont scindés ; compte tenu de la relation (\ref{equa:I1^-G3}),
on  trouve ainsi 12 types (table~\ref{tabl:types_5}). 
Finalement, les seuls éléments qui restent à examiner sont
$I_1R_2R_2$  et $X_5$. 

\paragraph{Cas où $R=I_1 \oplus R_2 \oplus R_2$ (ordre $2$).}
\label{para:R=I1R2R2}
Pour traiter les dimensions $n \geq 5$, il  est utile de considérer
plus généralement $R=I_q \oplus R_2 \oplus R_2$ où $q=n-4 \geq 1$.
Soit $\al=(1,1)$. Avec les notations du  \S \ref{subs:typalg1-4},
on a alors
\begin{equation}
\label{equa:IqR2R2}
{\small 
P=\left(
\begin{array}{ccc}
U &  v \al & w \al\\
\al' x'&  V  & W \\
\al' y'&  X  & Y 
\end{array}
\right)
\esp
{\rm et}
\esp
F=\left(
\begin{array}{ccc}
A &  b \al & c \al\\
\al' b'&  d(I_2+R_2) & f I_2 + g R_2\\
\al' c'&  g I_2 + f R_2 & e(I_2+R_2)\\
\end{array}
\right),}
\end{equation}
avec $A,~U \in M_q(\Z)$ et $A$ symétrique, $b$, $c$, $v$, $w$, $x$,
$y \in \Z^q$ (vecteurs colonnes), $d$, $e$, $f,~g \in \Z$ 
et $V$, $W$, $X$, $Y$  polynômes de degré au plus 1  en $R_2$.
Soit $\De=I_q \oplus \De_4$ (voir équation~(\ref{equa:De3_De4})), de sorte
que $\De R \De^{-1} = I_{q+2,2}$.
On vérifie que $\De C_R \De^{-1}$ coïncide
avec le sous-groupe $\Ga^{n-2,2}$ de $\GL_n(\Z)$ défini dans le lemme
suivant.

\begin{lemm}
\label{lemm:Ga_0^kl}
Pour $0<l<k$, notons $\Ga_0^{k,l}$ le groupe  des %
$\smallmat{A}{B}{C}{D} \in \GL_k(\Z)$ où
$C \in M_{l,k-l}(\Z)$  est congrue à $0$ modulo~2. Posons
également $\Ga_0^{k,k}=\GL_k(\Z)$. Soient $r_1$ et $r_2$ $\geq 2$
et soit  $\Ga^{r_1,r_2}$ le groupe  des matrices
$Q_1 \oplus Q_2  \in \GL_{r_1+r_2}(\Z)$
avec $Q_j=\smallmat{A_j}{B_j}{C_j}{D_j} \in \Ga_0^{r_j,2}$ ($j=1,2$) et 
$D_1 \equiv D_2$ modulo~2. Alors les morphismes \og composantes\fg\
de $\Ga^{r_1,r_2}$ dans $\Ga_0^{r_j,2}$ ($j=1,2$) sont surjectifs.
\end{lemm}

\pr Étant donné $Q_1=\smallmat{A_1}{B_1}{C_1}{D_1} \in \Ga_0^{r_1,2}$,
il existe $D_2 \in M_2(\Z)$ dont les coefficients valent 0 ou 1 et telle
$D_2 \equiv D_1$ modulo  2. Mais $D_1$ est inversible modulo~2, donc
$\det(D_2)=\pm 1$. Par construction, on a  
$Q_1 \oplus I_{r_2-2} \oplus D_2 \in \Ga^{r_1,r_2}$.

\medskip
Posons maintenant $\De \cdot F = A_F \oplus B_F$. 
On  a alors $B_F=2(f-g)J_2$ et 
\begin{equation}
\label{equa:AF_IqR2R2}
{\small 
A_F=\left(
\begin{array}{ccc}
A  &  2b & 2c \\
2b' &  4d  & 2f+2g \\
2c' & 2f+2g & 4e
\end{array}
\right).}
\end{equation}
Ainsi $\det(A_F)$ est divisible par 4 ; vu que
$(f-g)^2 \det(A_F)=4\det(F)=\pm 4$, on a ${(f-g)}^2=1$ et
$\det(A_F)=\pm 4$. En mettant 2 en facteur dans les  deux dernières
 colonnes de $A_F$, on voit que $\det(A)$ est impair. Si $M$ est
 une matrice carrée symétrique et  inversible, on la relation
\begin{equation}
\label{equa:dec_orth}
{\small
\left(
\begin{array}{cc}
I  &  0 \\
-B'M^{-1} &  I \\
\end{array}
\right)
\cdot
\left(
\begin{array}{cc}
M  &  B \\
B' &  C \\
\end{array}
\right)  
=
\left(
\begin{array}{cc}
M  &  0 \\
0 &  C-B'M^{-1}B \\
\end{array}
\right).}
\end{equation}
Comme $A$ est inversible sur l'anneau $\Z_2$ des entiers 2-adiques,
on en déduit que la décomposition de
Jordan 2-adique de $A_F$ est de la forme $f_1 \oplus 2 f_2$, avec
$f_1$ de rang $q$ et  $f_2$ {\em paire} de rang 2. De plus, 
si $A_F$ est équivalente par un élément $Q_1 \in \GL_{n-2}(\Z)$
à  une matrice $A^0$ de la forme (\ref{equa:AF_IqR2R2}),
un calcul modulo~2 montre que nécessairement $Q_1 \in \Ga_0^{n-2,2}$.
Grâce au lemme \ref{lemm:Ga_0^kl}, il existe $P\in C_R$
tel que  $A_{P \cdot F}=A^0$. En particulier la classe
entière de $A_F$ caractérise le type $[F]$.
\par
Revenons au cas $q=1$. D'après ce qui précède, $A_F$ appartient au
genre ${\rm I}_{r,s}(2^2_{\rm II})$, qui est nécessairement indéfini
en rang $r+s=3$  (cf. \cite[pp. 386-387]{Conw-Sloa1999a}). On sait
de plus que les genres indéfinis de
rang $\geq 3$ et de petit déterminant ne contiennent qu'une seule
classe d'équivalence entière. Donc $A_F$ est équivalente à 
$\pm I_1\oplus 2 U_2$ qui représente le genre
${\rm I}_{r,s}(2^2_{\rm II})$ suivant que  $(r,s)=(2,1)$ ou $(1,2)$.
On peut donc supposer (voir ci-dessus)  que $A_F=\pm I_1\oplus 2 U_2$
et finalement $F$ est équivalente à un type décomposé
$\pm (I_1 \oplus F_4)$.

\paragraph{Cas où $R=X_5$ (ordre $8$).}  
\label{para:R=X5}
Considérons plus généralement $R=I_q \oplus X_5$ ($q \geq 0$). 
Soit $\eps$ une racine primitive 8-ième de l'unité et soit $\mcO$
l'anneau des entiers du corps cyclotomique $\Q[\eps]$.
Pour $z \in \mcO$, notons $M(z)$ la matrice de la multiplication
par $z$ dans la base $(1,\eps,\eps^2,\eps^3)$. Soit enfin
$\De = I_q \oplus \De_5$ et $\Psi=I_q \oplus \Psi_5$ avec
$$
{\small 
\De_5=\left(
\begin{array}{ccccc}
2 & 1 & 1 & 1 & 1\\
0 & 1 & 0 & 0 & 0\\
0 & 0 & 1 & 0 & 0\\
0 & 0 & 0 & 1 & 0\\
0 & 0 & 0 & 0 & 1
\end{array}
\right)
\esp {\rm et} \esp
\Psi_5=\left(
\begin{array}{ccccc}
2 & 1 & 1 & 1 & 1\\
1 & 0 & 0 & 0 & 1\\
1 & 0 & 0 & 0 & 0\\
1 & 1 & 0 & 0 & 0\\
1 & 1 & 1 & 0 & 0
\end{array}
\right).}
$$
Soit $P\in C_R$. On a $\De P \De^{-1}  = Q_1 \oplus M(\tht)$ où
$Q_1 \in \Ga_0^{q+1,1}$ (groupe défini au  lemme \ref{lemm:Ga_0^kl}) et 
$\tht = v +w \eps +x \eps^2+y \eps^3 \in  \mcO^\times$.
Noter que $Q_1[q+1,q+1] \equiv v+w+x+y \equiv 1$ (mod~2).
Posons $P=P(Q_1,\tht)$.
On vérifie alors que $\Psi^\vee \Psi= P(I_{q+1},-\eps^3)$
et que $P(Q_1,\tht)^\sig = \Psi^{-1} P(Q_1,\tht)' \Psi=
P(Q_1^\tau,\bar{\tht})$, où $\tau$ est l'involution de 
$\Ga_0^{q+1,1}$ ($q>0$) définie par
$$
{\small 
{\left(
\begin{array}{ll}
U & v \\
2w' & x
\end{array}
\right)}^\tau
=
\left(
\begin{array}{ll}
U & w \\
2v' & x
\end{array}
\right),}
$$
ou l'identité si $q=0$. En posant $F\Psi=P(A,\mff)$, l'équation
$F=RF'$ ramenée sous la forme (\ref{equa:G=RxiGsig}) s'écrit
$P(A,\mff) = P(A^\tau,\bar{\mff})$.
\par
Si $q=0$, les types associés à $X_5$ sont
donc paramétrés par les couples $(a,\mff)$ avec $a = \pm 1$ 
et $\mff$ unité réelle de $\Q[\eps]$ modulo multiplication
par $\tht \bar{\tht}$ ($\tht \in \mcO^\times$). On sait que
$\nu = 1+\eps-\eps^3$ ($=1\pm \sqrt{2}$ suivant le choix de $\eps$)
est une unité fondamentale à la fois pour $\Q[\sqrt{2}]$ et 
pour $\Q[\eps]$. Par suite, on peut supposer que $f=\pm 1$ ou  $\pm \nu$
et on trouve 8 types. Les types $F_5$,  $G_5$, $H_5$ et $K_5$
(équation~(\ref{equa:indec_rg5}) et table~\ref{tabl:types_5}),
indécomposables, 
correspondent respectivement à $a=1$ et $\mff = 1$, $-1$,
$\nu$ et $-\nu$.
\par
Si $q>0$, considérons $\De \cdot F = A_F \oplus B_F$. On a
$\De \cdot \Psi^{-1} = I_q \oplus 2I_1 \oplus M(\eps^3-\eps^2)$. La
relation  $\De \cdot F = (\De P(A,\mff) \De ^{-1} )(\De \cdot \Psi^{-1})$
montre alors que $A_F=A(I_q \oplus 2I_1)$ et que
$B_F=M(\mff(\eps^3-\eps^2))$. Par suite $\det(A_F)= \pm 2$ et
$A$ est de la forme $\smallmat{C}{2b}{2b'}{2d}$. 
D'après l'argumentation donnée au \S \ref{subs:typalg1-4}
(cas $R=I_1\oplus X_3$, p.~\pageref{para:R=I_1X_3}),
pour $1 \leq q+1 \leq 6$, il existe  $Q_1\in \GL_{q+1}(\Z)$ tel que
$Q_1 \cdot A = B \oplus (\pm 2I_2)$ et nécessairement 
$Q_1 \in \Ga_0^{q+1,1}$. Comme $P(Q_1,1) \in C_R$, on peut
supposer que $A = B \oplus (\pm 2I_2)$ et finalement  tous les types 
associés à $R$ sont scindés.

\medskip
\begin{table}[!h]
\centering
\renewcommand{\arraystretch}{1.2}
\renewcommand{\tabcolsep}{0.6 em}%
{\small
\begin{tabular}{llll}
\hline
$d$  & $F F^\vee$ & $F$ & type réel\\ \hline
1   & $I_5$ & $I_{p,5-p}$ & $I_{p,5-p}$\\
\hline
2 & $I_{3,2}$ & $I_{p,3-p}  J_2$ &$I_{p,3-p}  J_2$\\
  & $I_{1,4}$ & $\pm I_1 J_4$   &$\pm I_1 J_2^2$\\ 
  & $I_1 R_2 R_2$ & $\pm I_1 F_4$  &$I_{2,1}J_2$\\ 
\hline
4 & $I_1 J_2  J_2$ & $\pm I_1 K_4$ & $\pm I_1 P_{2;4,1}P_{2;4,1}^-$ \\ 
  & & $\pm I_1 L_4$ & $\pm I_1 P_{2;4,1}^2$ \\
  & & $\pm I_1^- L_4$  & $\pm I_1^- P_{2;4,1}^2$ \\
  & $I_2  X_3$ & $\pm I_2 G_3$   & $\pm I_3 P_{2;4,1}^-$\\
  && $\pm I_2 G_3^-$   & $\pm I_{2,1} P_{2;4,1}$\\
  && $\pm I_2 H_3$   & $\pm I_{3} P_{2;4,1}$\\
  && $\pm I_2^- H_3$ & $\pm I_{1,2} P_{2;4,1}$\\
  && $\pm U_2   G_3$   & $\pm I_{2,1} P_{2;4,1}^-$\\
  && $\pm U_2   H_3$   & $\pm I_{2,1} P_{2;4,1}$\\
  & $I_2^- X_3$ & $\pm J_2 G_3$  & $\pm J_2 I_1 P_{2;4,1}^-$ \\
  && $\pm J_2 H_3$  & $\pm J_2 I_1 P_{2;4,1}$ \\ 
\hline
8   & $X_5$ & $\pm F_5^*$    & $\pm I_1 P_{2;8,1}^- P_{2;8,3}$ \\ 
  && $\pm G_5^*$   & $\pm I_1 P_{2;8,1}P_{2;8,3}^-$ \\  
  && $\pm H_5^*$  & $\pm I_1 P_{2;8,1}^-P_{2;8,3}^-$ \\
  && $\pm K_5^*$  & $\pm I_1 P_{2;8,1}P_{2;8,3}$ \\
\hline
\end{tabular}}
\caption{Types isoduaux principaux de rang $5$ (44 principaux)}
\label{tabl:types_5}
\end{table}

À nouveau, les types sont distingués par leurs composantes 
réelles et au besoin par la parité de
leur composante symétrique  (voir proposition~\ref{prop:dec_can}). 
Par exemple, $I_2G_3^-$ et
$U_2H_3$ ont le même type réel $[I_{2,1}P_{2;4,1}]$ mais
leur composantes symétriques  équivalent respectivement à  
$(\Z^3,x^2+y^2 - 2z^2)$ et à  $(\Z^3,2xy+2z^2)$.

\subsection{Types algébriques principaux de rang 6}
\label{subs:typalg6}
À conjugaison près dans $\GL_6(\Z)$, il y a 39 éléments $R$ 
d'ordre $2^k$ et de déterminant 1 (table~\ref{tabl:of5-7}).
On peut écarter les $R$ qui contiennent un seul facteur $R_2$,
$J_2$, $X_3^-$, $X_4$, $Z_4$, $Z'_4$, $Z''_4$ ou $X_5^-$ car
$\det(F)$ est divisible par 2 pour toute solution de $F=RF'$
(certains cas relèvent aussi de la proposition~\ref{prop:dec1} et ex.
\ref{exem:epsxi}). De même on peut éliminer $J_2J_2J_2$ et $X_6$.   
Les types associés à $I_{p,q}$, $I_{p,q}J_2J_2$ ou $I_{1,2}X_3$
sont connus ou scindés. Il en est  de même pour $I_3X_3$
(\S \ref{subs:typalg1-4}, cas $I_1X_3$, p.~\pageref{para:R=I_1X_3})
qui donne 10 types 
et pour $I_1X_5$ (\S \ref{subs:typalg5}, cas $X_5$,
p.~\pageref{para:R=X5})  avec 2 relations 
\begin{equation}
[I_1 \oplus (-F_5)] = [(-I_1) \oplus G_5]
\esp {\rm et} \esp
[I_1 \oplus (-H_5)] = [(-I_1) \oplus K_5].
\end{equation}
Ces relations sont obtenues grâce à $P(Q,1)$ où
 $Q=\smallmat{1}{1}{2}{1}$ (notation du cas $R=X_5$) et 
on trouve seulement 12 types (table~\ref{tabl:types_6}).
Il ne reste plus qu'à examiner $I_2^\pm R_2R_2$, $X_3X_3$,
$X_3^-X_3^-$ et $Y_6$.

\paragraph{Cas où $R=I_2 \oplus R_2 \oplus R_2$ (ordre $2$).} 
Soit $F \in \GL_6(\Z)$ une solution de $F=RF'$ et soit 
$A_F$ la composante symétrique de $F$ (équations~(\ref{equa:IqR2R2})
et~(\ref{equa:AF_IqR2R2})). On sait  (\S \ref{subs:typalg5})
que la décomposition de Jordan 2-adique de $A_F$  est de la forme
$f_1 \oplus 2 f_2$, avec $f_2$ {\em paire} de rang 2. La forme
$A_F$ appartient donc l'un des genres ${\rm I}_{r,s}(2^2_{\rm II})$,
${\rm II}_{r,s}(2^2_{\rm II})$ ou ${\rm II}_{r,s}(2^{-2}_{\rm II})$ 
avec $r+s=4$ (cf. \cite[pp. 386-387]{Conw-Sloa1999a}). Les deux
premiers genres sont indéfinis, respectivement représentés par l'unique
forme $I_{p,q}\oplus 2 U_2$ ou  $U_2 \oplus 2 U_2$ à équivalence près.
Le genre ${\rm II}_{r,s}(2^{-2}_{\rm II})$ ($r+s=4$) est défini
car $r-s \equiv 0$ (mod 8); dans ce cas, d'après
\cite{Kork-Zolo1872a}, la forme $A_F$ réalise le maximum de l'invariant
d'Hermite en dimension 4 et est équivalente $\pm D_4$, avec
$$
{\small 
D_4=\left(
\begin{array}{ccccc}
2 & 1 & 0 & 0 \\
1 & 2 & 0 & 2 \\
0 & 0 & 4 & 2 \\
0 & 2 & 2 & 4 
\end{array}
\right).}
$$
Grâce au lemme \ref{lemm:Ga_0^kl}, on peut supposer (voir \S
\ref{subs:typalg5}) que $A_F= A \oplus 2U2$ avec $\det(A)=\pm 1$
auquel cas le type $[F]$ est scindé,  ou bien que $A_F= \pm D_4$
qui conduit au type indécomposable $\pm F_6$ (voir
 table~\ref{tabl:types_6}).

\paragraph{Cas où $R=(-I_2) \oplus R_2 \oplus R_2$ (ordre $2$).}
\label{para:R=I2^-R2R2}
Plus généralement, considérons $R=I_q\oplus (-I_2) \oplus R_2 \oplus R_2$
($q \geq 0$). En posant $\al=(1,1)$, $\be=(1,-1)$, on a 
$$
{\small 
P=\left(
\begin{array}{cccc}
T & 0  & v \al & w \al \\
0 &  U & v_1 \be & w_1 \be\\
\al' x' & \be' x'_1&  V  & W \\
\al' y' & \be' y'_1&  X  & Y 
\end{array}
\right)
~~
{\rm et}
~~
F=\left(
\begin{array}{cccc}
A & 0 & b \al & c \be\\
0 & \Om &  b_1 \be & c_1 \be\\
\al' b' & -\be' b'_1&  d(I_2+R_2) & f I_2 + g R_2\\
\al' c' &-\be' c'_1&  g I_2 + f R_2 & e(I_2+R_2)\\
\end{array}
\right),}
$$
avec $A,~T \in M_q(\Z)$ et $A$ symétrique, $\Om,~U \in M_2(\Z)$
et $\Om$  antisymétrique, $b$, $c$, $v$, $w$, $x$, $y \in \Z^q$,
$b_1$, $c_1$, $v_1$, $w_1$, $x_1$, $y_1 \in \Z^2$, $d$, $e$,
$f,~g \in \Z$  et $V$, $W$, $X$, $Y$  polynômes en $R_2$. La matrice 
$\De=[I_q \oplus (U_2 \otimes I_2) \oplus I_2](I_{q+2} \oplus \De_4)$
(voir équation~(\ref{equa:De3_De4})) vérifie $\De R \De^{-1} = I_{q+2,4}$,
$\det(\De)=4$, $\De C_R \De^{-1} = \Ga^{q+2,4}$ (voir lemme
\ref{lemm:Ga_0^kl}) et $\De \cdot F= A_F \oplus B_F$ avec
\begin{equation}
\label{equa:AFBF}
{\small
A_F=\left(
\begin{array}{ccc}
A  &  2b & 2c \\
2b' &  4d & 2(f+g) \\
2c' & 2(f+g) &  4e \\
\end{array}
\right)
~~ {\rm et} ~~ 
B_F=\left(
\begin{array}{ccc}
\Om  &  -2b_1  & 2c_1\\
2b'_1 &  0  & 2g-2f\\
-2c'_1  & 2f-2g  & 0\\
\end{array}
\right).}
\end{equation}
Les entiers $\det(A_F)$ et $\det(B_F)$ sont divisibles par 4, donc
valent $\pm 4$. Comme au \S \ref{subs:typalg5} (cas $I_q
\oplus R_2 \oplus R_2$, p.~\pageref{para:R=I1R2R2}) on prouve 
que $\det(A)$ est impair et que la $A_F$ appartient à l'un des 
genres ${\rm I}_{r,s}(2^2_{\rm II})$,
${\rm II}_{r,s}(2^2_{\rm II})$ ou ${\rm II}_{r,s}(2^{-2}_{\rm II})$. 
\par
Si $q=0$, la classe de $B_F$ détermine $[F]$.
En effet, il existe $Q_2 \in \GL_4(\Z)$ tel que
$Q_2 \cdot B_F = J_2 \oplus 2 J_2$ 
(\cite[\S 5 n$^{\rm o}$ 1]{BourbakiAlg9-1959a}). 
Nécessairement $Q_2 \in \Ga_0^{4,2}$ et (lemme \ref{lemm:Ga_0^kl})
il existe $Q_1 \in \GL_2(\Z)$ tel que
$\De^{-1} (Q_1 \oplus Q_2) \De \in C_R$.
Ainsi on peut supposer que $B_F=  J_2 \oplus 2 J_2$~:  $F$ est
scindé et $[F]=[J_2 \oplus F_4]$.
 
\paragraph{Cas où $R=X_3 \oplus X_3$ (ordre $4$).} 
Pour $\De=(I_1\oplus R_3 \oplus I_2)(\De_3 \oplus \De_3)$
(voir équation~(\ref{equa:De3_De4})) on a
$\De R \De^{-1} =I_2 \oplus J_2 \oplus J_2$.  Si $P \in C_R$,
alors $\De P \De^{-1}=Q_1 \oplus Q_2$ où $Q_1\in \GL_2(\Z)$ et
$Q_2$ s'identifie à un élément $\mfP \in \GL_2(\Z[i])$ (après
choix de $(1,i)$ comme $\Z$-base de $\Z[i]$) ; de plus
$Q_1 \equiv \re \mfP + \im \mfP$ (mod 2). Dans ces conditions
nous posons $P=P(Q_1,\mfP)$.  Cela étant, la matrice $\Psi
=\Psi_3 \oplus \Psi_3$ avec
\begin{equation}
\label{equa:Psi3}
{\small
\Psi_3 =\left(
\begin{array}{rrr}
2 &  -1 &  1 \\
-1 &  0 &  -1 \\
1 &  0 &  0
\end{array}
\right)}
\end{equation}
vérifie $\Psi R = R^\vee \Psi$,  $\Psi^\vee \Psi = R$
et $P(Q_1,\mfP)^\sig = \Psi^{-1} P(Q_1,\mfP)' \Psi = P(Q_1',\mfP^*)$.
Soit $F$ une solution inversible de $F=RF'$. En posant
$F\Psi=P(A,\mfF) \in C_R$, l'équation (\ref{equa:G=RxiGsig}) se traduit
par $A=A'$ et $\mfF +\mfF^*=0$ ; l'action de 
$P=P(Q_1,\mfP)$ est donnée par $P P(A,\mfF) P^\sig =
P(Q_1\cdot A, \mfP \cdot \mfF)$ où $\mfP \cdot \mfF=\mfP \mfF \mfP^*$.
On a $i\mfF=\smallmat{a}{\be}{\bar{\be}}{d}$ avec $(a,d) \in \Z^2$
et $ad-|\be|^2=\pm 1$.  Si $ad=0$, alors $i\mfF$ est équivalente
-- comme forme hermitienne  -- à $U_2$ ou à $I_{1,1}$ selon que 
$a+d$ est pair ou impair (voir par exemple la preuve du lemme
\ref{lemm:antidiag} pour le premier cas). Sinon, en supposant 
les termes diagonaux toujours non nuls, on se ramène par un
procédé standard de réduction (voir \S \ref{subs:typalg1-4}, 
cas $R=J_2\oplus J_2$) à $\sqrt{2} |\be| \leq |a| \leq |d|$,
\cadab à $ |a| = |d| =1$ et $\be=0$. Finalement $i\mfF$ est 
équivalente à $\pm I_2$, $I_{1,1}$ ou $U_2$.

\begin{lemm} 
\label{lemm:GL2(Z[i])}
L'application qui à $\mfP \in \GL_2(\Z[i])$ associe 
$\re \mfP + \im \mfP$ {\rm (mod~2)}  définit un morphisme surjectif
$\Theta$ de $\GL_2(\Z[i])$  sur $\GL_2(\Z/2\Z)$.
\end{lemm} 

\pr Posons $\Theta_0(\mfP)= \re \mfP + \im \mfP$  pour 
$\mfP \in \GL_2(\Z[i])$. On a $\Theta_0(I_2)=I_2$ et
$\Theta_0(\mfP)\Theta_0(\mfQ)-\Theta_0(\mfP\mfQ) =2 \im\mfP \im\mfQ$
(noter que $|\det(\mfP)|^2 \equiv \det(\Theta_0(\mfP))^2$ (mod~2)).
Enfin, la surjectivité est évidente.

Il résulte de ce lemme que le morphisme naturel de $C_R$ dans
$\GL_2(\Z[i])$ est surjectif  : si $\mfP \in \GL_2(\Z[i])$, la matrice 
$Q_1$ des restes modulo 2 de $\Theta_0(\mfP)$ est inversible. On peut
donc supposer que $i \mfF=\pm I_2$, $I_{1,1}$ ou $U_2$. Puisque
$F\Psi=P(A,\mfF) \in C_R$, on a $A \equiv i \mfF$ (mod~2) et il
existe $Q_1 \in \GL_2(\Z)$ telle que $Q_1\cdot A = U_2$ si
$i \mfF = U_2$, $Q_1\cdot A = I_{p,q}$ ($p+q=2$) dans les autres
cas. On vérifie ensuite que le morphisme $\Theta$ du lemme 
\ref{lemm:GL2(Z[i])} est encore surjectif en restriction au
groupe $\mcI$ des isométries hermitiennes de $U_2$ (par exemple
$\smallmat{i}{1}{0}{i} \in \mcI$). Quand $A$ est paire ($i \mfF = U_2$),
la matrice $Q_1$ se complète ainsi en $P(Q_1,\mfP) \in C_R$ avec
$\mfP \in \mcI$ ; on trouve un seul type (indécomposable) représenté
par  $G_6=\smallmat{0}{G_3}{G_3}{0}$ vérifiant $G_6\Psi=P(U_2,-iU_2)$.
Quand $A$ est impaire, $Q_1$ se complète de même en
$P(Q_1,\mfP) \in C_R$ avec $\mfP$ isométrie hermitienne de $\mfF$,
par exemple $\mfP=\smallmat{1+i}{1}{1}{1-i}$ pour $i\mfF=I_{1,1}$
et $Q_1 \equiv U_2$ (mod~2) -- les autres cas sont évidents. On 
peut donc supposer que $A$ et $\mfF$ sont diagonales, ce qui  
donne 9 types (suivant les classes de $A$ et $\mfF$), tous
scindés (table \ref{tabl:types_6}). Cela implique l'existence
d'une  relation~:
\begin{equation}
\label{equa:H3H3^-=G3G3^-}
[H_3 \oplus (-H_3) ] = [G_3 \oplus  (-G_3)],
\end{equation}
l'équivalence étant donnée par $P(I_2,\mfP)$ avec  
$\mfP=\smallmat{1}{i-1}{~i+1}{i}$.

\paragraph{Cas où $R=(-X_3) \oplus (-X_3)$ (ordre $4$).} 
Par rapport au cas précédent, il suffit de changer le
signe de $R$. On cherche donc $F\Psi=P(A,\mfF) \in C_R$ avec
$A$ antisymétrique et $\mfF$ hermitienne, en particulier 
$A$ et $\mfF$ sont \og paires \fg. Nécessairement $A=\pm J_2$
et (lemme) on peut supposer que $A = J_2$.  Comme la réduction
modulo 2 de $\Sp_2(\Z)=\SL_2(\Z)$ dans $\GL_2(\Z/2\Z)$ est surjective,
on peut supposer aussi que $\mfF=U_2$. Cela conduit à un seul
type (indécomposable) représenté par $H_6=\smallmat{0}{H_3}{-H_3}{0}$
(pour lequel $A=J_2$ et $\mfF=\smallmat{0}{i}{-i}{0}$).

\begin{table}[!tbp]
\centering
\renewcommand{\arraystretch}{1.1}
\renewcommand{\tabcolsep}{0.6 em}%
{\small
\begin{tabular}{llll}
\hline
$d$  & $F F^\vee$ & $F$ & type réel\\ \hline
1 & $I_6$ & $I_{p,6-p}$ & $I_{p,6-p}$\\
  &       & $U_6$ & $I_{3,3}$\\
\hline
2 & $I_6^-$ & $J_6$  &$J_2^3$\\
  &$I_{4,2}$ & $I_{p,4-p} J_2$  & $I_{p,4-p} J_2$\\
  & & $U_4 J_2$ &  $I_{2,2} J_2$\\
  &$I_{2,4}$ & $I_{p,2-p} J_4$ & $I_{p,2-p} J_2^2$\\
  & & $U_2 J_4$ & $I_{1,1} J_2^2$\\
  & $I_2 R_2 R_2$ &  $I_{p,2-p} F_4$ & $I_{p+1,3-p}J_2$\\
  & & $U_2 F_4$ & $I_{2,2} J_2$ \\
  & & $\pm F_6^*$ &  $\pm I_4 J_2$ \\
  &$I_2^- R_2 R_2$ & $J_2 F_4$ &  $I_{1,1} J_2^2$\\
\hline
4 & $I_2 J_2 J_2$ & $I_{p,2-p} K_4$ & $I_{p,2-p}P_{2;4,1}P_{2;4,1}^-$\\
  & &$I_{p,2-p} L_4$ &  $I_{p,2-p}P_{2;4,1}^2$\\
  & &$I_{p,2-p} L_4^-$ & $I_{p,2-p}P_{2;4,1}^{-2}$\\
  & & $U_2 K_4$  & $I_{1,1}P_{2;4,1}P_{2;4,1}^-$\\
  & & $\pm U_2 L_4$ &  $\pm I_{1,1}P_{2;4,1}^2$\\
  & $I_2^-  J_2 J_2$ &  $J_2 K_4$ & $J_2P_{2;4,1}P_{2;4,1}^-$\\
  & & $\pm J_2 L_4$  & $\pm J_2 P_{2;4,1}^2$\\
  & $I_3  X_3$ & $I_{p,3-p}  G_3$  & $I_{p+1,3-p} P_{2;4,1}^-$ \\
  & &$I_{p,3-p}  G_3^-$  & $I_{p,4-p} P_{2;4,1}$\\
  & & $\pm I_3 H_3$  & $\pm I_4 P_{2;4,1}$ \\
  & $I_{2,1}^- X_3$ & $\pm J_2 I_1 G_3$  & $\pm J_2 I_2 P_{2;4,1}^-$\\
  & &$\pm J_2 I_1 H_3$ & $\pm J_2 I_2 P_{2;4,1}$\\
  & &$\pm J_2 I_1 G_3^-$  & $\pm J_2 I_{1,1} P_{2;4,1}$\\
   & $X_3  X_3$ &  $\pm(G_3  G_3)$ & $\pm I_2 P_{2;4,1}^{-2}$\\
  & &$\pm G_3  H_3$ & $\pm I_2 P_{2;4,1} P_{2;4,1}^-$\\
  & &$\pm H_3  H_3$  & $\pm I_2 P_{2;4,1}^2$ \\
  & &$\pm G_3  H_3^-$  & $\pm I_{1,1} P_{2;4,1}^{-2}$\\
  & &$G_3 G_3^-$ & $I_{1,1}P_{2;4,1}P_{2;4,1}^-$\\
  & &$G_6^*$  & $I_{1,1}P_{2;4,1}P_{2;4,1}^-$\\
   & $X_3^- X_3^-$ & $H_6^*$  & $J_2P_{2;4,1}P_{2;4,1}^-$\\
\hline
8 & $I_1 X_5$ & $\pm I_1 F_5$  & $\pm I_2 P_{2;8,1}^-P_{2;8,3}$\\
  & &$\pm I_1  F_5^-$ & $\pm I_{1,1} P_{2;8,1}P_{2;8,3}^-$\\
  & &$\pm I_1  G_5$ &  $\pm I_2 P_{2;8,1}P_{2;8,3}^-$\\
  & &$\pm I_1  H_5$ &  $\pm I_2 P_{2;8,1}^-P_{2;8,3}^-$\\
  & &$\pm I_1  H_5^-$ & $\pm I_{1,1} P_{2;8,1}P_{2;8,3}$\\
  & &$\pm I_1  K_5$ &  $\pm I_2 P_{2;8,1}P_{2;8,3}$\\
  & $Y_6$ & $\pm K_6^*$  & $\pm P_{2;4,1}P_{2;8,1}^-P_{2;8,3}^-$\\
  & &$\pm L_6^*$ &  $\pm P_{2;4,1}P_{2;8,1}P_{2;8,3}$\\
  & &$\pm M_6^*$ &  $\pm P_{2;4,1}P_{2;8,1}P_{2;8,3}^-$\\
  & &$\pm N_6^*$ &  $\pm P_{2;4,1}P_{2;8,1}^-P_{2;8,3}$\\
\hline
\end{tabular}}
\caption{Types isoduaux principaux de rang $6$ (88 principaux)}
\label{tabl:types_6}
\end{table}

\paragraph{Cas où $R=Y_6$ (ordre $8$).} 
Soit $\eps$ une racine primitive 8-ième de l'unité, soit $\mcO$
l'anneau des entiers du corps cyclotomique $\Q[\eps]$ muni de 
la $\Z$-base $(1,\eps,\eps^2,\eps^3)$ et soit $M(\al)$
la matrice de la multiplication par $\al \in \mcO$. 
Considérons
$$
{\small 
\De=\left(
\begin{array}{rrrrrr}
1 & -1 & ~0 & ~1 & 0 & -1 \\
1 & 1 & 1 & 0 & -1 & 0 \\
0 & 0 & 1 & 0 & 0 & 0 \\
0 & 0 & 0 & 1 & 0 & 0 \\
0 & 0 & 0 & 0 & 1 & 0 \\
0 & 0 & 0 & 0 & 0 & 1
\end{array}
\right)
\esp {\rm et} \esp
\Psi=\left(
\begin{array}{rrrrrr}
1 & -1 & 0 & 1 & 0 & -1 \\
1 & 1 & 1 & 0 & -1 & 0 \\
1 & 0 & 0 & 0 & 0 & -1 \\
0 & -1 & 0 & 0 & 0 & 0 \\
-1 & 0 & -1 & 0 & 0 & 0 \\
0 & 1 & 1 & -1 & 0 & 0
\end{array}
\right).}
$$
On a $\Psi R =R^\vee \Psi$, $\det(\De)=2$, 
$\De R \De^{-1}=J_2 \oplus X_4$ et 
$\De P  \De^{-1}=\smallmat{p}{-q}{q}{p} \oplus M(\tht)$ pour
$P\in C_R$, avec $p+iq \in \Z[i]^\times$ et $\tht \in \mcO^\times$.
Posons $P=P(p+iq,\tht)$. Le commutant $C_R$ est isomorphe à
$\Z[i]^\times \times \mcO^\times$ et  l'involution $\sig$ correspond
à la conjugaison :  $P(z,\tht)^\sig=P(\bar{z},\bar{\tht})$.
Sachant que $\Psi^\vee \Psi = P(i,-\eps^3)$, la relation 
(\ref{equa:G=RxiGsig}) page~\pageref{equa:G=RxiGsig} montre que 
les solutions inversibles de $F=RF'$ sont données par
$F\Psi = P(a,\mff)$ avec $a=\pm 1$ et $\mff$ unité réelle de $\mcO$.
Comme pour le cas de $X_5$ (voir \S \ref{subs:typalg5}), on se ramène à 
$\mff= \pm 1$ ou $\pm \eta = \pm ( 1-\eps+\eps^3)$, ce qui conduit
à 8 types indécomposables $\pm K_6$,  $\pm L_6$, $\pm M_6$,
$\pm N_6$ (voir équation~(\ref{equa:indec_rg6}), annexe) 
correspondant respectivement à $a=\pm 1$ et $\mff = 1$, $-1$,
$\eta$ et $-\eta$.

\subsection{Types algébriques principaux de rang 7}
\label{subs:typalg7}
Dans $\GL_7(\Z)$ les classes de conjugaison d'ordre $2^k$  et de
déterminant 1  sont au nombre de 72 (table \ref{tabl:of5-7}).
On a $\det(F) \equiv 0$ (mod~2) si $F=RF'$ dans les cas
suivants : $R$ contient un seul facteur $R_2$, $J_2$, $X_3^-$,
ou contient un facteur $X_4$, $Z_4$, $Z'_4$, $Z''_4$, $X_5^-$,
$R_2R_2R_2$, $J_2J_2J_2$, $X_6$, $X_7$, $Y_7$, $Z'_7$ ou 
$Z''_7$ (voir aussi  la proposition~\ref{prop:dec1} et 
l'exemple~\ref{exem:epsxi}).
Les types associés à $I_{p,q}$ ($q>0$), $I_{p,q}J_2J_2$, $I_{2,2}X_3$,
$I_4^-X_3$, $I_1X_3^-X_3^-$, $I_2^-X_5$ et $I_1Y_6$ sont 
scindés (proposition~\ref{prop:dec1}), de même que $I_4X_3$
(\S \ref{subs:typalg1-4}, cas $I_1X_3$, p.~\pageref{para:R=I_1X_3})
et $I_2X_5$ (\S \ref{subs:typalg5}, cas $X_5$,
p.~\pageref{para:R=X5}).
Il reste à traiter $I_3R_2R_2$, $I_{1,2}R_2R_2$, $I_1X_3X_3$,
$R_2R_2X_3$, $J_2J_2X_3$ et $Z_7$.

\paragraph{Cas où $R=I_3 \oplus R_2 \oplus R_2$ (ordre $2$).} 
Reprenons les notations du \S \ref{subs:typalg5} 
p.~\pageref{para:R=I1R2R2} qui traite partiellement  ce cas.
La composante symétrique $A_F$ de $F$ est donnée par l'équation
(\ref{equa:AF_IqR2R2})), dans laquelle $\det(A)$ est impair.
Comme $A$ est une matrice symétrique $3 \times 3$, cela entraîne que
l'un des termes diagonaux de $A$ est impair. La forme $A_F$
est donc de type impair et appartient au genre 
${\rm I}_{r,s}(2^2_{\rm II})$ (la décomposition de Jordan
2-adique de $A_F$ est précisée au \S \ref{subs:typalg5}).
Si $A_F$ est définie, d'après \cite{Kork-Zolo1877a}
($\ga_5=2^{3/5}$) elle représente $\pm 1$ sauf si elle est
équivalente à la forme paire $D_5$, ce qui est exclu. Dans le cas
indéfini, le genre ${\rm I}_{r,s}(2^2_{\rm II})$ ($r+s=5$)
ne contient qu'une classe entière représentée par
$I_{r-1,s-1}\oplus 2 U_2$. Dans tout les cas $A_F$ représente
$\pm 1$ et on peut supposer (\S \ref{subs:typalg5} {\it ibid.})
que $A_F= (\pm I_1) \oplus B$~: les types associés à $R$ sont
tous scindés. Noter l'existence d'une relation, à savoir 
\begin{equation}
\label{}
[I_3 \oplus F_4] = [(-I_1) \oplus F_6],
\end{equation}
qui provient de l'équivalence entre les formes $I_3 \oplus U_2$
et $(-I_1) \oplus D_4$.

\paragraph{Cas où $R=I_{1,2} \oplus R_2 \oplus R_2$ (ordre $2$).} 
D'après l'étude faite au \S \ref{subs:typalg6} (cas
$I_2^-\oplus R_2 \oplus R_2$, p.~\pageref{para:R=I2^-R2R2})
la composante symétrique $A_F$ de $F$ est de la forme
(\ref{equa:AFBF}) avec~$A$ scalaire impair. Ainsi, on a 
$A_F \in {\rm I}_{r,s}(2_{\rm II}^2)$, genre  forcément indéfini pour
$r+s=3$ (cf. \cite[p. 387]{Conw-Sloa1999a}). Il existe donc
un élémént $Q_1 \in \GL_3(\Z)$ tel que $Q_1 \cdot A_F = (\pm I_1)\oplus 2 U_2$
qui représente ${\rm I}_{r,s}(2_{\rm II}^2)$ ($r+s=3$). Nécessairement on~a 
$Q_1 \in \Ga_0^{3,2}$ et $Q_1$ se complète en un élément de $\Ga^{3,4}$
(lemme \ref{lemm:Ga_0^kl}), groupe isomorphe au commutant $C_R$.
Finalement on peut supposer $A_F = (\pm I_1)\oplus 2 U_2$ et
tous les types associés à $R$ sont scindés.
 
\paragraph{Cas où $R=I_1 \oplus X_3 \oplus X_3$ (ordre $4$).} 
Soit $\De=(I_2\oplus R_3 \oplus I_2)(I_1 \oplus \De_3 \oplus \De_3)$
(voir équation~(\ref{equa:De3_De4})).  On a
$\De R \De^{-1} =I_3 \oplus J_2 \oplus J_2$ ; les éléments
de $C_R^\De$ sont de la forme $Q_1 \oplus Q_2$ où $Q_1 \in \Ga_0^{3,2}$,
$Q_2$ s'identifie à un élément $\mfP \in \GL_2(\Z[i])$ et 
le bloc $2 \times 2$ inférieur droit de $Q_1$ est congru à
$\re \mfP + \im \mfP$ modulo~2. Soit $A_F$ la composante
symétrique de $F$. Comme au \S \ref{subs:typalg5}  (cas
$I_1\oplus R_2 \oplus R_2$), on prouve que $A_F$ appartient à l'un
des genres ${\rm I}_{r,s}(2_{\rm I}^2)$ ou ${\rm I}_{r,s}(2_{\rm II}^2)$,
$r+s=3$. En utilisant la constante d'Hermite $\ga_3$ et 
\cite[pp. 386-387]{Conw-Sloa1999a} on voit que  $A_F$ représente $\pm 1$
dans tous les cas. Par suite il existe $Q_1 \in \Ga_0^{3,2}$
tel que $Q_1 \cdot A_F =(\pm I_1) \oplus 2 B$ ($\det(B) = \pm 1$).
D'après le  lemme \ref{lemm:GL2(Z[i])}, $Q_1$ se complète en un élément 
de $C_R^\De$ et tous les types associés à $R$ sont scindés.

\paragraph{Cas où $R=R_2 \oplus R_2 \oplus X_3$ (ordre $4$).} 
On prend $\De=(I_2\oplus R_3 \oplus I_2)(\De_4 \oplus \De_3)$
(équation (\ref{equa:De3_De4})), matrice  qui vérifie $\De R \De^{-1}=
I_{3,2} \oplus J_2$ et $\det(\De)=8$. Posons 
$\De \cdot F = A_F \oplus B_F \oplus C_F$. On a alors
\begin{equation}
\label{equa:F_AF}
{\small
\begin{array}{l}
F=\left(
\begin{array}{cccc}
a(I_2+R_2) &  bI_2+cR_2 & e \al' & 0 \\
cI_2+bR_2 &  d(I_2+R_2) & f \al' & 0 \\
e \al     & f \al	& g	& h \be \\
0 & 0 & h J_2 \be' & h(I_2+J_2)
\end{array}
\right)
~~ {\rm et}\\ \\
A_F= 2 \left(
\begin{array}{ccc}
2a & b+c   &  2e  \\
b+c &  2d & 2f \\
2e & 2f & 2g-h
\end{array}
\right),
\end{array}
}
\end{equation}
avec $\al =(1,1)$, $\be=(1,0)$ et $a,b,c,d,e,f,g,h \in \Z$.
De plus $B_F= 2 (c-b)J_2$ et $C_F= h(I_2 + J_2)$. On voit
ainsi que $\det(A_F/2)=\pm 1$, $b-c=\pm 1$ et $h = \pm 1$.
Un calcul modulo~2 montre que toute équivalence entre $A_F$
et $A^0$ de la forme (\ref{equa:F_AF}) s'écrit
$Q_1=\smallmat{U}{V}{W}{x} \in GL_3(\Z)$ avec $V$ et $W$
divisibles par 2. Il existe alors $Q_2 \in \GL_2(\Z)$
congrue à $U$ modulo~2 et dans ces conditions, on vérifie
que  $Q_1 \oplus Q_2 \oplus I_2 \in C_R^\De$. On pourra
donc fixer arbitrairement $A_F$, pourvu qu'elle soit de la forme
(\ref{equa:F_AF}). Si $A_F$ et indéfinie, on prend $A^0/2=
U_2 \oplus (\pm I_1)$ qui conduit à 4 types scindés
$\pm (F_4 \oplus G_3)$, $\pm (F_4 \oplus H_3)$. Si $A_F$ est définie
(disons positive quitte à changer le signe de $F$), posons
$$
{\small
A^0=2 \left(
\begin{array}{ccc}
2 & 1 & 0\\
1 & 2 & 2\\
0 & 2 & 3
\end{array}
\right)
\esp {\rm et} \esp 
\Phi=\left(
\begin{array}{ccc}
0 & 1 & 0\\
1 & 0 & 0\\
2 & -2 & 1
\end{array}
\right).}
$$
Remarquons que $\Phi$ est une isométrie de $A^0$ telle que
$\Phi \oplus U_2 \oplus I_2 \in C_R^\De$ : on peut fixer 
$b-c$ sans changer $A^0$.  En fixant $A_F=A^0$ (ci-dessus) et
$b-c=-1$, on trouve deux matrices $F_7$ et $G_7$ selon que
$h$ vaut $1$ ou $-1$ (voir équation~(\ref{equa:indec_rg7}),
annexe). Le cas où $A_F$ est définie
conduit donc à 4 types indécomposables $\pm F_7$ et $\pm G_7$.

\begin{table}[!tbp]
\centering
\renewcommand{\arraystretch}{1.1}
\renewcommand{\tabcolsep}{0.6 em}%
{\small
\begin{tabular}{llll}
\hline
$d$  & $F F^\vee$ & $F$ &  type réel\\
\hline
1   & $I_7$ & $I_{p,7-p}$  & $I_{p,7-p}$\\
\hline
2   & $I_{5,2}$ & $I_{p,5-p} J_2$ & $I_{p,5-p} J_2$ \\
   & $I_{3,4}$ & $I_{p,3-p} J_4$  & $I_{p,3-p} J_2^2$ \\
   & $I_{1,6}$ & $I_{p,1-p} J_6$  & $I_{p,1-p} J_2^3$ \\
   & $I_3 R_2 R_2$  & $I_{p,3-p} F_4$ & $I_{p+1,4-p}J_2$\\
 & &$\pm I_1 F_6$ & $\pm I_5 J_2$ \\
   & $I_{1,2} R_2 R_2$ & $\pm I_1 J_2 F_4 $ & $\pm I_{2,1}J_2^2$ \\
\hline
4   & $I_3 J_2 J_2$ & $I_{p,3-p} K_4$  & $I_{p,3-p}P_{2;4,1}P_{2;4,1}^- $\\
 & & $I_{p,3-p} L_4$  & $I_{p,3-p}P_{2;4,1}^2$\\
 & & $I_{p,3-p} L_4^-$  & $I_{p,3-p}P_{2;4,1}^{-2}$\\
 & $I_{1,2} J_2 J_2$ & $\pm I_1 J_2 K_4$ & $\pm I_1J_2P_{2;4,1}P_{2;4,1}^-$\\
 & & $\pm I_1 J_2 L_4$   & $\pm I_1J_2P_{2;4,1}^2$\\
 & & $\pm I_1 J_2 L_4^-$   &$\pm I_1J_2P_{2;4,1}^{-2}$\\
   & $I_4 X_3$ & $I_{p,4-p} G_3$ & $I_{p+1,4-p}P_{2;4,1}^-$\\
 & &$I_{p,4-p} G_3^-$ & $I_{p,5-p}P_{2;4,1}$\\
 & &$\pm I_4 H_3$  &  $\pm I_5 P_{2;4,1}$ \\
 & &$\pm U_4 G_3$  &  $\pm I_{3,2} P_{2;4,1}^-$\\
 & &$\pm U_4 H_3$  &  $\pm I_{3,2} P_{2;4,1}$\\
   & $I_{2,2} X_3$ & $\pm I_2 J_2 G_3$ & $\pm I_3 J_2 P_{2;4,1}^-$\\
 & & $\pm I_2 J_2 G_3^-$ & $\pm I_{2,1} J_2 P_{2;4,1}$\\
 & & $\pm I_2 J_2 H_3$ &  $\pm I_3 J_2 P_{2;4,1}$ \\
 & & $\pm I_2 J_2 H_3^-$  & $\pm I_{2,1} J_2 P_{2;4,1}^-$\\
 & & $\pm U_2 J_2 G_3$ & $\pm I_{2,1} J_2 P_{2;4,1}^-$\\
 & & $\pm U_2 J_2 H_3$ & $\pm I_{2,1} J_2 P_{2;4,1}$\\
   & $I_4^- X_3$ & $\pm J_4 G_3$ & $\pm I_1J_4P_{2;4,1}^-$\\
 & & $\pm J_4 H_3$ & $\pm I_1J_4P_{2;4,1}$\\
   & $I_1 X_3 X_3$ &$\pm I_1 G_3 G_3$ & $\pm I_3P_{2;4,1}^{-2}$\\
 & &$\pm I_1 G_3 G_3^-$  & $\pm I_{2,1}P_{2;4,1}P_{2;4,1}^-$\\
 & &$\pm I_1 G_3^- G_3^-$ & $\pm I_{1,2}P_{2;4,1}^2$ \\
 & &$\pm I_1 G_3 H_3$  & $\pm I_3P_{2;4,1}P_{2;4,1}^-$\\
 & &$\pm I_1 G_3^- H_3$  & $\pm I_{2,1}P_{2;4,1}^2$\\
 & &$\pm I_1 H_3 H_3$ &  $\pm I_3P_{2;4,1}^2$\\
 & &$\pm I_1 G_6$   & $\pm I_{2,1}P_{2;4,1}P_{2;4,1}^-$\\
   & $I_1 X_3^- X_3^-$ &$\pm I_1 H_6$ & $\pm I_1J_2P_{2;4,1}P_{2;4,1}^-$\\
 & $R_2 R_2 X_3$ & $\pm F_4 G_3$ & $\pm I_{2,1}J_2P_{2;4,1}^-$\\
 & & $\pm F_4 H_3$ & $\pm I_{2,1}J_2P_{2;4,1}$\\
 & & $\pm F_7^*$  & $\pm I_3 J_2 P_{2;4,1}$ \\
 & & $\pm G_7^*$  &$\pm I_3 J_2 P_{2;4,1}^-$ \\
 & $J_2 J_2 X_3$ & $\pm K_4 G_3$ & $\pm I_1P_{2;4,1}P_{2;4,1}^{-2}$\\
 &  & $\pm K_4 H_3$ &  $\pm I_1P_{2;4,1}^2P_{2;4,1}^-$\\
 &  & $\pm L_4^- G_3$ & $\pm I_1P_{2;4,1}^{-3}$\\
 &  & $\pm L_4 H_3$ &  $\pm I_1P_{2;4,1}^3$\\
\hline
\end{tabular}}
\caption{Types isoduaux principaux de rang $7$ ($t'=166$) : ordre 1,2,4}
\label{tabl:types_7_a}
\end{table}

\paragraph{Cas où $R=J_2 \oplus J_2 \oplus X_3$ (ordre $4$).} 
En posant  $\De=I_4 \oplus R_3^2 \De_3$  et $\Psi=I_4 \oplus \Psi_3$
(équations (\ref{equa:Psi3}) et (\ref{equa:De3_De4})), on a  
$\De R \De^{-1}= J_2 \oplus J_2\oplus J_2\oplus I_1$  et
$\Psi R = R^\vee \Psi$. Soit $\mcO=\Z[i]$. Si $P \in C_R$, alors
$\De P \De^{-1}= \mfP \oplus (\pm I_1)$ avec $\mfP \in \GL_3(\mcO)$ ;
de plus, les coefficients $\mfP_{3,1}$ et  $\mfP_{3,2}$ sont 
multiples  de  $1+i$. Comme d'habitude nous écrivons  $P=P(\mfP,\pm 1)$.
Les solutions inversibles de $F=RF'$ sont données par
$F\Psi=P(\mfF,\det F)$ avec
$$
{\small
\mfF = \left(
\begin{array}{cc}
\mfA & z \\
(1+i)z^* & i f 
\end{array}
\right)
\esp {\rm o\grave{u}} \esp 
\mfA=\left(
\begin{array}{cc}
a(1-i) & \be \\
-i\bar{\be} & d (1-i) 
\end{array}
\right),}
$$
$a,d$ et $f \in \Z$, $\be \in \mcO$ et $z\in \mcO^2$. En développant 
le déterminant de $\mfF$ qui vaut $\pm 1$, on remarque que 
$f |\be|^2$ est impair. On vérifie également que $z=0$ correspond
à un type scindé suivant  $J_2 \oplus J_2$ et $X_3$.
L'action de $P=P(\mfP,\pm 1) \in C_R$ sur $\mfF$
se traduit par $PFP' \Psi =P(\mfP \mfF \mfP^\tau,\det F)$, où
$\tau$ est l'anti-automorphisme définit comme suit :
$$
{\small
\mfP^\tau=\left(
\begin{array}{cc}
\mfQ^* & i w  \\
-(1+i)v^* & \bar{u} 
\end{array}
\right)
\esp {\rm pour} \esp
\mfP = \left(
\begin{array}{cc}
\mfQ & v \\
(1+i)w' & u 
\end{array}
\right),}
$$
avec $u \in \mcO$, $v$, $w \in \mcO^2$ et $\mfQ \in M_2(\mcO)$. Noter
que $\tau$ n'est pas involutif comme dans les autres cas
-- en fait $\tau$ est d'ordre~8. Remarquons d'abord que l'élément 
$\mfP=\smallmat{I_2}{~ifz}{0}{~1}$ transforme $z$ en $z(1-f^2)$~:
si $f=\pm 1$, le type $[F]$ est scindé. Nous allons montrer que
l'on peut toujours se ramener à ce cas. Pour cela, considérons 
$-i(\mfP\mfF\mfP^\tau)_{3,3}$ qui définit à $\mfF$ fixée
une forme hermitienne $h_F$ en les  variables $u \in \mcO$ et
$w=\smallvect{w_1}{w_2} \in \mcO^2$, à savoir
$$
h_F(u,w_1,w_2)=f |u|^2 + 2a |w_1|^2 + 2d |w_2|^2
+2 \re [(1+i)(uz^*\bar{w} + \be w_1 \bar{w_2})].
$$
On a $\det(h_F)=2 \det(\mfF)=\pm 2$. Par ailleurs, la relation
$h_F(u,w_1,w_2)=\pm 1$ entraîne $\pgcd(u,(1+i)w_1,(1+i)w_2)=1$
et dans ces conditions, comme $\mcO$ est principal,
 il existe $\mfP \in \GL_3(\mcO)$ dont
la troisième ligne est $(u,(1+i)w_1,(1+i)w_2)$. Pour conclure,
il suffit donc de prouver que $h_F$ représente $1$ ou $-1$.  
Une base du $\Z$-module $\mcO^3$ étant choisie, notons $H_F^\Z$
la matrice de la forme quadratique $h_F^\Z$ induite par $h_F$ 
sur ce $\Z$-module.
Sachant que $f|\be|^2$ est impair, on trouve une sous-matrice extraite
de $H_F^\Z$ de façon diagonale dont le déterminant est impair.
D'après la relation (\ref{equa:dec_orth}), $H_F^\Z$ est équivalente
sur $\Z_2$ à $A \oplus 2 B$ avec $A$ inversible de rang 4 et $B$
inversible  paire. Par suite $h_F^\Z$ appartient au genre
${\rm I}_{r,s}(2^2_{\rm II})$. Dans le cas indéfini, ce genre
comprend une seule  classe d'équivalence ; la forme 
$h_F^\Z$ est équivalente à $I_{p,q}\oplus 2 U_2$ avec
$(p,q)=(1,3)$ ou $(3,1)$ (puisque  $\det(h_F^\Z)=\det(h_F)^2=4$)
et $h_F$ représente $\pm 1$.

Supposons maintenant que $h_F$ soit définie, par exemple positive.
Nous estimons le minimum de $h_F$ grâce à un ensemble de Siegel
pour $\GL_n(\C)$. Soit $K=U(n)$, soit $A$ le sous-groupe des matrices
diagonales à coefficients réels $>0$ et soit $N$ le sous-groupe des
matrices triangulaires supérieures  unipotentes. On a
la décomposition d'Iwasawa $\GL_n(\C)=KAN$. Si $u$ et $t$ sont des
réels $>0$, on note $A_t$ l'ensemble
des $(a_{j,k}) \in A$ telles que $a_{j,j} \leq t a_{j+1,j+1}$
($1 \leq j \leq n-1)$, $N_u$ l'ensemble des $(n_{j,k}) \in N$ 
telles que $|n_{j,k}| \leq u$  ($1 \leq j < k \leq n)$ et on 
pose $\mfS_{t,u}=KA_tN_u$ (ensemble de Siegel). Grâce à 
l'approximation des complexes par les éléments de $\mcO$, on  prouve
comme pour $\GL_n(\R)$ que $\GL_n(\C)= \mfS_{\sqrt{2},1/\sqrt{2}}
\GL_n(\mcO)$  (cf. \cite[pp. 14-15]{Borel1969a}). En utilisant 
l'action transitive (ici à droite) de $\GL_n(\C)$  sur les
matrice hermitiennes $>0$, on en déduit qu'il existe une base de
$\mcO^3$ dans laquelle les coefficients  diagonaux de la matrice de
$h_F$  sont de la forme
$$
\left\{
\begin{array}{l}
h_{1,1}= b_1  \\
h_{2,2}= b_2  + b_1 |n_{1,2}|^2 \\
h_{3,3}= b_3  + b_2 |n_{2,3}|^2 + b_1 |n_{1,3}|^2
\end{array}
\right.
$$
avec $0<b_1 \leq 2 b_2 \leq 4 b_3$, $|n_{j,k}|^2 \leq 1/2$
($j < k $) et $b_1b_2b_3 = \det(h_F)=2$. L'inégalité 
$b_1^3 \leq 16$ montre que $b_1$ vaut 1 ou 2. Si le 
minimum de $h_F$ était égal à~2, on aurait $b_1=2$,
$b_2b_3=1$, $b_3 \leq 1\leq b_2 \leq \sqrt{2}$, puis
$2\leq  h_{2,2} \leq \sqrt{2}+1$,  $2 \leq  h_{3,3} \leq 2 +\sqrt{2}/2$
et enfin  $h_{1,1}=h_{2,2}=h_{3,3}=2$, ce qui est absurde puisque $h_F$
est impaire. Donc $h_F$ représente 1 (en fait on pourrait montrer que
$H_F^\Z$ est équivalente à $I_2 \oplus D_4$).

Finalement, les types associés à $R$ sont tous scindés. On obtient
seulement 8 types (table~\ref{tabl:types_7_a}) car il y a deux relations :
$$
[K_4 \oplus G_3] = [(-L_4) \oplus H_3]
\esp {\rm et} \esp
[K_4 \oplus H_3] = [L_4 \oplus G_3].
$$
Notons $F_1,F_2,F_3,F_4$ ces 4 types et posons  $F_j \Psi =
P(\mfF_j,1)$ ($j=1,\ldots,4$). Alors $\mfF_1=\mfA \oplus (-i)I_1$, 
$\mfF_2=\mfB \oplus iI_1$,  $\mfF_3= \mfA \oplus i I_1$,
$\mfF_4=-\mfF_2$ où $\mfA=\smallmat{0}{~1}{-i}{~0}$ et
$\mfB=\smallmat{i-1}{~i}{-1}{~i-1}$. Les relations précédentes
proviennent de $\mfF_2=\mfP \mfF_1 \mfP^\tau$ et
$\mfF_4=\mfQ \mfF_3 \mfQ^\tau$ avec
$$
{\small
\mfP=\left(
\begin{array}{ccc}
i & 0 & i\\
0 & -i & -1\\
1+i & 1+i & 1
\end{array}
\right)
\esp {\rm et} \esp 
\mfQ=\left(
\begin{array}{ccc}
1 & 0 & 1\\
0 & 1 & i\\
-1-i  & 1+i & -1
\end{array}
\right).}
$$

\begin{table}[h]
\centering
\renewcommand{\arraystretch}{1.1}
\renewcommand{\tabcolsep}{0.6 em}%
{\small
\begin{tabular}{llll}
\hline
$d$  & $F F^\vee$ & $F$  & type réel\\
\hline
8   & $I_2 X_5$ & $I_{p,2-p} F_5$  & $I_{p+1,2-p} P_{2;8,1}^- P_{2;8,3}$\\
 & & $I_{p,2-p} F_5^-$ & $I_{p,3-p} P_{2;8,1} P_{2;8,3}^-$\\
 & & $I_{p,2-p} H_5$  & $I_{p+1,2-p} P_{2;8,1}^- P_{2;8,3}^-$\\
 & & $I_{p,2-p} H_5^-$ & $I_{p,3-p} P_{2;8,1} P_{2;8,3}$ \\
 & &$\pm I_2   G_5$  & $\pm I_3 P_{2;8,1} P_{2;8,3}^-$\\
 & &$\pm I_2   K_5$  & $\pm I_3 P_{2;8,1} P_{2;8,3}$\\
 & &$\pm U_2   F_5$   & $\pm I_{2,1} P_{2;8,1}^- P_{2;8,3}$\\
 & &$\pm U_2   G_5$   & $\pm I_{2,1} P_{2;8,1} P_{2;8,3}^-$\\
 & &$\pm U_2   H_5$   & $\pm I_{2,1} P_{2;8,1}^- P_{2;8,3}^-$\\
 & &$\pm U_2   K_5$   & $\pm I_{2,1} P_{2;8,1} P_{2;8,3}$\\
 & $I_2^- X_5$ &$\pm J_2 F_5$ & $\pm J_2 I_1 P_{2;8,1}^- P_{2;8,3}$\\
 & &$\pm J_2   G_5$  & $\pm J_2 I_1 P_{2;8,1} P_{2;8,3}^-$\\
 & &$\pm J_2   H_5$  & $\pm J_2 I_1 P_{2;8,1}^- P_{2;8,3}^-$\\
 & &$\pm J_2   K_5$  & $\pm J_2 I_1 P_{2;8,1} P_{2;8,3}$\\
 & $I_1 Y_6$ &$\pm I_1 K_6$ & $\pm I_1 P_{2;4,1} P_{2;8,1}^- P_{2;8,3}^-$\\
 & &$\pm I_1 K_6^-$ &  $\pm I_1P_{2;4,1}^-P_{2;8,1}P_{2;8,3}$\\
 & &$\pm I_1 L_6$   & $\pm I_1P_{2;4,1}P_{2;8,1}P_{2;8,3}$\\
 & &$\pm I_1 L_6^-$ & $\pm I_1P_{2;4,1}^-P_{2;8,1}^-P_{2;8,3}^-$\\
 & &$\pm I_1 M_6$ &  $\pm I_1P_{2;4,1}P_{2;8,1}P_{2;8,3}^-$\\
 & &$\pm I_1 M_6^-$  & $\pm I_1P_{2;4,1}^-P_{2;8,1}^-P_{2;8,3}$\\
 & &$\pm I_1 N_6$   & $\pm I_1P_{2;4,1}P_{2;8,1}^-P_{2;8,3}$\\
 & &$\pm I_1 N_6^-$ & $\pm I_1P_{2;4,1}^-P_{2;8,1}P_{2;8,3}^-$\\
  & $Z_7$ &$\pm H_7^*$  & $\pm I_1 P_{2;4,1}^- P_{2;8,1}^- P_{2;8,3}$\\ 
 & &$\pm K_7^*$   & $\pm I_1 P_{2;4,1}^- P_{2;8,1} P_{2;8,3}^-$\\
 & &$\pm L_7^*$   & $\pm I_1 P_{2;4,1} P_{2;8,1}^- P_{2;8,3}^-$\\
 & &$\pm M_7^*$   & $\pm I_1 P_{2;4,1} P_{2;8,1} P_{2;8,3}$\\
\hline
\end{tabular}}
\caption{Types isoduaux principaux de rang $7$ ($t'=166$),
suite : ordre 8}
\label{tabl:types_7_b}
\end{table}

\paragraph{Cas où $R=Z_7$ (ordre $8$).} 
Soit $\eps$ une racine primitive 8-ième de l'unité, soit
$\mcO=\Z[\eps]$ et pour $\tht \in \mcO$, soit $M(\tht)$
la matrice de la multiplication par $\tht$ dans la base
$(1,\eps,\eps^2,\eps^3)$. Posons
$$
{\small 
\De=\left(
\begin{array}{rrrrrrr}
4 &  -2 &  ~2 &  ~1 &  ~1 &  1 &  1 \\
0 &  2 &  0 &  0 &  1 &  0 &  -1 \\
0 &  0 &  2 &  1 &  0 &  -1 &  0 \\
0 &  0 &  0 &  1 &  0 &  0 &  0 \\
0 &  0 &  0 &  0 &  1 &  0 &  0 \\
0 &  0 &  0 &  0 &  0 &  1 &  0 \\
0 &  0 &  0 &  0 &  0 &  0 &  1
\end{array}
\right)
\esp {\rm et} \esp
\Psi=\left(
\begin{array}{rrrrrrr}
4 &  -2 &  2 &  1 &  1 &  ~1 &  ~1 \\
-2 &  0 &  -2 &  -1 &  -1 &  0 &  0 \\
2 &  0 &  0 &  0 &  1 &  1 &  0 \\
1 &  0 &  0 &  0 &  0 &  1 &  0 \\
1 &  -1 &  0 &  0 &  0 &  0 &  1 \\
1 &  -1 &  1 &  0 &  0 &  0 &  0 \\
1 &  0 &  1 &  1 &  0 &  0 &  0
\end{array}
\right).}
$$
On a $\Psi R =R^\vee \Psi$, $\det(\De)=16$ et 
$\De R \De^{-1}=I_1 \oplus J_2 \oplus X_4$. Si $P\in C_R$, alors
$\De P  \De^{-1}=uI_1  \oplus \smallmat{v}{-w}{w}{v} \oplus M(\tht)$
avec  $u=\pm 1$, $\al=v+iw \in \Z[i]^\times$ et  
$\tht =x +y \eps + z \eps^2 +t \eps^3  \in \mcO^\times$ ; de plus
les relations  $v \equiv x+z$ (mod 2), $w \equiv y +t$ (mod 2) et
$u \equiv v+w+2(y+z)$ (mod 4) caractérisent $C_R^\Delta$ comme 
sous-groupe de $\Z^\times \times \Z[i]^\times \times \mcO^\times$.
On voit que la projection de $C_R^\Delta$ sur la composante 
$\mcO^\times$ est surjective. Posons $P=P(u,\al,\tht)$.
L'involution $\sig$ s'écrit  
$P(u,\al,\tht)^\sig=P(u,\bar{\al},\bar{\tht})$. 
Sachant que $\Psi^\vee \Psi = R = P(1,-i,\eps)$, la relation 
(\ref{equa:G=RxiGsig}) page~\pageref{equa:G=RxiGsig} montre que 
$F\Psi = P(a,bi,\mff)$ avec $a=\det(F)=\pm 1$, $b^2= 1$  et
$\mff \in \mcO^\times$ tel que $\mff = \eps^2 \bar{\mff}$. En utilisant
l'unité fondamentale $\eta = 1-\eps+\eps^3$ et l'action de $C_R$,
on se ramène à  $\mff= \pm \eps$ ou $\pm \eps \eta^{-1} =
\pm(1-\eps +\eps^2)$,
ce qui conduit à 8 types indécomposables : $H_7=P(1,-i,\eps)$, 
$K_7=P(1,-i,-\eps)$, $L_7=P(1,i,\eps \eta^{-1})$, 
$M_7=P(1,i,-\eps \eta^{-1})$ et leurs opposés
(équation (\ref{equa:indec_rg7}), annexe).

\section{Classification des types géométriques $(n \leq 7)$}
\label{sect:class_geom}

\subsection{Types réels. Aspects différentiels et métriques}
\label{subs:types_reels_aspects_diff}
Pour déterminer les types géométriques, \cadab les matrices
de Gram associées à un type algébrique, il sera commode 
d'utiliser des décompositions {\em réelles} des types isoduaux.
Par extension du cas entier, nous introduisons la notion de
{\em type réel} -- on définirait de même la 
notion de type rationnel. 
Soit $\GL_n^{\pm}(\R)$ le sous-groupe
de $\GL_n(\R)$ formé des matrices de déterminant $\pm 1$. 

\begin{defi}\label{defi:type_reel}
\indent (1) Un élément $F\in\GL_n^\pm(\R)$ {\em représente un type 
algébrique réel} si $FF^\vee$ est d'ordre fini. La classe 
d'équivalence de $F$ sous l'action de $\GL_n^\pm(\R)$ est alors 
appelée {\em type algébrique réel}, noté $[F]$.\\
\indent (2) Le {\em type géométrique réel} associé à un type 
algébrique réel $[F]$ est  la classe d'isométrie 
de $V_F=\{A\in P_n; AF^\vee A=F \}$ dans l'espace $P_n$ des matrices de 
Gram.
\end{defi}

Si $F \in \GL_n(\Z)$ représente un type isodual,  
les propriétés différentielles  et géométriques
(dimension, courbure, etc) de $V_F$ ne dépendent que du 
type réel associé à $F$. Nous examinons ici les premières
propriétés de $V_F$ en considérant plus généralement $F$ réel. 
La structure géométrique de $V_F$  sera précisée
au \S \ref{subs:composantes}, théorème~\ref{theo:geom_VF}.
\par
Considérons maintenant l'ensemble  $P_n$ des matrices symétriques  
positives de déterminant~1. Rappelons que la métrique
\begin{equation}
\label{equa:met_Pn}
ds^2=\tr (A^{-1}dA)^2\esp (A \in P_n).
\end{equation}
munit $P_n$ d'une structure d'espace symétrique riemannien. 
De plus, le groupe  $\GL_n^\pm(\R)$ agit isométriquement sur $P_n$ par  
\begin{equation}\label{equa:PA} 
P \cdot A = PAP'  \esp (P \in GL_n^\pm(\R), A \in P_n).
\end{equation}
Nous appellerons {\em sous-espace symétrique de de $P_n$} toute
sous-variété totalement géodésique, complète et connexe de $P_n$.

\begin{prop}
\label{prop:som_prod_se_sym}
Soient $V$ et $W$ des sous-espaces symétriques de $P_m$ et $P_n$
respectivement. Alors \\
\indent {\rm a)} $V\oplus W =\{A\oplus B; A \in V, B\in W\}$ est 
un sous-espace symétrique de $P_{n+m}$ isométrique au produit riemannien
$V\times W$.\\
\indent {\rm b)}  $V\otimes W =\{A\otimes B; A \in V, B\in W\}$ est 
un sous-espace symétrique de $P_{mn}$ isométrique à la variété
 produit $V\times W$ munie de la métrique
$$ds^2_{V\otimes W}=n \, ds^2_V + m \,  ds^2_W ,$$
où $ds^2_V$ et $ds^2_W$ désignent les métriques de $V$ et de $W$.
\end{prop}

\pr On sait qu'une partie non vide de $P_N$ 
est un sous-espace symétrique \ssi elle est fermée, connexe
et %
stable par toutes les symétries centrées en ses points 
(voir  \cite[lemme~2.10]{Bavard2005a}).
La symétrie $s_A$ centrée en $A \in P_N$ a pour expression 
$s_A(B)=AB^{-1}A$ ($B\in P_N$).  
Le cas de $V\oplus W$ est immédiat (propriétés topologiques claires
et calcul sur les  blocs diagonaux). Ensuite, rappelons 
que $A \otimes B=(a_{ij}B)_{1\leq i,j\leq m}$ (matrice diagonale 
par blocs $n \times n$ , où $A=(a_{ij})_{1\leq i,j\leq m}$).
Il s'agit évidemment d'une version matricielle du produit tensoriel
 des formes bilinéaires. On  a les relations suivantes :
\begin{equation}
\label{equa:prod_tens}
\renewcommand{\arraystretch}{1.2}
\begin{array}{l}
 (A\otimes B)(C\otimes D)=AC \otimes BD, ~~
~~ I_m\otimes I_n=I_{mn}\\
(A\otimes B)'=A'\otimes B',~~
(A\otimes B)[u\otimes v] = A[u]B[v],\\
\tr(A\otimes B)=\tr A \tr B, ~~
\det(A\otimes B)=(\det A)^n (\det B)^m, 
\end{array}
\end{equation}
avec $A,C \in M_m(\R)$, $B,D \in M_n(\R)$ et $(u,v)\in \R^m\times \R^n$.
Par suite $V\otimes W$ est bien inclus dans $P_{mn}$ et 
\og symétrique\fg. Montrons que
$V\otimes W$ est fermé dans  $P_{mn}$. Soit
$A_k\otimes B_k = (a_{ij;k}B_k)_{ij}$ une suite qui converge vers 
$M=(M_{ij})_{ij} \in P_{mn}$ ($1\leq i,j\leq m$). En prenant le 
déterminant de chaque bloc, on trouve $\lim_k a_{ij;k}^n= \det M_{ij}$.
Puisque $M\in P_{mn}$ et $A_k\in P_m$, on a $\det M_{1,1} >0$ et 
$a_{1,1;k}>0$ ; par suite, $(a_{1,1;k})_k$ converge vers un réel $>0$
et $B_k$ converge vers une matrice $B\in W$ ($W$ est fermé dans 
$M_n(\R)$). Il en résulte que $A_k$ converge vers une matrice
$A\in V$ ($V$ est fermé dans $M_m(\R)$). Comme $V\otimes W$ est 
clairement connexe, on conclut que c'est un sous-espace symétrique.
Concernant l'aspect métrique, on a d'après \eqref{equa:prod_tens} 
$$(A\otimes B)^{-1}d(A\otimes B)= 
A^{-1}dA \otimes I_n + I_m \otimes B^{-1}dB,$$
d'où le résultat (on rappelle que $\tr(A^{-1}dA)=\tr(B^{-1}dB) = 0$).

\begin{prop}
\label{prop:VF_symetrique}
Soit $[F]$ un type algébrique réel (définition~\ref{defi:type_reel}). 
Alors l'ensemble
$V_F=\{A\in P_n;AF^\vee A=F\}$
  est non vide et c'est un sous-espace symétrique de $P_n$.
\end{prop}

\pr La non vacuité de $V_F$ résulte de l'argument donné en 
\ref{subs:carac_types} (théorème du point fixe d'Élie Cartan).
On a une application $\tau$ de $P_n$ dans un sous-espace 
symétrique $W_n$ de $P_{2n}$ (voir~\ref{subs:carac_types}
et~\ref{subs:finitude_types}) et $V_F$ correspond 
{\it via} $\tau$ au lieu fixe d'un groupe fini d'isométries de
$W_n$. De plus (voir~équation~(\ref{equa:Wn}), p.~\pageref{equa:Wn}), 
on a $\tr(\tau_A^{-1}d\tau_A)^2 = 2 \tr(A^{-1}dA)^2$ ($A\in P_n$). Le 
plongement $\tau$ est donc  isométrique (à un facteur près), d'où 
le résultat.  
\par 
Si $[F]$ est un type algébrique réel, on pose 
$$\OO(F)=
\{P\in\GL_n(\R); P'FP=F\}$$
 (sous-groupe de $\GL_n^{\pm}(\R)$) et 
$\SO_F(\R) =  \OO_F(\R) \cap\SL_n(\R)$.

\begin{prop}
\label{prop:VF_orbite}
Soit $[F]$ un type algébrique réel. Alors $V_F$ est une orbite 
du groupe $\SO^t_F(\R)$ (transposé du  groupe $\SO_F(\R)$).
\end{prop}

\pr Il est immédiat de vérifier que $V_F$ est stable 
par l'action de $\OO^t_F(\R)$. 
La transitivité de $\SO^t_F(\R)$ résulte du 
lemme suivant.
\begin{lemm}
\label{lemm:conj_On}
Soient $U_1,U_2\in\OO(n)$. On suppose qu'il existe $P\in\SL_n(\R)$
tel que $U_2=PU_1 P'$ Alors $U_1$ et $U_2$ sont conjuguées par
un élément de $\SO(n)$. 
\end{lemm}

En effet soit $A=PP'\in P_n$ ($P\in\SL_n(\R)$). Il est clair
que $A\in V_F$ \ssi la matrice $U_P=P^{-1}FP^\vee$ est orthogonale
(\cadab $U_P=U_P^\vee$). Fixons un point $A_0=P_0 P'_0 \in V_F$ et 
posons $U_0=U_{P_0}$, puis
$$\mcP_{P_0,F}=\{P\in\SL_n(\R); PU_0P'=F\} = 
\{P\in\SL_n(\R); PP_0^{-1}\in \SO^t_F(\R)\}.$$
L'ensemble $\{PP';P \in\mcP_{P_0,F}\}$ est l'orbite de $A_0$ sous 
l'action de  $\SO^t_F(\R)$, incluse dans  $V_F$ (d'après le
début de la preuve). Maintenant si  $A=PP' \in V_F$, 
on a $U_P=(P^{-1}P_0)U_0(P^{-1}P_0)'.$
D'après le lemme~\ref{lemm:conj_On} il existe $U\in \SO(n)$ tel que 
$U_P=UU_0U^{-1}$. Posons $Q=PU$ : alors $PP'=QQ'$ et 
$QU_0Q'=PU_PP'=F$, donc $Q\in\mcP_{P_0,F}$. Ce qui montre que
$V_F=\{PP';P \in\mcP_{P_0,F}\}$; c'est donc une  orbite de $\SO^t_F(\R)$.

{\it Preuve du lemme~\ref{lemm:conj_On}.}
Rappelons que deux éléments de $\OO(n)$ sont conjugués dans $\OO(n)$
\ssi ils ont le même polynôme caractéristique. En particulier 
deux matrices orthogonales conjuguées dans $\GL_n(\R)$ le sont 
aussi dans $\OO(n)$. 
Comme  $U_i\in\OO(n)$ ($i=1,2$), on a
$U_2^2=U_2U_2^\vee = PU_1 U_1^\vee P^{-1} = P U_1^2 P^{-1}$.
Par suite, il existe $U\in \OO(n)$ tel que $U_2^2=U U_1^2U^{-1}$.
Le polynôme caractéristique des $U_i^2$ s'écrit
$$\chi_{U_1^2} = \chi_{U_2^2} = 
(X-1)^a(X+1)^b\prod_{j=1}^k(X^2-2 \cos(\theta_j)X
+1)^{m_j},$$
avec $a+b+2\sum_{j=1}^k m_j=n$, les $\theta_j$ étant mutuellement 
distincts et  choisis dans  $]0,\pi[$. 
On a $U_1^2=(U^{-1}P)U_1^2(U^{-1}P)^{-1}$ et 
$U_1^2=(U^{-1}P)'U_1^2(U^{-1}P)^\vee$ ;
les sous-espaces 
$E_+=\ker(U_1^2-I)$, $E_-=\ker(U_1^2+I)$ et 
$E_j = \ker(U_1^4-2\cos(\theta_j)U_1^2+I))$ ($1\leq j \leq k$)
sont donc stables par $U^{-1}P$ et par $(U^{-1}P)'$.
Ils sont également stables
par $U_3=U^{-1}U_2U$, qui vérifie $U_3=(U^{-1}P)U_1(U^{-1}P)'$. 
En remplaçant $U_2$ par $U_3$ et $P$ par $U^{-1}P$, on est donc
ramené au cas où  $\R^n=E_+$, $E_-$ ou $E_j$ (pour un 
$j=1,\ldots,k$). 

{\it Cas où $\R^n=\ker(U_1^2-I)$}. On a $U_i^2=I$, donc
$U'_i=U_i$ ($i=1,2$). Soient $p_i=\dim(\ker(U_i-I)$ et 
$q_i=\dim(\ker(U_i+I)$ ($i=1,2$). Le couple $(p_i,q_i)$ est la 
signature de $U_i$ comme forme bilinéaire symétrique. 
Puisque $U_2=PU_1P'$,
on a $(p_1,q_1)=(p_2,q_2)$ et les $U_i$ ($i=1,2$) sont conjuguées
dans $\OO(n)$. On peut choisir la conjugaison directe.

{\it Cas où $\R^n=\ker(U_1^2+I)$}. Les $U_i$ ($i=1,2$) sont
antisymétriques de  carré $-I$. Comme  $\chi_{U_1}=\chi_{U_2}=
(X^2+1)^{n/2}$, il existe un élément $V\in \OO(n)$ tel que
$U_2=VU_1V^{-1}=VU_1V'$. En comparant les Pfaffiens 
$$\pf (U_2) = \det(V)\pf(U_1)=\det(P)\pf(U_1),$$
on voit que $\det(V) = \det(P)$.

{\it Cas où $\R^n=\ker(U_1^4-2\cos(\theta)U_1^2+I)$ ($\theta\in
]0,\pi[$)}. Le polynôme caractéristique de $U_i$ est de la forme
$$\chi_{U_i}= (X^2-2\cos(\theta/2) X +1)^{p_i} 
(X^2 + 2\cos(\theta/2) X +1)^{q_i}\esp (i=1,2).$$
Le couple $(2p_i,2q_i)$ s'interprète comme la signature de la 
forme bilinéaire symétrique $U_i+U'_i$. Puisque $U_2+U'_2 = P(U_1+U'_1)P'$,
on a $(p_1,q_1)=(p_2,q_2)$ et il existe $V\in \OO(n)$ tel que
$U_2=VU_1V^{-1}=VU_1V'$. Les matrices antisymétriques 
$U_i- U'_i$ ($i=1,2$) sont de rang maximal ($\theta/2\in ]0,\pi/2[$)
et vérifient $U_2-U'_2=V(U_1-U'_1)V'=P(U_1-U'_1)P'$. En comparant les 
Pfaffiens (voir ci-dessus), on trouve $\det(V) = \det(P)$.

\begin{prop}[critère de Mahler pour les réseaux isoduaux]
\label{prop:VF_Mahler}
Soit $[F]$ un type algébrique {\em entier} et soit 
$\SO_F^t(\Z)=\OO_F^t(\R) \cap \SL_n(\Z)$. Alors l'application
$$\SO_F^t(\Z)\backslash V_F \to \SL_n(\Z)\backslash P_n$$
 est  propre et à fibres finies. En particulier on a un critère 
de Mahler dans l'espace  $\SO_F^t(\Z)\backslash V_F$ des 
réseaux isoduaux de type $[F]$.
\end{prop}

\pr C'est une application directe de la proposition~2.1 
de \cite{Bavard2005a} p.~222. Sachant que 
$V_F$ est une orbite de $\SO_F^t(\R)$ (proposition~\ref{prop:VF_orbite}),
il suffit de vérifier
que ce groupe est \og pseudo-algébrique\fg\ au sens de 
\cite[définition 2.1]{Bavard2005a}, ce qui est immédiat.
En effet soit  $\ph : \SL_n(\C) \to \GL(M_n(\C))$ la 
représentation usuelle ($\ph(P)\cdot X = PXP'$).  Le réseau 
$M_n(\Z)$ est stable par $\ph(\SL_n(\Z))$ et par définition
 $P\in \SO^t(\R)$ \ssi $\ph(P)$ fixe le vecteur $F\in M_n(\Z)$.
\medskip

Pour $P\in \GL_n^{\pm}(\R)$, notons $\Phi_P$ l'endomorphisme 
de $\sym_n(\R)$ (matrices  symétriques) défini par 
$$\Phi_P(X)= P \cdot X = P X P'  \esp (X \in \sym_n(\R)).$$
Soit $[F]$ un type algébrique réel et soit $A\in V_F$. 
L'espace tangent $T_AV_F$ est l'ensemble  des $X\in\sym_n(\R)$ 
 tels que $A F^\vee X + XF^\vee A  = 0$, c'est-à-dire
$$T_AV_F = \ker(\Phi_{AF^\vee}+\id).$$
On remarquera que, d'après la proposition~\ref{prop:VF_orbite},
les endomorphismes $\Phi_{AF^\vee}$ quand $A$ décrit $V_F$ sont tous
conjugués. En effet pour $A,B\in V_F$, il existe $P\in \SL_n(\R)$ tel
que $P\cdot F= F$ et $P\cdot A= B$ ;  on a alors 
$\Phi_{BF^\vee} = \Phi_P  \Phi_{AF^\vee} \Phi_P^{-1}$.
Noter également que $\Phi_{AF^\vee}$ coïncide sur $V_F$ avec la symétrie
géodésique centrée au point $A$ (laquelle est donnée par 
$\sig_A(B) = AB^{-1}A$ pour $B\in P_n$). Ainsi, les symétries de l'espace symétrique
$V_F$ sont données par des applications linéaires. 
\par
Soit $u\in \R$ et soit  $l_u^F$ la restriction à $V_F$ de la fonction longueur 
$l_u$  définie sur $P_n$  par $l_u(A)=A[u]$ ($A\in P_n$).

\begin{prop}[Espace tangent. Gradients des longueurs]
\label{prop:esp_tangent}
Soit $[F]$ un type algébrique réel tel que $\dim V_F \geq 1$, soit $A\in V_F$
 et soit $u\in \R$. Alors\\
\indent {\rm a)} le polynôme minimal de $\psi_{A,F}$
se met sous la forme $(x+1)f(x)$ avec $f\in \R[x]$ tel que $f(-1)\neq 0$ et 
on a  
$$T_AV_F=\im f(\Phi_{AF^\vee}),$$
\indent {\rm b)}  le gradient (pour la métrique \eqref{equa:met_Pn}) de $l_u^F$ au point 
$A$  est donné  par 
$$\nabla_A l_u^F = \frac{1}{f(-1)}f(\Phi_{AF^\vee})(Auu'A).$$
\end{prop}

\pr Observons  d'une part  que $\Phi_{AF^\vee}$ est d'ordre fini 
(car $\Phi_{AF^\vee}^2= \Phi_{FF^\vee}$)  et 
d'autre part que $\Phi_{AF^\vee}$ est une isométrie pour le produit scalaire
défini par 
$\langle X,Y\rangle _A = \tr (A^{-1}XA^{-1}Y)$ sur $\sym_n(\R)$. 
Le polynôme minimal de $\Phi_{AF^\vee}$ est donc à racines simples, d'où 
l'assertion~a). Ensuite, on a la décomposition orthogonale
$\sym_n(\R)=T_AV_F \oplus \im(\Phi_{AF^\vee}+\id)$. En écrivant 
$f=(x+1)g +f(-1)$ où $g\in \R[x]$, on voit que la projection
orthogonale (pour $\langle \cdot,\cdot \rangle _A$) sur $T_AV_F$ 
vaut $f(\Phi_{AF^\vee})/f(-1)$,  d'où l'assertion~b).

\begin{rema}%
\label{rema:espace_tangent}
Le polynôme $f$ ne dépend que de $F$ puisque tous les $\Phi_{AF^\vee}$ ($A\in V_F$)
sont conjugués. De plus, on voit que les éléments donnés par la 
proposition~\ref{prop:esp_tangent} (paramétrage de l'espace 
tangent et gradients) se calculent sans avoir à expliciter la 
sous-variété $V_F$. 
\end{rema}

\subsection{Composantes des types géométriques}
\label{subs:composantes}

 Soit $\mcE$ l'ensemble des couples d'entiers
$(k,l) \in \N^{*2}$  tels que $\pgcd(k,l)=1$
et $1 \leq l < \max(3,k)/2$. Considérons les polynômes
\begin{equation}
\label{equa:polynomes_Psi}
\left\{
\begin{array}{l}
\Psi_{1,1}=x-1,~~\Psi_{2,1}=x+1 ~~\mathrm{et}\\
\Psi_{k,l}= x^2-2\cos(2 l \pi/k)x+1 ~~((k,l) \in \mcE,~k \geq 3).
\end{array}
\right.
\end{equation}

\begin{prop} 
\label{prop:dec_can_R}
{\rm (décomposition canonique d'un type algébrique  réel)}\\
Soit \mbox{$F \in \GL_n^\pm(\R)$} tel que $R= FF^\vee$ soit d'ordre
fini~$d$. Pour tout  $(k,l)\in \mcE$ avec~$k$ diviseur de~$d$ 
on pose $W_{k,l}(F)=\ker \Psi_{k,l}(R^\vee)$. Alors\\
\indent {\rm a)} la somme des $W_{k,l}(F)$ est directe et égale
à $\R^n$,\\
\indent {\rm b)} cette somme directe $\bigoplus W_{k,l}(F)$ est
$F$-orthogonale {\em bilatère},\\
\indent {\rm c)} les classes d'équivalence des $\R$-modules 
bili\-né\-ai\-res $(W_{k,l}(F),F)$ sont des invariants du type
réel $[F]$.
\end{prop}

\pr On procède comme pour la preuve du cas entier (proposition
\ref{prop:dec_can}, p.~\pageref{prop:dec_can}).

\medskip
\begin{rema}
Si $F \in \GL_n(\Z)$ représente un type isodual, ses composantes
réelles  sont évidemment des invariants de la 
classe entière $[F]$.
\end{rema}

Soit $[F]$ un type algébrique réel. Par un choix convenable d'une base  
adaptée à la décomposition canonique $\R^n=\bigoplus_{k,l} W_{k,l}(F)$
(proposition~\ref{prop:dec_can_R}), on voit  que $F$ est 
$\GL_n^\pm$-équivalente à une matrice décomposée $\oplus_{k,l} F_{k,l}$
avec $|\det F_{k,l}|=1$. Les types réels $[F_{k,l}]$ ne dépendent
que de $[F]$ et  seront appelés {\em composantes canoniques de $[F]$}.
On a  
\begin{equation}
\label{equa:m_kl}
[F]=\bigoplus_{k,l}[F_{k,l}] 
\esp \mathrm{avec} \esp 
\chi_{F_{k,l}F^\vee_{k,l}}= \Psi_{k,l}^{m_{k,l}(F)}
~~
((k,l)\in\mcE),
\end{equation}
où $m_{k,l}(F)$ est la multiplicité de $\Psi_{k,l}$ comme facteur 
de $\chi_{FF^\vee}$. Les $m_{k,l}(F)$ sont presque tous nuls.

\begin{theo}[structure des types géométriques]
\label{theo:geom_VF}
Soit $[F]$ un type algébrique réel et soient $(F_{k,l})$ des 
représentants des composantes canoniques de $[F]$, comme 
ci-dessus. %
Alors l'espace symétrique $V_F$ est isométrique au  produit 
riemannien
$\prod_{k,l} V_{F_{k,l}}$ (indexé par les %
$(k,l)\in\mcE$ tels que $m_{k,l}(F)>0$) et il existe
$(p,q) \in\N^2$, $g\in\N$ et $(p_{k,l},q_{k,l})\in \N^2$
($(k,l)\in\mcE,~k\geq 3$) tels que
$$
\left\{
\renewcommand{\arraystretch}{1.2}
\begin{array}{ll}
V_{F_{1,1}}\simeq \SO_0(p,q)/\SO(p)\times\SO(q) & (p+q=m_{1,1}(F)),\\
V_{F_{2,1}}\simeq \Sp_{2g}(\R)/\U(g)  & (2g=m_{2,1}(F)),\\
V_{F_{k,l}}\simeq \SU(p_{k,l},q_{k,l})/\SSS(\U(p_{k,l})\times \U(q_{k,l})) & 
(k\geq 3, ~p_{k,l}+q_{k,l} =m_{k,l}(F)).
\end{array} 
\right.
$$
\end{theo}

\begin{defi}
Soit $\mcE^H=\{(k,l)\in\mcE;k\geq 3\}$. 
Si $[F]$ est un type algébrique réel de rang $n$. La famille
$$((p,q);g;(p_{k,l},q_{k,l})_{(k,l)\in\mcE^H})$$
exhibée dans le  
théorème~\ref{theo:geom_VF} est appelée {\em signature} de
$[F]$. Noter que le rang $n$ de $[F]$ vérifie la relation
$n=p+q+2g+2 \sum_{(k,l)\in\mcE^H} (p_{k,l}+q_{k,l}).$
 Les espaces
symétriques $V_{F_{1,1}}$, $V_{F_{2,1}}$ et $V_{F_{k,l}}$ ($k\geq 3$) du
théorème~\ref{theo:geom_VF} sont appelés respectivement
{\em composante symétrique, composante alternée et composantes
hermitiennes} de $V_F$. 
\end{defi}

\begin{coro}
\label{coro:dim_VF}
Soit $[F]$ un type  réel de signature 
$((p,q);g;(p_{k,l},q_{k,l})_{(k,l)\in\mcE^H})$. Alors la dimension
de l'espace symétrique $V_F$ est donnée par 
$$
\dim V_F= pq + g(g+1) + 2 \sum_{(k,l)\in\mcE^H} p_{k,l} q_{k,l}.
$$
\end{coro}

{\it Preuve du théorème~\ref{theo:geom_VF}.} On peut supposer
que $F=\oplus_{k,l}F_{k,l}$, somme indexée par les couples
$(k,l)\in\mcE$ tels que $m_{k,l}=m_{k,l}(F)>0$. 
Par suite (voir preuve de l'assertion b), proposition~\ref{prop:dec1} 
p.~\pageref{prop:dec1}), on a 
$V_F=\{\oplus A_{k,l};A_{k,l} \in V_{F_{k,l}}\}$.
Comme l'ensemble des matrices diagonales par blocs
$$\bigoplus_{j=1}^N P_{n_j}= \{\oplus_{j=1}^N A_j; A_j\in P_{n_j}, 
1\leq j \leq  N\} \subset P_n
\esp (\sum_{j=1}^N n_j = n)$$
donne un plongement isométrique du produit riemannien
$\prod _{j=1}^N P_{n_j}$ dans $P_n$ (proposition~\ref{prop:som_prod_se_sym}-a), 
on voit que $V_F$ est 
isométrique au produit  des $V_{F_{k,l}}$. On est donc ramené à 
l'étude des composantes, c'est-à-dire au cas où le 
polynôme minimal de $R=FF^\vee$ vaut $\mu_R=\Psi_{k,l}$ pour un 
$(k,l)\in\mcE$.\\
\indent Si $\mu_R=x+1$, alors $F$ est symétrique et on peut 
prendre  $F=I_{p,q}$
($p+q=n$). On a $\OO_F(\R)=\OO^t_F(\R)=\OO(p,q)$ 
et $V_F$ est une orbite de la composante neutre $\SO_0(p,q)$ 
(prop~\ref{prop:VF_orbite}). En considérant la matrice identité
$I_n\in V_F$, on trouve que $V_F\simeq \SO_0(p,q)/
\SO(p)\times\SO(q)$. Si $\mu_R=x+1$, alors $F$ est équivalente à 
$J_{2g}$ ($2g=n$) et un argument analogue montre que
$V_F\simeq \Sp_{2g}(\R)/\U(g)$.\\
\indent Supposons maintenant que $\mu_R=\Psi_{k,l}$ pour un 
$(k,l)\in\mcE$ avec $k\geq 3$. Quitte à conjuguer, on peut supposer 
que  $R=R(2\pi l/k)\oplus \ldots \oplus R(2\pi l/k)$ 
($m_{k,l}$ termes) avec
$$
R(\theta)=\left(
\begin{array}{cc}
\cos \theta  & - \sin \theta \\
\sin \theta  & \cos \theta 
\end{array}
\right).
$$
La sous-algèbre $\mcC_R$ de $M_n(\R)$ des matrices qui commutent 
avec $R$ coïncide avec le commutant de la structure complexe
$R(\pi/2)\oplus \ldots \oplus R(\pi/2)$. Soit  $m=n/2$. Il
est clair que l'application $\kappa : M_m(\C) \to M_n(\R)$
définie par 
\begin{equation}
\label{equa:kappa}
\kappa(Z)=
\left(
\left(
\begin{array}{cc}
\re z_{i,j} & -\im z_{i,j}\\
\im z_{i,j} & \re z_{i,j}
\end{array}
\right)
\right)_{1\leq i,j\leq m}
~~ (Z=(z_{i,j})_{1\leq i,j\leq m} \in M_m(\C))
\end{equation}
(matrice formée de blocs $2 \times 2$) induit un isomorphisme 
de $\R$-algèbres de $M_m(\C)$ sur $\mcC_R$, Pour tout $Z\in M_m(\C)$,
en posant $Z^*=\overline{Z'}$, on a 
$$
\kappa(Z^*)=\kappa(Z)',~
\det (\kappa(Z)) =|\det(Z)|^2
~\mathrm{et}~
\tr(\kappa(Z)) = \tr(Z) + \tr(Z^*).$$
Revenons à la détermination de $V_F$. Compte tenu de $R^\vee=R$,
on vérifie aisément que $F\in\mcC_R$ et que $V_F$, $\OO_F(\R)$
et $\OO_F^t(\R)$ sont inclus dans $\mcC_R$ (en particulier
$\OO_F(\R)=\SO_F(\R)$). Pour $M\in\mcC_R$, posons 
$M_\C=\kappa^{-1}(M)\in M_m(\C)$. La relation $F=RF'$ s'écrit 
dans $M_m(\C)$ sous la forme
$
e^{- i \pi l/k}F_\C =(e^{- i \pi l/k}F_\C)^*.
$
Par suite, il existe  $(p,q)\in \N^2$ ($p+q=m$) et $Q\in \GL_m(\C)$ 
(avec $|\det Q|=1$) tels que 
$Q (e^{- i \pi l/k}F_\C)Q^* = I_{p,q}$. On peut donc supposer que 
$F_\C=e^{i \pi l/k} I_{p,q}$ , c'est-à-dire 
\begin{equation}
\label{equa:I_pq_herm}
F=I_{p,q}\otimes R_{\pi l/k}=R_{\pi l/k} \oplus \ldots  \oplus R_{\pi l/k}
\oplus (-R_{\pi l/k})\ldots  \oplus (-R_{\pi l/k}).
\end{equation}
Dans ces conditions, on a $P\in \SO^t_F(\R)$ \ssi
$P_\C\in \U(p,q)$ et le type géométrique est alors donné, d'après 
la proposition~\ref{prop:VF_orbite}, par
$$V_F=\{\kappa(UU^*); U\in \U(p,q)\}
\simeq \SU(p,q)/\SSS(\U(p)\times \U(q)).$$

\smallskip
\begin{rema}
La preuve précédente (voir aussi proposition~\ref{prop:dec1}, assertion c))
montre que la composante neutre du groupe $\OO_F(\R)$ est isomorphe
à un produit de groupes de la forme $\SO_0(p,q)$, 
$\Sp_{2g}(\R)$ ou $\U(p,q)$. On notera que {\em la dimension du noyau 
de l'action de $\OO^t_F(\R)$ sur $V_F$ est égale au nombre de 
composantes hermitiennes de $V_F$}.
\end{rema}

Une autre conséquence de la preuve du théorème~\ref{theo:geom_VF}
est la détermination des types algébriques réels (voir notamment
(\ref{equa:I_pq_herm})).  Posons
\begin{equation}
\label{equa:type_Pkl}
P_{2;k,l}=  \left(
\begin{array}{cc}
\cos \pi l/k  & - \sin \pi l/k \\
\sin \pi l/k  & \cos \pi l/k 
\end{array}
\right)
\esp ((k,l)\in \mcE^H).
\end{equation}

\begin{prop}[classification des types algébriques réels]
\label{prop:types_alg_reels}
Soit $[F]$ un type algébrique réel et soit 
$((p,q);g;(p_{k,l},q_{k,l})_{(k,l)\in\mcE^H})$ la signature de $[F]$.
On a alors la décomposition 
$$[F]=p[I_1]\oplus q [-I_1]\oplus g [J_2]
\bigoplus_{(k,l)\in\mcE^H} 
\left(p_{k,l}[P_{2;k,l}] \oplus q_{k,l}[-P_{2;k,l}]\right).$$
Deux types algébriques réels sont équivalents \ssi ils ont la 
même signature. 
\end{prop}

En particulier tout type réel $F$ de rang 2 avec $FF^\vee$ d'ordre $k$ 
est équivalent à 
$\pm I_2$ ou $I_{1,1}$ si $k=1$, à $J_2$ si $k=2$ et à 
 $\pm P_{2;k,l}$ ($(k,l)\in\mcE^H$) si $k\geq 3$. 
Ces types réels
admettent des représentants entiers uniquement pour $k=1$, $2$
et $6$ ($P_{2;6,1}$ est équivalent à $F_2$).

\subsection{Paramétrage des types géométriques}
\label{subs:param}
Soit $\G$ l'un des groupes $\Sp_{2g}(\R)$, $\SO_0(p,q)$ ou 
$\SU(p,q)$ et soit $\K$ un compact maximal de $\G$. 
Il est bien connu (voir ci-dessous) 
que l'espace symétrique $\G/\K$ associé à $\G$  est isomorphe
comme espace homogène à un  ouvert  $\Om_\G$ d'un espace numérique.
Nous obtenons ainsi,   {\it via} le 
théorème~\ref{theo:geom_VF}, un paramétrage 
des matrices de Gram des réseaux isoduaux. En effet,
soit $V_{F_{k,l}}$ une composante isomorphe à $\G/\K$
(voir théorème~\ref{theo:geom_VF}).
Alors le choix de  points bases respectifs dans $\Om_\G$ et $V_{F_{k,l}}$ 
définit, grâce à l'action de $\G$, un plongement équivariant 
d'image $V_{F_{k,l}}$ de $\Om_\G$ dans l'espace des matrices de Gram.
\par
La  composante alternée sera paramétrée par l'espace de 
Siegel $\mfH_g$. %
Rappelons que   $\mfH_g= \{Z = X + iY; X,Y \in\sym_g(\R), Y>0\}$
et  que le groupe symplectique $\Sp_{2g}(\R)$ agit sur $\mfH_g$ par 
homographies :
\begin{equation}
\label{equa:homographies} 
{\small 
P\cdot Z = (AZ+B)(CZ+D)^{-1} 
\esp 
(P=
\left(
\begin{array}{cc}
A & B\\
C & D
\end{array}
\right)
\in\Sp_{2g}(\R),~Z\in\mfH_g).
}
\end{equation}
Soit $\mfS_g = V_{J_{2g}} = \{A\in P_{2g};AJ_{2g}A=J_{2g}\}$. 
On obtient  un plongement équivariant $\sig_g : \mfH_g\to P_{2g}$ 
d'image $\mfS_g$ en posant $\sig_g(P\cdot iI) = PP'$~: 
\begin{equation}
\label{equa:Sig_g}
\small{
\sig_g(X+iY)=
\left(
\begin{array}{cc}
Y+XY^{-1}X & XY^{-1}\\
Y^{-1}X & Y^{-1}
\end{array}
\right)
~~(X+iY\in\mfH_g),
\esp
\mfS_g =\sig_g(\mfH_g).
}
\end{equation}
Noter que la métrique induite sur $\mfH_g$ par la métrique
riemannienne (\ref{equa:met_Pn}) de $P_{2g}$ est donnée par 
\begin{equation}
\label{equa:met_Hg}
ds^2= 2\tr(Y^{-1}dXY^{-1}dX +Y^{-1}dYY^{-1}dY)
\esp (X+iY \in\mfH_g).
\end{equation}

\par
Passons à la  composante symétrique. Posons 
$V_{p,q}=V_{ I_{p,q} }$. Soit $(p,q)$ avec 
$p>0$, $q>0$ et soit $\R^{p,q}$ l'espace $\R^n = \R^p\times\R^q$
 ($n=p+q$)  muni de la forme usuelle de signature $(p,q)$.
 Le groupe $\SO_0(p,q)$ agit sur  la grassmannienne $\mcG_{p,q}^-$
des sous-espaces négatifs de dimension $q$ de $\R^{p,q}$ et on a
$\mcG_{p,q}^-\simeq \SO_0(p,q)/\SO(p)\times\SO(q)$.
En remarquant que tout élément de $\mcG_{p,q}^-$ admet une unique
base de la forme $\smallmat{X}{\!}{I_q}{\!}$
avec $X'X -I_q <0$, on voit que  $\mcG_{p,q}^-$ est isomorphe
à 
\begin{equation}\label{equa:Om_pq}
{\mcV}_{p,q}=\{X\in M_{p,q}(\R); X'X-I_q< 0\}
\end{equation}
muni de l'action de $ \SO_0(p,q)$ par homographies analogue 
à (\ref{equa:homographies}). En choisissant comme points bases
$0\in\mcV_{p,q}$ (correspondant à $\{0\}\times\R^q\in\mcG_{p,q}^-$)
et $I_n\in V_{p,q}$, on trouve un plongement équivariant 
$\ph_{p,q}$ de $\mcV_{p,q}$ dans $P_n$ d'image $V_{p,q}$~: 
\begin{equation}\label{equa:Vpq}
\renewcommand{\arraystretch}{1.1}
\small{
\left\{
\begin{array}{l}
\ph_{p,q}(X)=
\left( 
\begin{array}{cc}
(I+XX')(I-XX')^{-1} & 2X(I-X'X)^{-1}\\
2(I-X'X)^{-1}X'& (I+X'X)(I-X'X)^{-1}
\end{array}
\right)~~~
(X\in\mcV_{p,q}),\\
V_{p,q}=\ph_{p,q}(\mcV_{p,q}).
\end{array}
\right.
}
\end{equation}
La métrique riemannienne induite par ce plongement est donnée 
par 
\begin{equation}\label{equa:met_Xpq}
ds^2= 8 \tr [(1-XX')^{-1}dX(1-X'X)^{-1}dX'] \esp
(X\in\mcV_{p,q}).
\end{equation}

\begin{rema}
Pour $q=1$, à partir de 
$(I-XX')^{-1}=I+XX' (1-X'X)^{-1}$, on voit que 
la métrique induite $ds^2$  sur $\mcV_{n-1,1}$ 
vérifie
\begin{equation}\label{equa:met_Klein}
\frac{ds^2}{8} = \frac{1}{(1-|X|^2)^2} 
((1-|X|^2)|dX|^2 + (X'dX)^2)
\esp (|X|^2= X'X < 1).
\end{equation}
Il s'agit de la métrique hyperbolique exprimée dans le modèle de 
Klein. La métrique induite sur $V_{n-1,1}$ est à courbure 
constante $-1/8$.
\end{rema}

Il peut être commode de paramétrer $V_{n-1,1}$ par le demi-espace
supérieur en composant le plongement 
$\ph_{n-1,1}$ avec une isométrie entre le demi-espace et le modèle
de Klein (et de même pour d'autres modèles). Par exemple, pour $n=3$ :
\begin{equation}
\renewcommand{\arraystretch}{1.1}
\label{equa:H2_V21}
{\small
\frac{1}{(\Im z)^2}
\left(
\begin{array}{ccc}
2(\Re z)^2  + (\Im z)^2 &  (|z|^2-1)\Re z &  (|z|^2+1)\Re z\\
(|z|^2-1)\Re z & \frac{1}{2}(|z|^4+1) -(\Re z)^2  &  
\frac{1}{2}(|z|^4 - 1)\\
(|z|^2+1)\Re z & \frac{1}{2}(|z|^4 - 1) & 
\frac{1}{2}(|z|^4+1)  + (\Re z)^2  
\end{array}
\right)
~~~ (\Im z >0),
}
\end{equation}
paramétrage que l'on peut aussi construire directement à partir
de l'isomorphisme classique entre $\PSL_2(\R)$ et $\SO_0(2,1)$.
\par

Enfin pour le groupe $\SU(p,q)$ et les composantes hermitiennes,
on a une description analogue à $\SO_0(p,q)$ à partir de la
grassmannienne des sous-espaces négatifs de dimension $q$ de 
l'espace hermitien standard de signature $(p,q)$.
Comme dans le cas réel, cette grassmannienne s'identifie 
à 
\begin{equation}\label{equa:mcWpq}
\mcW_{p,q}= \{Z\in M_{p,q}(\C); Z^*Z -I_q < 0\}
\end{equation}
muni de l'action de $\SU(p,q)$ par homographies et on trouve,
{\it via} les matrices hermitiennes et le morphisme d'algèbres 
$\kappa$ (voir équation (\ref{equa:kappa})),
un plongement équivariant de $\mcW_{p,q}$ dans $P_n$  ($n=2(p+q)$) 
d'image   $W_{p,q}=\kappa (\{PP^*;P\in\SU(p,q)\})$, avec  
\begin{equation}\label{equa:Wpq}
\small{
W_{p,q}=%
\kappa \left(
\left\{
\left( 
\begin{array}{cc}
(I+ZZ^*)(I-ZZ^*)^{-1} & 2Z(I-Z^*Z)^{-1}\\
2(I-Z^*Z)^{-1}Z^*& (I+Z^*Z)(I-Z^*Z)^{-1}
\end{array}
\right);
Z \in \mcW_{p,q}
\right\}
\right).
}
\end{equation}
Quant à la métrique induite sur $\mcW_{p,q}$, elle s'exprime
(compte tenu notamment de $\tr\kappa(W)=\tr W +\tr W^*$) 
par
\begin{equation}\label{equa:met_Ypq}
ds^2 = 16 \tr [(1-ZZ^*)^{-1}dZ (1-Z^* Z)^{-1}dZ^*].
\end{equation}
Nous utiliserons essentiellement le cas particulier $(p,q)=(1,1)$, soit
$$
\small{
W_{1,1}= \left\{
\frac{1}{1-|z|^2}
\left(
\begin{array}{cccc}
1+|z|^2 &  0 & 2\Re z & -2\Im z \\ 
0 & 1+|z|^2 & 2\Im z & 2\Re z \\ 
2\Re z &  2\Im z &  1+|z|^2 & 0 \\ 
-2\Im z & 2\Re z & 0 & 1+|z|^2
\end{array}
\right) ; 
|z| < 1 \right\},
}
$$
que l'on peut aussi paramétrer par le demi-plan de Poincaré
$\mfH_1$  : 
\begin{equation}
\label{equa:W11}
\small{
W_{1,1}= \left\{
\frac{1}{2 y}
\left(
\begin{array}{cccc}
1+|z|^2 &  0 & |z|^2-1 & 2x\\ 
0 & 1+|z|^2 & -2x & |z|^2-1 \\ 
|z|^2-1  &  -2x  &  1+|z|^2 & 0 \\ 
2x & |z|^2-1  & 0 & 1+|z|^2
\end{array}
\right) ; 
z=x+iy\in\mfH_1 \right\}.
}
\end{equation}

\begin{rema}
\label{rema:plans_hyperboliques}
En petite dimension, il existe des isomorphismes entre certains
groupes symplectiques réels, orthogonaux ou unitaires. 
Les types  géométriques correspondants, 
muni de la métrique induite par $P_n$, sont homothétiques
car copies du même espace symétrique irréductible.
Ainsi, grâce aux isomorphismes $\PSL_2(\R) \simeq \SO_0(2,1)
\simeq \PSU(1,1)$, 
chacune des familles $\mfS_g$, $V_{p,q}$ et $W_{p,q}$ contient 
un plan à courbure constante négative :
$\mfS_1$, $V_{2,1}$ et $W_{1,1}$, de courbure 
respective  $-1/2$, $-1/8$ et $-1/4$
(voir \eqref{equa:met_Hg}, \eqref{equa:met_Klein} et
\eqref{equa:met_Ypq}). Les types $V_{3,2}$ et  $\mfS_2$ (isomorphisme 
$\SO_0(3,2)\simeq  \PSL_4(\R)$) sont homothétiques (avec un rapport
de $\sqrt{2}$) mais pas isométriques. 
Enfin,   les types  $W_{2,2}$ et $V_{4,2}$ (isomorphisme  
$\SO_0(4,2) \simeq  \SU(2,2)/\{\pm I_4\}$) sont isométriques, 
au moins abstraitement. 
\end{rema}

Pour paramétrer un type géométrique $V_F$, il suffit de scinder
$F$ comme type réel, c'est-à-dire trouver $P\in\GL_n^\pm(\R)$ 
tel que $F = P\cdot(\oplus_{k,l}F_{k,l})$ avec 
$F_{1,1}=I_{p,q}$, $F_{2,1}=J_{2g}$ et $F_{k,l}=p_{k,l} P_{2;k,l}\oplus
q_{k,l}(-P_{2;k,l})$ pour $(k,l)\in\mcE^H$.
 On a alors 
\begin{equation}
\label{equa:param}
V_F=P\cdot(\bigoplus_{(k,l)\in\mcE} V_{F_{k,l}})
\end{equation}
où les blocs diagonaux $V_{F_{k,l}}$ sont de la forme
(\ref{equa:Sig_g}), (\ref{equa:Vpq}) ou (\ref{equa:Wpq}).

\subsection{Automorphismes et inclusion}
Soit $F\in\GL_n(\Z)$ représentant un type algébrique. 
Nous appellerons {\em groupe
d'automorphismes de $V_F$}, noté  $\Ga_F$, le sous-groupe de
de  $\GL_n(\Z)$ formé des  éléments qui fixent $V_F$ point par point.
Il s'agit d'un groupe fini et on  a toujours $FF^\vee\in\Ga_F$.
De plus, la classe de conjugaison de $\Ga_F$ dans $\GL_n(\Z)$ est un 
invariant du type géométrique représenté par $V_F$.
\par
Pour la détermination de $\Ga_F$, il sera souvent utile de remarquer
qu'un élément   $T\in\GL_n(\Z)$ (et plus généralement
$T\in\GL_n^\pm(\R)$) fixe  $V_F$ point par point
\ssi $T$ fixe  un point $A\in V_F$ ainsi que  tous les vecteurs 
de l'espace  tangent $T_AV_F$.

\begin{prop}[critère d'inclusion]
\label{prop:inclusion}
Soit $F\in\GL_n(\Z)$  représentant un type algébrique, soit 
$A\in V_F$  et soit $G\in\GL_n(\Z)$.
Alors $[G]$ est un type algébrique tel que $V_F\subset V_G$ 
(resp. tel que $A\in V_G$) \ssi $FG^\vee\in \Ga_F$
(resp. \ssi $FG^\vee\cdot A = A$).
\end{prop}

\pr Pour $A\in V_F$, on a 
$F^\vee(AG^\vee A G^{-1})F' = A^{-1}FG^\vee A (FG^\vee)'.$
Par suite $AG^\vee A =G$ \ssi $(FG^\vee)\cdot  A =A$, d'où le 
résultat (noter que~$G$ représente un type algébrique \ssi 
$V_G\neq  \emptyset$).
 
Le critère d'inclusion permet d'ordonner les types géométriques.
Un type $[F]$ étant donné,
l'ensemble des $G\in\GL_n(\Z)$ avec $GG^\vee$ d'ordre 
fini et $V_F\subset V_G$ coïncide avec
$$\{F'T^\vee;T\in\Ga_F\}.$$

\begin{prop}[maximalité des types orthogonaux et symplectique]
\label{prop:orth_symp_max}
Soient $[F]$ et $[G]$ des  types algébriques.
\begin{enumerate}
\item  Si $\Ga_F=\{\pm I_n\}$, alors $[F]$ est un type orthogonal ou 
symplectique (\cadab $FF^\vee=\pm I_n$) et $V_F$ est maximal. De plus,
$V_F$ est équivalent à $V_G$ \ssi $[G]=[\pm F]$.
\item Si $[F]$ est le type symplectique ou un type orthogonal de signature
distincte de $(n,0)$, $(0,n)$ ($n\geq 2$) et  $(1,1)$, alors 
$\Ga_F=\{\pm I_n\}$ et   $V_F$ est maximal.
\end{enumerate}
\end{prop}

\pr L'assertion (1) est immédiate car $FF^\vee \in\Ga_F$ et d'après 
ce qui précède les seuls types $[H]$ tels que $V_F\subset V_H$
sont donnés par $H=\pm F$.\\
\indent Pour la preuve de (2), on utilise les types réels. Il 
existe $P\in\GL_n^\pm(\R)$ tel que $P\cdot F = F_0$ avec 
$F_0=J_{2g}$  ou $I_{p,q}$
(bien sûr, si le type est symplectique ou symétrique indéfini impair,
on peut prendre $P\in\GL_n(\Z)$).  
Il suffit de prouver (dans les conditions de l'énoncé)
 que tout $T \in\GL_n^\pm(\R)$ qui fixe point 
par point $V_{F_0}= \mfS_g$ ou $V_{p,q}$ vaut $\pm I_n$. Un tel  
élément $T$ doit vérifier $T\cdot I_n = TT'=I_n$ (\cadab $T\in
\OO(n)$ et $TXT'=X$ pour tout $X\in T_{I_n}V_{F_0}$. Mais, d'après la
remarque~\ref{rema:espace_tangent},  cet espace
tangent est donné par 
$$ T_{I_n}V_{F_0} = \{X-F_0XF'_0; X\in\sym_n(\R)\}.$$ 
En prenant $X=uu'$ ($u\in\R^n$), on a que  
$Tuu'-TF_0 u u'F'_0 = uu'T - F_0 u u'F'_0 T$. Posons
$\langle u,v\rangle = u'v$ ($u,v \in\R^n$). Pour tout $u\in\R^n$ tel que 
$F_0[u]=0$ la relation précédente entraîne
$$\langle u,u \rangle  Tu  = \langle u,Tu \rangle  u - 
\langle F_0 u ,Tu \rangle  F_0 u.$$
En multipliant par $(F_0u)'$, on trouve 
$ 2 \langle u,u\rangle
\langle Tu,F_0 u\rangle = 0$ d'où (si $u\neq 0$) 
$ \langle Tu,F_0 u\rangle = 0$. Finalement, $T$ doit laisser 
stable toute droite isotrope de~$F_0$. Dans le cas symplectique,
cette condition entraîne évidemment que $T$ est une homothétie. 
On vérifie facilement qu'il en de même dans le cas symétrique indéfini, 
à l'exception de la signature $(1,1)$ qui ne comprend que 
deux droites isotropes (ce type n'est d'ailleurs pas maximal, voir
table~\ref{tabl:types_geom_1_3}).

\begin{prop}
\label{prop:aut_ordre_4}
Soient $[F]$ et $[G]$ des  types algébriques. On suppose que 
$FF^\vee\neq \pm I_n$ et que le cardinal du groupe $\Ga_F$ vaut 4. 
Alors le type  géométrique $V_F$ est maximal. De plus,
$V_F$ est équivalent à $V_G$ \ssi $[G]=[\pm F]$ ou $[G]=[\pm F']$.
\end{prop}

\pr On a $\Ga_F=\{\pm I_n,\pm FF^\vee\}$. Les seuls types $[H]$
tels que $V_F\subset V_H$ vérifient  donc $H=\pm F$ ou  $H=\pm F'$.

\subsection{Classification : méthode et notations}
Nous déterminerons les types géométriques à partir des composantes
réelles  des types algébriques 
(tables~\ref{tabl:types_1_3}-\ref{tabl:types_7_b}),
comme indiqué  aux \S\S~\ref{subs:composantes} et~\ref{subs:param}. 
Le scindement sur $\R$ d'un type algébrique $[F]$ permet d'expliciter
$V_F$ comme produit d'espaces symétriques plongé dans l'espace des 
matrices de Gram (voir en particulier la 
relation~(\ref{equa:param})). 
\par
Pour alléger les notations, nous
poserons souvent $(F)=V_F$. La relation évidente $(-F)=(F)$ 
en entraîne d'autres. Par exemple, pour un  type  scindé 
$F=F_1\oplus F_2$  dans les conditions de la 
proposition~\ref{prop:dec1} ($\chi_{F_1F_1^\vee}$ et $\chi_{F_2F_2^\vee}$ 
premiers entre eux), on a
$$(\pm F_1 \oplus \pm F_2)=(F_1\oplus F_2).$$
Plus généralement, si $F$ et $G$ se scindent sur $\R$ à l'aide
d'un même élément $P\in\GL_n^\pm(\R)$ en composantes
$F_i$ et $G_i$ telles que $F_i=\pm G_i$  ($i=1,\ldots,k$),
avec $(\chi_{F_1F_1^\vee})_{i=1,\ldots,k}$ mutuellement
 premiers entre eux, alors $(F)=(G)$.  Ainsi, les 
nombreuses coïncidences entre les types géométriques
réduisent significativement la liste de ces types. 
\par  
Pour chaque dimension $n$ examinée, nous commencerons par établir
la liste des types géométriques (seulement les types principaux 
si $n\geq 5$), lesquels  seront ensuite ordonnés (à équivalence près)
grâce au critère d'inclusion (proposition~\ref{prop:inclusion}).
Nous obtiendrons ainsi les types maximaux jusqu'en rang 7.
La détermination de toutes les relations  d'inclusion 
nécessite la connaissance des groupes 
d'automorphismes $\Ga_F$. Certains de ces groupes seront précisés dans 
les tables  \ref{tabl:types_geom_1_3}-\ref{tabl:types_geom_6}
et~\ref{tabl:types_geom_7}. En particulier nous noterons
$\G(n)$ le groupe d'automorphismes du réseau $\Z^n$. 
Il s'agit du produit semi-direct défini
par l'action naturelle du groupe symétrique $\mcS_n$ sur
$(\Z/2\Z)^n$ ; le cardinal de $\G(n)$ vaut donc $2^n n!$.
Le groupe cyclique  (resp. diédral) d'ordre $k$ sera noté
$\C_k=\Z/k\Z$ (resp. $\D_{k}$, $k$ pair). 
\par 
Un des buts de cet article est l'étude de la densité 
des réseaux isoduaux. 
Rappelons que $V_F$ satisfait un critère de compacité de Mahler
(proposition~\ref{prop:VF_Mahler}). L'invariant d'Hermite $\mu$ admet
donc une valeur maximale sur $V_F$ : 
$$\mu_F=\max_{A\in V_F} \mu(A).$$
Pour chaque type géométrique, nous chercherons à estimer (et si 
possible à déterminer) cette \og constante d'Hermite relative\fg.
Nous dirons qu'elle est atteinte {\em de façon unique 
à équivalence près dans $V_F$} si deux points de $V_F$ 
réalisant $\mu_F$ sont équivalents par un élément 
$P\in \GL_n(\Z)$ {\em qui laisse stable $V_F$} (par exemple
$P\in \OO^t_F(\Z)$). Pour les types symétriques indéfinis
$F=I_{p,q}$ et $F=U_{2p}$ ($pq>0$), les $\mu_F$ seront notées
simplement $\mu_{p,q}$  et  $\mu_{p,p}^{\rm II}$.
La constante $\mu_{n,1}$ du  type orthogonal lorentzien impair 
est déterminée jusqu'en rang 12 dans \cite{Bavard2007a}
(voir en particulier la table 1, p.~40). Pour de nombreux 
types géométriques scindés (définition ci-dessous), le 
calcul de  $\mu_F$ se ramène aux dimensions inférieures.
\begin{defi}
\label{defi:tg_scinde}
Soit $[F]$ un type algébrique (entier). Le type géométrique
$V_F$ est {\em scindé} (ou {\em décomposé} s'il existe deux 
types algébriques entiers $[F_1]$ et $[F_2]$ tels que 
$V_F=V_{F_1}\oplus V_{F_2}$.
\end{defi}

\begin{rema}
Soit $[F]$ un type algébrique entier de rang $n$ et soit $k$ la 
partie entière de $n/2$. Si $V_F$ est scindé alors $\mu_F\leq \ga_k^{\rm isod}$.
\end{rema}
 
\par 
Les tables \ref{tabl:types_geom_1_3}-\ref{tabl:types_geom_6}
et~\ref{tabl:types_geom_7}  résument la classification des types 
géométriques.  Nous y indiquons notamment 
\begin{itemize}
\item la structure métrique (ou \og type symétrique\fg, noté
t. s.) de $V_F$, 
\item la structure du groupe
 $\Ga_F$ ou à défaut son cardinal, 
\item une estimation de la  constante $\mu_F$,
\item  les types 
géométriques contenant $V_F$ à équivalence près. 
\end{itemize}
\par
Les constantes $\mu_F$ sont toutes déterminées  jusqu'en rang~4
(tables \ref{tabl:types_geom_1_3} et \ref{tabl:types_geom_4}).
En dimension supérieure, à défaut de la valeur exacte de $\mu_F$,
{\em nous donnons toujours une minoration par un maximum local de densité
sur $V_F$} (tables \ref{tabl:types_geom_5}, \ref{tabl:types_geom_6}
et~\ref{tabl:types_geom_7}). Cette dernière propriété variationnelle 
provient du fait
que le réseau correspondant est dans chaque cas parfait et eutactique 
relativement à $V_F$ (voir par exemple \cite[proposition~2.1-(1)]{Bavard1997a}),
ce qui se vérifie en déterminant les gradients des fonctions longueurs,
soit  à partir d'un paramétrage explicite du type $V_F$, soit
à partir de l'expression \og implicite\fg\ de la 
proposition~\ref{prop:esp_tangent}. 
\par
Concernant la structure métrique, nous noterons $\HH^n_a$ l'espace 
de dimension~$n$ 
complet et simplement connexe à courbure sectionnelle constante $-1/a$ ;
 on trouvera par exemple trois sortes de  
\og plans hyperboliques\fg :  $$\HH^2_2, ~ \HH^2_8
~\mathrm{ou}~ \HH^2_4,$$
à courbure $-1/2$, $-1/8$ ou $-1/4$  suivant qu'ils proviennent d'un 
type réel symplectique, symétrique ou hermitien (voir 
rem.~\ref{rema:plans_hyperboliques}). 
\par
 Les détails techniques (matrices de décomposition réelle
des types, détermination des inclusions, calculs des gradients
et vérification des  maxima locaux, etc) des différentes étapes 
ne seront pas  explicités. Cependant nous détaillerons certains cas 
particuliers intéressants : paramétrage remarquable d'un type, 
détermination non immédiate d'une constante $\mu_F$, types géométriques
maximaux non scindés dans l'annexe \ref{subs:ann_tg_non_scin} (pour les
dimensions 5, 6 et 7, autres que les symétriques ou symplectiques).

\subsection{Types géométriques de rang $n \leq 3$}
\label{subs:geom_1_2_3}
Pour $n=2$, on trouve 5 types géométriques 
$V_{I_2}=\{I_2\}$, $V_{1,1}$ (voir équation~(\ref{equa:Vpq})),
 $V_{1,1}^{II}=V_{U_2}$, 
$V_{F_2}=\{A_2\}$ avec 
$$
\renewcommand{\arraystretch}{1}
V_{1,1}^{II}= 
\small{\left\{
\left(
\begin{array}{cc}
u &  0\\
0 & 1/u\\
\end{array}
\right); u \in\R \right\}} ~~~
\mathrm{et}~~~
A_2=\small{
\frac{1}{\sqrt{3}}
\left(
\begin{array}{cc}
2 &  -1\\
-1 & 2\\
\end{array}
\right)}.
$$
ainsi qu'un unique type maximal : le type symplectique $\mfS_2=P_2$.

\begin{table}[b]%
\centering
\renewcommand{\arraystretch}{1.2}
\renewcommand{\tabcolsep}{0.3 em}
{\small 
\begin{tabular}{cccccccclcl}
\hline
$n$ & $d$  & $F F^\vee$ &$n^\circ$   & $F$ & $(F)$& dim. & t. s. &
$\Ga_F$ & $\mu_F$ &$\subset_{\rm eq}$ \\
\hline
1 & 1 & $I_1$ & 1 &$I_1$ & $\{I_1\}$ &  0 & $\R^0$ & $\C_2$ & 1 & max.\\
\hline
2 &1 & $I_2$ &1 &$I_2$ & $\{I_2\}$ &  0 & $\R^0$ & $\D_8$ & 1 & 1-4\\
  & & & 2&$I_{1,1}$ & $V_{1,1}$ & 1 & $\R$ & $\C_2^2$ & $2/\sqrt{3}$ & 2,4\\
  & & & 3&$U_2$&  $V_{1,1}^{\rm{II}}$& 1&$\R$  &$\C_2^2$ & 1 &3,4\\
  & 2&$-I_2$ &4&$J_2$&$P_2$ &2 & $\HH^2_2$& $\C_2$ & $2/\sqrt{3}$& max.\\
  & 6&$W_2$&5 &$F_2$ &$\{A_2\}$ &0 &$\R^0$ &$\D_{12}$& $2/\sqrt{3}$& 
2,4,5\\
\hline
3 & 1&$I_3$ &1&$I_3$&$\{I_3\}$ &  0 & $\R^0$ & $\G(3)$ & 1 & 1-3\\
 & & &2&$I_{2,1}$ & $V_{2,1}$ & 2 & $\HH^2_8$ & $\C_2$&
$\al$ & max.\\
& 2&$I_{1,2}$&3&$I_1 J_2$&$P_1  P_2$ & 2 &
 $\HH^2_2$ & $\C_2^2$&1&max.\\
& 4 & $X_3$&4&$G_3$ & $\{\La_3\}$&0 & $\R^0$ & $\C_2 \D_8$ & 
$\al$ &2,4\\
& 6 &$I_1 W_2 $ &5&$ I_1 F_2$& $\{I_1  A_2\}$ 
&  0 &$\R^0$& $\C_2 \D_{12}$ &1 & 2,3,5\\
\hline
\end{tabular}}
\caption{Types géométriques isoduaux de rang  $n\leq 3$}
\label{tabl:types_geom_1_3}
\end{table}

Pour $n=3$, on constate d'abord les relations 
$(F_3)=(I_3)=\{I_3\}$, $(I_1 \oplus F_2)=(-I_1 \oplus F_2)=
\{I_1\oplus A_2\}$ et $(G_3)=(H_3)=\{\La_3\}$ avec 
\begin{equation}\label{equa:L3}
\La_3 = \small{\frac{1}{2}\left( \begin{array}{ccc}
1 +\sqrt{2}  & \sqrt{2} & -\sqrt{2} \\
\sqrt{2} & 2\sqrt{2} & 0  \\
-\sqrt{2} & 0 & 2\sqrt{2}  \end{array} \right)}.
\end{equation}
Outre ces trois types de dimension $0$, on trouve deux types 
maximaux :
$$
\begin{array}{ll}
\renewcommand{\arraystretch}{1}
V_{2,1} =
\small{
\left\{
\frac{1}{1-|z|^2}
\left(
\begin{array}{ccc}
1+\Re (z^2)  & \Im (z^2) & 2 \Re z\\
\Im (z^2) & 1 -\Re (z^2)  & 2 \Im z \\
2\Re z & 2 \Im z & 1+|z|^2
\end{array}
\right); |z|^2 < 1  \right\}} ~~~
\mathrm{et}~~~ \\
\\
V_{I_1\oplus J_2}= P_1\oplus P_2 = \small{
\left\{
\frac{1}{\Im z}
\left(
\begin{array}{ccc}
\Im z &  0 & 0\\
0  & |z|^2 & \Re z \\
0 & \Re z & 1
\end{array}
\right) ; 
\Im z >0
\right\}
}
\end{array}
$$
(comparer avec \cite[théorème 1]{Conw-Sloa1994a}). Métriquement, 
ces deux types sont
des \og plans hyperboliques\fg\ de courbure $-1/8$ et $-1/2$ 
respectivement (voir rem.~\ref{rema:plans_hyperboliques}).
On sait (\cite[théorème~2]{Conw-Sloa1994a}) que la constante d'Hermite \og
isoduale\fg\ en dimension  3 vaut
\begin{equation}\label{equa:cte_alpha}
\ga_3^\mathrm{isod} = \mu(\La_3) =  \alpha 
\esp\mathrm{avec}\esp
\alpha = \frac{1+\sqrt{2}}{2}
\simeq 1,2071,
\end{equation}
atteinte par un unique réseau (à isométrie près) dont une 
matrice de  Gram est $\La_3$ (équivalente à  un point de $V_{2,1}$ mais
qui n'appartient pas à $V_{2,1}$). Rappelons (\cite{Vinberg1972b})
que l'action induite
par $\OO_{2,1}(\Z)$ sur l'ensemble des  droites  négatives de 
$\R^{2,1}$ (\cadab sur $V_{2,1}$) est celle d'un groupe de 
réflexions du plan hyperbolique, de type triangulaire
$(2,4,\infty)$. On retrouve 
facilement  la densité maximale dans $V_{2,1}$ 
par un simple  découpage du triangle fondamental
(\cite[p. 45]{Bavard2007a}).
\par 
La situation en rang $n\leq 3$ est résumée dans la table 
\ref{tabl:types_geom_1_3}, où l'entier $d$ désigne l'ordre de
$FF^\vee$.
\subsection{Types géométriques de rang $4$}
\label{subs:geom_4}
En dimension 4, nous trouvons 23 types géométriques  dont 4 types
maximaux : les types symétriques (excepté $I_4$) et le type symplectique. 
Les constantes $\mu_F$ sont toutes explicitées 
(voir table~\ref{tabl:types_geom_4}). 
Deux réseaux remarquables réalisent 
le maximum de~$\mu$ sur les types maximaux. Il s'agit de
$$ 
{\small
D_4=\frac{1}{\sqrt{2}}
\left( \begin{array}{rrrr}
2& 0& 1& 1\\
0& 2& -1& 1\\
1& -1& 2& 0\\
1& 1& 0& 2 \end{array} \right)
\esp\mathrm{et}\esp
\La_4=\left( \begin{array}{rrrr}
2  & 1 & 1 & 1\\
1 & 4/3 & 1/3 & 1/3 \\
1 & 1/3 & 4/3 & 1/3 \\
1 & 1/3 & 1/3 & 4/3 \end{array} \right).}
$$
Il est bien connu que $\mu(D_4)=\sqrt{2}=\ga_4^\mathrm{isod}=\ga_4$ 
(constante d'Hermite en dimension~4),
maximum de $\mu$ dans $\mfS_2$ et $V_{2,2}$, 
et on a $\mu(\La_4)=4/3$, maximum de $\mu$ dans $V_{3,1}$
(\cite{Bavard2007a}) et dans $V_{2,2}^{\rm{II}}$ (voir plus bas).
Le groupe d'automorphismes de $\La_4$ est résoluble non nilpotent
d'ordre 144 et d'abélianisé $(\Z/2\Z)^3$.
\par 
D'après la table~\ref{tabl:types_4}, il y a 28 classes de types
algébriques au signe près. Les relations suivantes réduisent à 
23 le nombre de types géométriques ($\simeq$ désigne l'équivalence
entière) :
\renewcommand{\arraystretch}{1.1}
$$\begin{array}{l}
(I_1\oplus F_3)=(I_4)=\{I_4\},\\
(I_1\oplus H_3)=(I_1\oplus G_3)=\{I_1\oplus \La_3 \},\\
(-I_2\oplus F_2)=(I_2\oplus F_2)=\{I_2\oplus A_2 \},\\
(N_4) \simeq (M_4), ~~ (O_4) \simeq (L_4) = \{D_4\}.
\end{array}$$
Ces 23 types  se distinguent par leur dimension, leur
structure métrique, la valeur de $\mu_F$ et le cardinal du groupe
$\Ga_F$ (voir table~\ref{tabl:types_geom_4}).
\par
Notons $V_{2,2}^{\rm{II}}=(U_4)$ le type symétrique pair
 ($U_4=U_2\otimes I_2$). Les types  $V_{2,2}$ et $V_{2,2}^{\rm{II}}$
sont  évidemment  paramétrés par $\mcV_{2,2}$ 
(équation~(\ref{equa:Om_pq})) {\it via}~(\ref{equa:Vpq}), mais 
aussi par $\mfH_1\times \mfH_1$. Un paramétrage
particulièrement agréable pour $V_{2,2}^{\rm{II}}$ s'obtient comme suit.
Considérons $M_2(\Z)$ muni de la base 
$\mcB=(E_{1,1},E_{1,2},E_{2,2},-E_{2,1})$
(les $E_{i,j}$ étant les matrices élémentaires) et de 
la forme quadratique $f$ définie par 
$$f(A)=2\det(A) \esp A\in M_2(\Z)),$$ 
forme dont la matrice vaut dans la base $\mcB$ vaut $U_4$. 
La multiplication (à droite ou à gauche) par un élément 
de $\SL_2(\R)$ est une isométrie réelle de $f$. Ainsi l'application
$\phi$ définie par 
\begin{equation}
\label{equa:SO22}
 \phi(A,B)(X)=AXB^{-1} \esp
((A,B)\in \SL_2(\R)\times \SL_2(\R),~X\in M_2(\Z))
\end{equation}
induit un isomorphisme entre les groupes $\SL_2(\R)\times
\SL_2(\R)/\{\pm(I_2,I_2)\}$ et $\SO_0(f,\R)$. En explicitant 
les matrices dans $\mcB$ de multiplication par les éléments
de $\SL_2(\R)$, on trouve facilement
un paramétrage équivariant du type symétrique pair :
\begin{equation}
\label{equa:V22_pair}
\small{
V_{2,2}^{\rm{II}}= \left\{
\frac{1}{\Im z \Im w}
\left(
\begin{array}{cccc}
|z|^2 &  -|z|^2 \Re w  & -\Re z \Re w & - \Re z \\ 
-|z|^2 \Re w & |zw|^2 & \Re z |w|^2 & \Re z \Re w  \\ 
-\Re z \Re w  &  \Re z |w|^2 &  |w|^2  & \Re w \\ 
-\Re  z & \Re z \Re w  & \Re w & 1
\end{array}
\right) ; 
(z,w)\in \mfH_1\times \mfH_1 \right\}.
}
\end{equation}
La métrique induite par ce plongement s'écrit 
$ds^2=4(|dz|^2/\Im z +|dw|^2/\Im w)$ et les types $V_{2,2}^{\rm{II}}$
et $V_{2,2}$ sont donc  isométriques au produit $\HH^2_{4}\times
\HH^2_{4}$ de deux plans à courbure constante $-1/4$.

\begin{rema} Noter également la description 
$V_{2,2}^{\rm{II}}= (-J_2\oplus I_2)\cdot(\mfS_2\otimes \mfS_2)$ 
qui redonne immédiatement les propriétés  métriques 
(proposition~\ref{prop:som_prod_se_sym}-b).
\end{rema}

\begin{table}[htb]
\centering
\renewcommand{\arraystretch}{1.2}
\renewcommand{\tabcolsep}{0.3 em}
{\small 
\begin{tabular}{ccccccclcl}
\hline
 $d$  & $F F^\vee$ &$n^\circ$   & $F$ & $(F)$& dim. & t. s. &
$\Ga_F$ & $\mu_F$ &$\subset_{\rm eq}$ \\
\hline
1&$I_4$ &1&$I_4$&$\{I_4\}$ &  0 & $\R^0$ &$\G(4)$ &1 & 1-10,13\\
&&2&$I_{3,1}$ & $V_{3,1}$ & 3 & $\HH^3_8$ & $\C_2$& 4/3& max.\\
&&3&$I_{2,2}$ & $V_{2,2}$ & 4 & $\HH^2_4  \HH^2_4$  
& $\C_2$& $\sqrt{2}$& max.\\
&&4&$U_4$ & $V_{2,2}^{\rm{II}}$& 4  & $\HH^2_4  \HH^2_4$  
& $\C_2$& 4/3& max.\\
\hline
2&$I^-_4$&5&$J_4$ & $\mfS_2$ & 6  & $\mfH_2$ & $\C_2$&$\sqrt{2}$&max.\\
&$I_{2,2}$&6&$I_2  J_2$ & $\{I_2\}  P_2$& 2 & $\HH^2_2$ 
&$\D_8\C_2$ & 1 &5-8\\
& &7&$I_{1,1}  J_2$ & $V_{1,1}  P_2$& 3 &
$\R  \HH^2_2$ & $\C_2^3$ &$2/\sqrt{3}$ &5,7\\
& &8&$U_2  J_2$ & $V_{1,1}^{\rm{II}}  P_2$& 3 &
$\R  \HH^2_2$ & $\C_2^3$ &1 &5,8\\
&$R_2  R_2$&9& $F_4$ &$(F_4)$ & 3 &  $\R  \HH^2_2$
& $\C_2^3$ & $\sqrt{2}$&5,9 \\
\hline
3& $I_1 R_3$ &10& $I_1^-  F_3$ & $(I_1^-  F_3)$ &
1 & $\R$ & $\D_{12}$& 4/3 & 2,3,10\\
& & 11&$G_4$ & $\{\La_4\}$ &0 & $\R^0$ & 144 & 4/3 & 2-5,9-12,22\\
&$V_2  V_2$ &12& $H_4$ & $(H_4)$&2&$\HH^2_4$ & $\D_{12}$ &4/3&4,5,12\\
\hline
4 & $J_2  J_2$ & 13& $K_4$ & $(K_4)$ &2&$\HH^2_4$
& $\D_{8}$& $\sqrt{2}$ & 5,13\\
 & & 14& $L_4$ & $\{D_4\}$ &0 & $\R^0$ & 1152 & $\sqrt{2}$ &
3,5,9,13,14,16,22\\
 & $I_1  X_3$ & 15&$I_1  G_3$ & $\{I_1 \La_3\}$ &
0 &$\R^0$ & $\C_2^2 \D_8$ & 1 & 2,3,15,16\\
 & & 16&$I_1^-  G_3$ & $(I_1^-  G_3)$&1 &$\R$ &$\C_2 \D_8$ & 
$\sqrt{2}$ & 3,16\\
\hline
6 & $I_2  W_2$ & 17&$I_2  F_2$ & $\{I_2  A_2\}$&
0 & $ \R^0$ & $\D_8  \D_{12}$& 1 & 2,3,5-8,17-20 \\
& & 18&$I_{1,1}  F_2$ &  $V_{1,1}  \{A_2\} $ &
1 &$\R$ & $\C_2^2  \D_{12}$& $2/\sqrt{3}$ & 3,5,7,18,20\\
& & 19&$U_2  F_2$ &  $V_{1,1}^{\rm{II}}  \{A_2\} $ &
1 &$\R$ & $\C_2^2  \D_{12}$& 1 & 3,5,7,8,19,20\\
 & $I^-_2  W_2$ & 20&$J_2  F_2$ &  $P_2  \{A_2\} $ &
2 & $\HH^2_2$ &  $\C_2  \D_{12}$& $2/\sqrt{3}$ & 5,7,20\\
 & $W_2  W_2$ &  21&$F_2  F_2$ &  $\{A_2  A_2\} $ &
0 & $\R^0$ & 288 & $2/\sqrt{3}$ & 3-5,7,9,12,18,20-22\\
& & 22 & $F_2^-  F_2$ &  $(F_2^-  F_2)$ &
2 & $\HH^2_4$ & $\D_{12}$& $\sqrt{2}$ & 5,22\\
\hline
10 & $W_4$ &23& $M_4$ & $\{\La'_4\}$ & 0 & $\R^0$& $\D_{20}$ 
& $\be$ & 3,5,23\\
\hline
\end{tabular}}
\caption{Types géométriques isoduaux de rang   $4$
(23 types, 4 max.)}
\label{tabl:types_geom_4}
\end{table}

Nous déterminons les constantes d'Hermite relatives  en 
commençant par les types maximaux. Pour les types $\mfS_2$ et
$V_{2,2}$,  la constante $\mu_F$ vaut $\sqrt{2}$ car ces types 
contiennent $D_4$
à équivalence près : en posant $X=\smallmat{1}{-1}{1}{1}$ on a 
$$ 
D4 \simeq \sig_2\left( X/2  + i I_2/\sqrt{2}\right) \in \mfS_2
\esp  \mathrm{et}\esp 
D4 = \ph_{2,2}\left((1-1/\sqrt{2}) X\right) \in V_{2,2}
$$
(on rappelle que les plongements $\sig_g$ et $\ph_{p,q}$ sont
définis par (\ref{equa:Sig_g}) et (\ref{equa:Vpq}) respectivement). 
Le cas du type lorentzien $V_{3,1}$ est traité dans 
\cite[\S 1.3]{Bavard2007a} : on a $\mu_{3,1}=4/3$ et $V_{3,1}$
contient un point équivalent à $\La_4$. 
Il reste à examiner le type symétrique
pair $V_{2,2}^{\rm II}$. 
\begin{prop}\label{prop:V22pair}
Le maximum de l'invariant d'Hermite 
$\mu$ sur $V_{2,2}^{\rm II}$
vaut $\mu_{2,2}^{\rm II} =  4/3$, réalisé par une unique forme à
équivalence près dans $V_{2,2}^{\rm II}$.
\end{prop}
\pr Nous utiliserons le modèle ci-dessus. 
La formule (\ref{equa:SO22}) montre que l'action diagonale de 
$\SL_2(\Z)\times \SL_2(\Z)$ sur $\mfH_1\times \mfH_1$ correspond
à une action par équivalence entière sur le type $V_{2,2}^{\rm II}$.
Il est clair que le carré de 
\begin{equation}\label{equa:dom_fond_SL2}
\Delta=\{z\in\mfH_1; |\Re z| \leq 1/2, |z|\geq 1\}
\end{equation}
est un domaine fondamental pour $\SL_2(\Z)\times
\SL_2(\Z)$. Si $(z,w)\in \Delta\times \Delta$, on~a 
$(\Im z \Im w)^{-1}\leq 4/3$ avec égalité pour 
$\Im z= \Im w=\sqrt{3}/2$ ; comme ces derniers points correspondent
 à des matrices équivalentes à $\La_4$ (donc de minimum $4/3$),
on a finalement  $\mu_{2,2}^{\rm II} =  4/3$, voir 
(\ref{equa:V22_pair}).  
\par
Pour la détermination des autres constantes $\mu_F$, nous
pouvons écarter de la discussion les types décomposés
(diagonaux par blocs),  
dont les constantes $\mu_F$ sont données par celles des types 
de rang plus petit~: $(I_4)$,
$(I_2J_2)$, $(I_{1,1}J_2)$,  $(U_2J_2)$, $(I_1G_3)$, $(I_2F_2)$,
$(I_{1,1}F_2)$, $(U_2F_2)$, $(J_2F_2)$ et $(F_2F_2)$.
\par
Ensuite, il y a les types inclus (à équivalence près) à $\mfS_2$ ou 
à  $V_{2,2}$ et qui contiennent un point équivalent à $D_4$, pour
lesquels $\mu_F=\sqrt{2}$. Ces types 
correspondent aux chaînes 
$$(L_4) \subsetsim  (F) \subsetsim \mfS_2
\esp \mathrm{ou}\esp
 (L_4) \subsetsim (F)\subsetsim  V_{2,2}$$
que l'on trouve dans la dernière colonne de la 
table~\ref{tabl:types_geom_4}. Il s'agit de 
$(F_4)$, $(K_4)$, $(L_4)=\{D_4\}$, $(I_1^-G_3)$ et $(F_2^- F_2)$.
Les inclusions peuvent être explicitées dans chaque cas. Par 
exemple $(F_4)$ est inclus dans $\mfS_2$ car 
$$(F_4)=\left\{\sig_2\left(U_2(iu-z)/2+I_2(iu+z)/2\right); 
z \in\mfH_1, u >0\right\},$$
et la forme obtenue pour $z=1+i/\sqrt{2}, u=1/\sqrt{2}$ est 
équivalente à $D_4$. 
\par 
De même, pour les types $(F)$ vérifiant  $(G_4) \subsetsim (F)
\subsetsim V_{3,1}$
ou $(G_4)\subsetsim (F) \subsetsim V_{2,2}^{\rm II}$,
on a  $\mu_F=4/3$. Il s'agit des  types $(I_1^-F_3)$,
$(G_4)=\{\La_4\}$ et $(H_4)$ (voir table~\ref{tabl:types_geom_4}).
\par
Il reste un seul type à examiner, $(M_4)=\{\La'_4\}$ avec 
$${\small
\La'_4 = \frac{\sqrt{5+2\sqrt{5}}}{5}
\left( \begin{array}{cccc}
2& 1-\sqrt{5}& 1& \sqrt{5}-3
\\1-\sqrt{5}& 2\sqrt{5}-2& 3-2\sqrt{5}& 1\\
1& 3-2\sqrt{5}& 2\sqrt{5} -2 & 1-\sqrt{5}\\
\sqrt{5} - 3 & 1& 1-\sqrt{5}& 2
\end{array} \right),
}
$$
forme  de minimum $\be=\frac{2}{5}\sqrt{5+2\sqrt{5}} \simeq 1,2310.$

\subsection{Types géométriques principaux de rang $5$}
\label{subs:geom_5}
En rang $5$, nous trouvons 14 types géométriques principaux dont
5 maximaux : $I_{4,1},I_{3,2},I_{2,1}\oplus J_2, I_1\oplus J_4$ et
$I_1\oplus F_4$ ; nous déterminons les constantes $\mu_F$ à l'exception
de $\mu_{3,2}$ (voir table~\ref{tabl:types_geom_5}). 
Il existe 22 types  algébriques principaux de rang 5 au signe
près (table~\ref{tabl:types_5}), parmi lesquels on constate~8 
relations :
\renewcommand{\arraystretch}{1.1}
$$
\begin{array}{l}
(I_1 \oplus -L_4)=(I_1 \oplus L_4)=\{I_1\oplus D_4\},\\
(I_2 \oplus H_3)=(I_2 \oplus G_3)=\{I_2 \oplus \La_3\},~~
(-I_2 \oplus H_3)=(-I_2 \oplus G_3),\\
(U_2 \oplus H_3)=(U_2 \oplus G_3), ~~ (J_2 \oplus H_3)=(J_2 \oplus G_3),\\
(G_5) = (F_5) = \{\La'_5\},~~(H_5) \simeq (G_5),
~~ (K_5) = (H_5),
\end{array}
$$
où la forme $\La'_5$ est explicitée ci-dessous. Finalement, 
il y a 14 types géométriques principaux qui se différencient 
par leur dimension, leur structure métrique, la valeur de 
$\mu_F$ et le cardinal de $\Ga_F$, à l'exception des types
$(-I_2\oplus G_3)$ et \mbox{$(U_2\oplus G_3)$} distinguables 
à l'aide du critère d'inclusion (voir table~\ref{tabl:types_geom_5}).
On compte 8 types géométriques décomposés (pour lesquels on 
détermine facilement $\mu_F$) et 6 seulement non décomposés :
$V_{4,1}$, $V_{3,2}$, $(I_1\oplus F_4)$, $(-I_2\oplus G_3)$, 
$(U_2\oplus G_3)$ et $(F_5)$. Ce dernier type est réduit à un 
point $(F_5)= \{\La'_5\}$, défini par 
\begin{equation}\label{equa:cte_gamma}
{\small
\La'_5 = \frac{1}{2\ga}
\begin{pmatrix}
2+\ga-\ga^2 & - 1 & \ga^2-1 & \ga^2-1 & -1\\	
-1 & 2 & -\ga^2 & 0 & \ga^2\\	
\ga^2-1 & -\ga^2 & 2 & -\ga^2 & 0\\	
\ga^2-1 & 0 & -\ga^2 & 2 & -\ga^2\\	
-1 & \ga^2 & 0 & -\ga^2 & 2
\end{pmatrix}
}
~~~\mathrm{o\grave{u}}~~~
\gamma=(2-\sqrt{2})^{1/2}.
\end{equation}
Le minimum de $\La'_5$  vaut
$1/\ga=(1+1/\sqrt{2})^{1/2}\simeq 1,3065$.
\par
Passons à l'étude des constantes $\mu_F$, déjà déterminées pour
les types décomposés (voir table~\ref{tabl:types_geom_5})
et pour $(F_5)$. Il reste à examiner les types $V_{4,1}$, $V_{3,2}$,
$(I_1\oplus F_4)$, $(-I_2 \oplus G_3)$ et $(U_2\oplus G_3)$.
D'après \cite[théorème 1]{Bavard2007a}, on a $\mu_{4,1}=7/5$
réalisée de façon unique à équivalence  près dans 
$V_{4,1}$ par la forme $\La_5$ ci-dessous. Nous introduisons
également deux formes $\La_{5}^b$ et  $\La_{5}^c$ équivalentes
à $\La_5$ et appartenant respectivement à $V_{3,2}$ et 
$(I_1 \oplus F_4)$ :  
\begin{equation}\label{equa:L5_abc}
\renewcommand{\arraycolsep}{2pt}
\renewcommand{\arraystretch}{0.9}
{\small
\La_5 =\frac{1}{5}
\begin{pmatrix}
7 & 2 & 2 & 2 & 6\\
2 & 7 & 2 & 2 & 6\\
2 & 2 & 7 & 2 & 6\\
2 & 2 & 2 & 7 & 6\\
6 & 6 & 6 & 6 & 13
\end{pmatrix},
~
\La_{5}^b =\frac{1}{5}
\begin{pmatrix}
7 & 2 & 2 & 4 & 4 \\
2 & 7 & 2 & 4 & 4 \\
2 & 2 & 7 & 4 & 4 \\
4 & 4 & 4 & 8 & 3 \\
4 & 4 & 4 & 3 & 8
\end{pmatrix},
~ 
\La_{5}^c =\frac{1}{5}
\begin{pmatrix}
7 & -3 & -3 & 2 & 2 \\
-3 & 7 & 2 & -3 & 2 \\
-3 & 2 & 7 & 2 & -3 \\
2 & -3 & 2 & 7 & -3 \\
2 & 2 & -3 & -3 & 7
\end{pmatrix}.
}
\renewcommand{\arraycolsep}{3pt}
\renewcommand{\arraystretch}{1}
\end{equation}
\par

\begin{prop}\label{prop:maxI1F4}
Le maximum de $\mu$ sur le type géométrique $(I_1 \oplus F_4)$ vaut $7/5$,
atteint  de façon unique par $\La_5^c$ à équivalence près 
dans $(I_1\oplus F_4)$.
\end{prop}

\pr Pour $A=\smallmat{a}{b}{c}{d} \in
\SL_2(\R)$ et 
$B=\smallmat{\al}{\be}{\ga}{\de}\in\SL_2(\R)$ posons
\begin{equation}
\label{equa:theta_I1F4}
{\small
\theta_1(A,B)=\frac{1}{2}
\begin{pmatrix}
2 a d+2 b c & -2 a c & -2 a c & 2 b d & 2 b d\\
-2 a b & a^2+\alpha & a^2-\alpha & -b^2-\beta & -b^2+\beta\\
-2 a b & a^2-\alpha & a^2+\alpha & -b^2+\beta & -b^2-\beta\\
2 c d & -c^2-\gamma & -c^2+\gamma & d^2+\delta & d^2-\delta\\
2 c d & -c^2+\gamma & -c^2-\gamma & d^2-\delta & d^2+\delta
\end{pmatrix}.
}
\end{equation}
Cette application induit un isomorphisme entres les groupes 
$\PSL_2(\R)\times \SL_2(\R)$  et la composante neutre de
$\OO_{I_1F_4}^t(\R)$. Elle permet également d'écrire un paramétrage 
équivariant du  type $(I_1\oplus F_4)$ par le produit
 $\mfH_1\times \mfH_1$ sous la 
forme $A_{w,z}=P_{w,z}P'_{w,z}$ avec $w=s+it\in \mfH_1$, 
$z=x+iy \in\mfH_1$ et 
\begin{equation}
\label{equa:tg_I1F4}
{\small
P_{w,z}=\frac{1}{2t\sqrt{y}}
\begin{pmatrix}
2 t \sqrt{y} & 0 & 0 & 2 s \sqrt{y} & 2 s \sqrt{y}\\
-2 s t \sqrt{y} & t^2 \sqrt{y} + t y & t^2 \sqrt{y} - t y & 
-s^2 \sqrt{y}-x t & -s^2 \sqrt{y}+x t\\
-2 s t \sqrt{y} & t^2 \sqrt{y} - t y & t^2 \sqrt{y} + t y & 
-s^2 \sqrt{y}+x t & -s^2 \sqrt{y}-x t\\
0 & 0 & 0 & \sqrt{y}+t & \sqrt{y}-t\\
0 & 0 & 0 & \sqrt{y}-t & \sqrt{y}+t
\end{pmatrix}.
}
\end{equation}
Noter que $A_{i,i}=I_5$. D'après \eqref{equa:theta_I1F4}, le 
groupe $\{(A,B)\in\SL_2(\Z)^2;A\equiv B (\mathrm{mod}~2)\}$ agit
{\it via} $\theta_1$ par équivalence entière sur $(I_1\oplus F_4)$.
Soit $\Gamma(2)$ le noyau du morphisme de réduction modulo~2
de  $\SL_2(\Z)$ dans $\SL_2(\Z/2\Z)$. Le quadrilatère 
hyperbolique de sommets $-1,0,1,\infty$ est un domaine 
fondamental classique pour l'action de $\Gamma(2)$ sur $\mfH_1$.
Et faisant agir $(A,A)$ ($A\in\SL_2(\Z)$), puis 
$\{I_2\} \times \Gamma(2)$ et enfin les matrices
$I_{4,1}$ et $I_3\oplus U_2$ sur le type (qui changent $s$ en $-s$ et 
$x$ en $-x$, respectivement), on se ramène à 
$$
(w,z)\in\mcD =\{0\leq \Re w \leq \frac{1}{2}, |w|^2\geq 1\}
\times \{0\leq \Re z \leq 1, |z-\frac{1}{2}|^2\geq \frac{1}{4}\}.
$$

On rappelle la notation $A[u]=u'Au$ ($A \in P_n, u\in \R^n$).
Notons  $(\eps_i)_{1\leq i \leq 5}$ la base naturelle de $\Z^5$
et posons
\begin{equation}
\label{equa:lgr_I1F4}
\begin{array}{ll}
l_1=A_{w,z}[\eps_1]=\frac{2s^2+t^2}{t^2}, &
 l_2=A_{w,z}[\eps_2]=\frac{1}{2}(\frac{|w|^2}{\Im w})^2+
\frac{|z|^2}{2\Im z},\\
l_3=A_{w,z}[\eps_4]=\frac{1}{2(\Im w)^2}+\frac{1}{2\Im z}, &  
l_4=A_{w,z}[\eps_1+\eps_2-\eps_5]=
\frac{1}{2}(\frac{|w-1|^2}{\Im w})^2+\frac{|z-1|^2}{2\Im z}.
\end{array}
\end{equation}
Ces fonctions ont toutes une interprétation géométrique. Ainsi
$l_1=\cosh (2d)$ où~$d$ est la distance à la géodésique 
$\Re w=0$ dans le plan hyperbolique des $w$. Les autres 
$l_k$ ($k=2,3,4$) dépendent d'expressions  de la forme 
$|w-p|^2/\Im w$ ($p\in\Q)$, $1/\Im w$ et de quantités analogues
pour $z$.
Rappelons que $\log (|w-p|^2/\Im w)$ (resp. $\log (1/\Im w)$
est une fonction de Busemann du plan hyperbolique $\mfH_1$, mesurant
la \og distance \fg\ à la pointe $p\in \R$ 
(resp. à la pointe $\infty$); ses sous-niveaux sont les horoboules 
centrées en $p$ (resp. en $\infty$). Ainsi l'horoboule 
$\mcH^p(w_0)$ centrée
en $p\in \partial \mfH_1=\R\cup \{\infty\}$ dont le bord 
contient $w_0\in\mfH_1$ est donnée par  
\begin{equation}
\label{equa:horoboules}
\mcH^p(w_0)=\left\{\frac{|w-p|^2}{\Im w}\leq 
\frac{|w_0-p|^2}{\Im w_0}\right\}~~
(p\in \R) ~~\mathrm{ou}~~
\mcH^\infty(w_0)=\{\Im w \geq \Im w_0\}.
\end{equation}
Nous allons donc majorer $\mu_{I_1 F_4}$ en recouvrant le 
domaine $\mcD$ par une famille {\it ad hoc} de  produits d'ensembles
géométriques (bandes et horoboules). 
\par
\begin{figure}[h]
\labellist
\small\hair 2pt
\pinlabel $t$  at -10 190
\pinlabel $y$  at 280 190
\pinlabel $t_0$  at -10 142
\pinlabel $w_0$  at 85 142
\pinlabel $w_e$  at 85 160
\pinlabel $s$  at 172 10
\pinlabel $x$  at 462 10
\pinlabel $T$  at 57 107
\pinlabel $y_0$  at 280 112
\pinlabel $T'$  at 335 95
\pinlabel $\sigma_0$  at 355 122
\pinlabel $z_0$  at 390 105
\pinlabel $z_e$  at 345 57
\endlabellist
\centering
\includegraphics[scale=0.58]{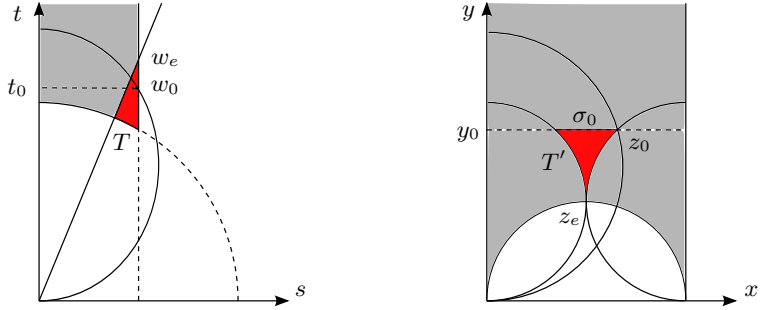}
\caption{Type géométrique $(I_1\oplus F_4)$}
\label{figu:I1F4}
\end{figure}

Soit $(w_e,z_e)\in \mcD$ défini par $w_e=\frac{1}{2}(1+i\sqrt{5})$ et 
$z_e= \frac{1}{2}(1+i)$ (voir figure~\ref{figu:I1F4}). On vérifie
que $\La_5^c=A_{w_e,z_e}$ et que le minimum de cette forme vaut $7/5$, 
valeur  commune des $(l_k)$ en $(w_e,z_e)$.
Si $w$ appartient à la bande $|s/t|\leq \sqrt{5}$, alors
$l_1(w,z)\leq 7/5$ pour tout $z\in \mfH_1$. On peut donc supposer
que $w\in T$ avec 
$$T=\{w\in\mfH_1;0\leq s\leq 1/2, |w|\geq 1 ~\mathrm{et}~
s \leq \sqrt{5}t\}.$$
Remarquer que $T$ est inclus dans les horoboules 
$\mcH^0(w_e)$ et $\mcH^1(w_e)$  (voir \eqref{equa:horoboules}).
Sur $T\times \mcH^0(z_e)\subset \mcH^0(w_e)\times \mcH^0(z_e)$
(resp. sur $T\times \mcH^1(z_e)\subset  \mcH^1(w_e)\times \mcH^1(z_e)$),
on~a  $l_2 \leq 7/5$ (resp. $l_4 \leq 7/5$). Pour 
$w\in T$ et $y\geq 15/22$, on a $l_3(w,z)\leq 7/5$ (car
$t\geq \sqrt{3}/2$). Il reste à examiner le cas de 
$(w,z)\in T\times T'$, où $T'$ est le triangle bordé
par les horosphères $\partial \mcH^1(z_e)$, $\partial \mcH^0(z_e)$
et $y =  15/22$. Pour conclure, on découpe $T'$ en tranches horizontales. 
Soit $y_0\in [1/2,15/22]$, soit $\sig_0=T'\cap \{y=y_0\}$ et
soit $t_0>0$ tel que $t_0^{-2}=14/5-y_0^{-1}$. On a évidemment
$l_3\leq 7/5$ sur $\sig_0\times (T\cap \{t\geq t_0\})$. Ensuite, on 
majore $l_2$ sur $\sig_0\times (T\cap \{t\leq t_0\})$ en utilisant 
le produit d'horoboules $\mcH^0(w_0)\times \mcH^0(z_0)$ 
($w_0=1/2+it_0$, $z_0$ extrémité droite de $\sig_0$,
voir figure~\ref{figu:I1F4}). Un calcul élémentaire conduit à 
$$l_2(w,z)\leq \frac{67}{80}+\frac{15}{32y_0}
+\frac{5y_0}{28y_0-10}-\left( \frac{1}{y_0}-1 \right)^{1/2}
~~ (z,w)\in \sig_0\times (T\cap \{t\geq t_0\}).
$$
Le second membre de cette inégalité étant une fonction strictement
décroissante de $y_0$ sur $[1/2,15/22]$, on a que 
$l_2\leq 7/5$ sur $\sig_0\times (T\cap \{t\leq t_0\})$, d'où
finalement $\mu_{I_1F_4}=7/5$. Concernant
le cas d'égalité, on observera que toutes les inégalités
$l_k\leq 7/5$ données plus haut sont strictes l'intérieur des domaines 
considérés ; pour les points du bord où $l_k=7/5$, on peut
toujours obtenir l'inégalité stricte  $\mu < 7/5$ en changeant 
de fonction $l_k$, à l'exception du  point $(w_e,z_e)$. 
Par exemple $l_3=7/5$ sur  $\sig_0\times (T\cap \{t = t_0\})$ 
mais d'après ce qui précède $l_2< 7/5$ sur cet ensemble,
sauf pour $y_0=1/2$. Ce qui achève la preuve.

\begin{prop}\label{prop:max-I2G3}
Le maximum de $\mu$ sur le type géométrique $(-I_2 \oplus G_3)$ 
vaut $4/3$, réalisé  de façon unique à équivalence près 
dans $(-I_2 \oplus G_3)$.
\end{prop}

\pr Soit $M=\smallmat{a}{b}{c}{d} \in \SL_2(\R)$. On pose
$q(a,b,c,d)=a^2+b^2-c^2-d^2$. Alors l'application 
$\theta_2:\SL_2(\R)\to \SL_5(\R)$ définie par
$$
\renewcommand{\arraystretch}{1.3}
\renewcommand{\arraycolsep}{2pt}
{\small
\theta_2(M)=
\begin{pmatrix}
ad +bc & ac -bd & \sqrt{2}(bd+ac) &  -\frac{\sqrt{2}}{2}(bd+ac) & 
\frac{\sqrt{2}}{2}(bd+ac) \\
ab-cd &  \frac{1}{2}q(a,d,b,c) &
\frac{\sqrt{2}}{2}q(a,b,c,d) & 
-\frac{\sqrt{2}}{4}q(a,b,c,d) & 
\frac{\sqrt{2}}{4}q(a,b,c,d) \\ 
\frac{1}{2}(ab+cd) & 
\frac{\sqrt{2}}{4}q(a,c,b,d) & 
\frac{\tr (MM')}{2} &  
\frac{1}{2}- \frac{\tr (MM')}{4} & %
\frac{\tr (MM')}{4} -\frac{1}{2} \\
0 & 0 & 0 & 1 & 0\\
0 & 0 & 0 & 0 & 1
\end{pmatrix}
}
\renewcommand{\arraycolsep}{3pt}
$$
induit un isomorphisme de $\PSL_2(\R)$ sur la composante neutre
du groupe $\OO^t_F(\R)$, $F=-I_2\oplus G_3$. Posons
\begin{equation}
\label{equa:Mz}
{\small
M_z=\frac{1}{\sqrt{\Im z}}
\begin{pmatrix}\Im z  & \Re z \\ 0 & 1 \end{pmatrix}
\esp (z\in\mfH_1).
}
\end{equation}
 La formule $A_z=\theta_2(M_z)\cdot (I_2\oplus \La_3)$   donne un 
paramétrage équivariant avec $A_i=I_2\oplus \La_3$
(voir~\eqref{equa:L3}  pour la définition de $\La_3$).
Pour $\al=u+v/\sqrt{2}\in \Z[1/\sqrt{2}]$, on pose
$\al^c=u-v/\sqrt{2}$. 
Soit $\Gamma$ le sous-groupe de $\SL_2(\R)$ formé des matrices
$$
\begin{pmatrix}
\al & -\eps \be^c \\
\be & \eps \al^c
\end{pmatrix}
\in \SL_2(\Z[1/\sqrt{2}])
~~ \mathrm{avec}~~
\left\{
\begin{array}{l}
\al\al^c+\be\be^c=\eps=\pm 1,\\
\al-\al^c+\be-\be^c \in 2\sqrt{2}\Z .
\end{array}
\right.
$$
Une vérification facile montre  que $\theta_2(\Ga)$ est inclus dans 
$\SL_5(\Z)$. En  particulier~$\Ga$ est discret et ses 
orbites sur $\mfH_1$ ne s'accumulent pas -- pour la topologie des
pointes -- sur les pointes du groupe (points fixes des paraboliques 
de $\Ga$). Autrement dit, si $p\in \partial \mfH_1$ est 
une pointe de $\Ga$, alors toute 
orbite de $\Ga$ sur $\mfH_1$ contient un point $z$ à 
\og distance minimale\fg\ de $p$, c'est-à-dire tel que 
\begin{equation}\label{equa:dist_mini_pointe}
|z-p|^2 \leq \frac{\Im z}{\Im \ga(z)} |\ga(z)-p|^2
=  |cz+d|^2 |\ga(z)-p|^2
\end{equation}
pour tout $\ga=\smallmat{a}{b}{c}{d}\in \Ga$
(voir~\eqref{equa:horoboules}). De plus, si $t\in \Ga$ est un
parabolique fixant~$p$, tous les $t^k(z)$ ($k\in \Z$) vérifient 
encore la condition~\eqref{equa:dist_mini_pointe}. 
Considérons les transformations suivantes :
$$
t = \frac{1}{2}
\begin{pmatrix} 
2-\sqrt{2} &  2+\sqrt{2} \\
-2 +\sqrt{2} & 2+\sqrt{2}
\end{pmatrix}\in \Ga
\esp \mathrm{et} \esp
\ga= 
\begin{pmatrix} 
5+ 2 \sqrt{2} &  - 3 \sqrt{2} \\
3 \sqrt{2} & -5 +2 \sqrt{2}
\end{pmatrix}
\in \Ga.
$$

\begin{figure}[b]
\labellist
\small\hair 2pt
\pinlabel $0$  at 128 -10
\pinlabel $1$  at 225 -10
\pinlabel $1\!-\!\sqrt{2}$  at 88 -10
\pinlabel $1\!+\!\sqrt{2}$  at 385 -10
\pinlabel $z_1$  at 252 62
\pinlabel $\mcE$  at 252 110
\pinlabel $z_0$  at 210 157
\endlabellist
\centering
\includegraphics[scale=0.58]{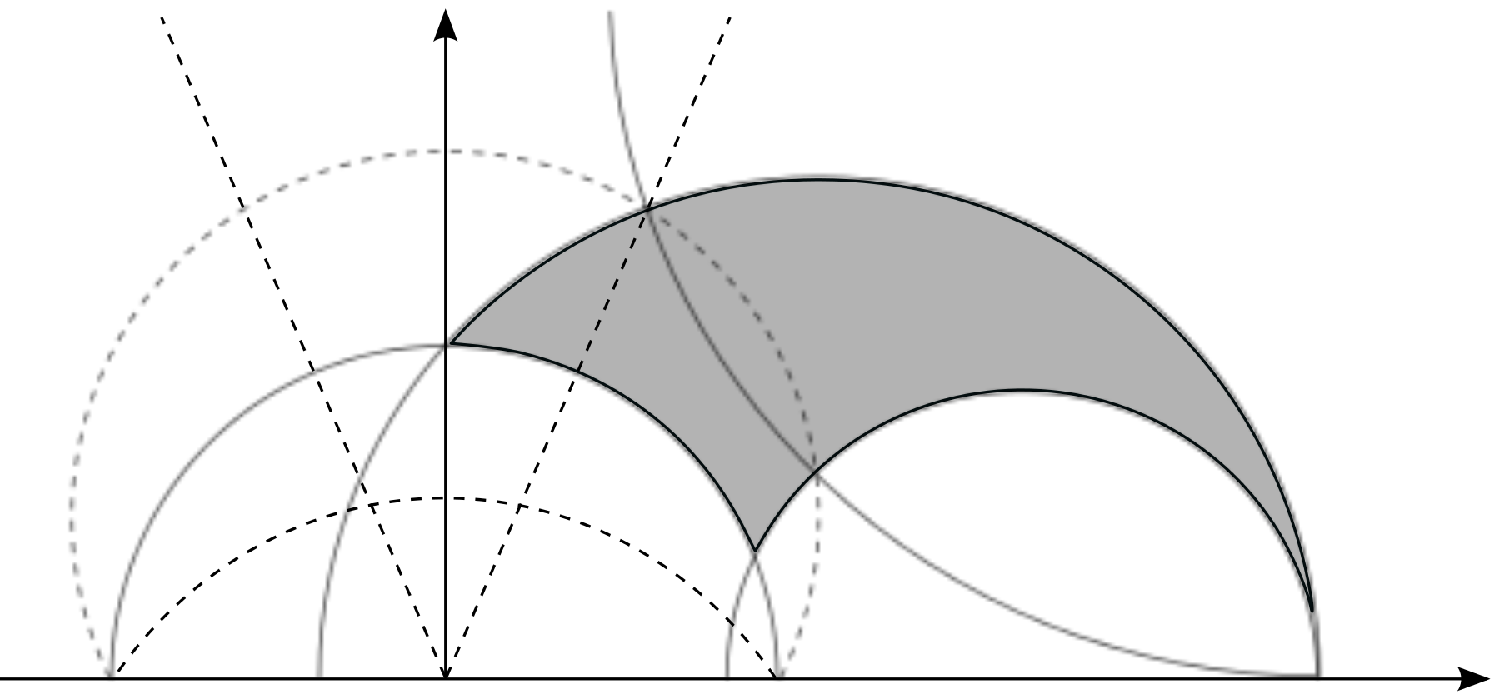}
\vspace{5pt}
\caption{Type géométrique $(-I_2 \oplus G_3)$}
\label{figu:I2mG3}
\end{figure}

L'élément $t\in \Ga$ est parabolique de point fixe $p=1+\sqrt{2}$ et
la condition~\eqref{equa:dist_mini_pointe} pour $p$ et $\ga$ se
réduit à $|z|^2\geq 1$. D'après ce qui précède, toute orbite de
$\Ga$ coupe le domaine $\mcE$ défini par les inégalités 
(voir figure~\ref{figu:I2mG3}) 
$$
|z|^2\geq 1,~~|z-1|^2\leq 2 ~~\mathrm{et}~~ 
|7z-3(1+2\sqrt{2})|^2 \geq (4+\sqrt{2})^2.
$$
Les fonctions  $l_1(z)=A_z[\eps_1+\eps_2-2\eps_3+\eps_4-\eps_5]$ et 
$l_2(z)=A_z[\eps_2]$ ($z\in\mfH_1$) s'explicitent en
$$l_1(z)=\frac{(\sqrt{2}-1)^2}{2}
\left(\frac{|z-(1+\sqrt{2})|^2}{\Im z}\right)^2
~~~ \mathrm{et}~~~
l_2(z)=1+\frac{1}{2}\left(\frac{|z|^2-1}{\Im z}\right)^2.
$$
D'après \eqref{equa:horoboules}, les sous-niveaux de $l_1$ sont les
horoboules centrées en $p=1+\sqrt{2}$. De plus, si $d$ désigne la 
distance à la géodésique $|z|^2-1=0$, alors $l_2=\cosh(2d)$. 
Par suite on a $\min(l_1,l_2)\leq 4/3$ sur $\mcE$ avec 
égalité seulement pour les points 
$z_0=\frac{1}{7}[1+2\sqrt{2}+i\sqrt{3}(4+\sqrt{2})]$ et
$z_1= t(z_0)$ (voir figure\ref{figu:I2mG3}).  On vérifie enfin que 
le minimum de la forme $A_{z_0}$ vaut $4/3$, d'où le résultat. 

\begin{table}[t]
\centering
\renewcommand{\arraystretch}{1.2}
\renewcommand{\tabcolsep}{0.3 em}
{\small 
\begin{tabular}{ccccccclcl}
\hline
 $d$  & $F F^\vee$ &$n^\circ$   & $F$ & $(F)$& dim. & t. s. &
$\Ga_F$ & $\mu_F$ &$\subset_{\rm eq}$ \\
\hline
1&$I_5$ &1&$I_5$&$\{I_5\}$ &  0 & $\R^0$ &$\G(5)$ &1 & 1-8,0\\
&&2&$I_{4,1}$ & $V_{4,1}$ & 4 & $\HH^4_8$ & $\C_2$& 7/5& max.\\
&&3&$I_{3,2}$ & $V_{3,2}$ & 6 & $V_{3,2}$  & $\C_2$& $\geq 7/5$& max.\\
\hline
2&$I_{3,2}$&4&$I_3J_2$ & $\{I_3\}  P_2$ & 2  & $\HH^2_2$ 
& $\G(3) \C_2$& 1 & 4-6,0\\
&&5&$I_{2,1}J_2$ & $V_{2,1}  P_2$ & 4  & $\HH^2_8  \HH^2_2$ 
& $\C_2^2$& $2/\sqrt{3}$& max.\\
&$I_{1,4}$&6&$I_1  J_4$ & $\{I_1\}  \mfS_2$& 6 & $\mfH_2$ 
&$\C_2^2$&1 &max.\\
&$I_1R_2R_2$&7&$I_1  F_4$ & $(I_1  F_4)$& 4 & 
$\HH^2_8  \HH^2_2$ &$\C_2^2$& 7/5 &max.\\
\hline
4& $I_1J_2J_2$ &8& $I_1 K_4$ & $\{I_1\}  (K4)$ &
2 & $\HH^2_4$ & $\C_2  \D_{8}$& 1 & 6,8\\
& & 9&$I_1L_4$&$\{I_1D_4\}$ &0 & $\R^0$ & $\C_2\Ga_{L_4}$
& 1 & 3,6-9,11,0\\
&$I_2X_3$ &10& $I_2G_3$ & $\{I_2 \La_3\}$ &0 & $\R^0$ &
$\C_2\D_8^2$ & 1 & 2,3,5,10-13\\
& &11& $I_2^-G_3$ & $(I_2^-G_3)$ &2 & $\HH^2_8$ &
$\C_2\D_8$ & 4/3 & 3,11\\
& &12& $U_2G_3$ & $(U_2G_3)$ &2 & $\HH^2_8$ &
$\C_2\D_8$ & 4/3 &3,12\\
& $I_2^- X_3$ & 13& $J_2G_3$ & $P_2 \{\La_3\}$ &2&$\HH^2_2$&
$\C_2^2\D_8$ & $2/\sqrt{3}$ & 5,13\\
\hline
8 & $X_5$ & 14& $F_5$ & $\{\La'_5\}$ &0&$\R^0$
& $\C_2  \D_{16}$& $1/\ga$ & 3,14\\
\hline
\end{tabular}}
\caption{Types géométriques principaux de rang $5$
(14 types, 5 max.)}
\label{tabl:types_geom_5}
\end{table}

\begin{prop}\label{prop:maxU2G3}
Le maximum de $\mu$ sur le type géométrique $(U_2 \oplus G_3)$ 
vaut $4/3$, réalisé  de façon unique à équivalence près 
dans $(U_2 \oplus G_3)$.
\end{prop}

\pr Pour $M=\smallmat{a}{b}{c}{d} \in \SL_2(\R)$ on  pose
$$
\renewcommand{\arraystretch}{1}
{\small
\theta_3(M)=
\begin{pmatrix}
d^2 & -c^2 & 2 c d & -c d & c d \\
-b^2 & a^2 & -2 a b & a b & -a b \\
b d & -a c & a d+b c &\frac{1}{2}(1-ad-bc) & \frac{1}{2}(ad+bc -1)\\
0 & 0 & 0 & 1 & 0 \\
0 & 0 & 0 & 0 & 1
\end{pmatrix}.
}
$$
Cette application induit un isomorphisme de $\PSL_2(\R)$ 
sur $\OO^t_{U_2G_3}(\R)_0$, ainsi qu'un paramétrage équivariant
du type  $(U_2\oplus G_3)$ par $B_z=\theta_3(M_z)\cdot (I_2\oplus
\La_3)$ (\mbox{$z\in\mfH_1$}, $M_z$ est définie  en \eqref{equa:Mz}).
On voit immédiatement que $\SL_2(\Z)$ agit {\it via} $\theta_3$ 
par équivalence entière sur $(U_2\oplus G_3)$. On peut donc 
supposer que \mbox{$z\in \Delta$} (domaine fondamental standard pour
$\SL_2(\Z)$, voir~\eqref{equa:dom_fond_SL2}). 
Noter également que $I_{3,2}$ agit par
équivalence sur  $B_z$ en changeant le signe 
de $\Re z$. La forme $B_z$ est facile~à expliciter, en particulier
on trouve $B_z[\eps_1]=1/(\Im z)^2$. Par suite on~a $\mu\leq 4/3$
sur~$\De$, avec égalité seulement pour 
$z=\frac{1}{2}(\pm 1 + \sqrt{3})$. Les formes correspondantes sont
de minimum $4/3$ (et équivalentes à la forme $A_{z_0}$ ci-dessus).

\medskip
\begin{rema}
Les formes maximisant $\mu$ sur   $(-I_2\oplus G_3)$
et $(U_2\oplus G_3)$ sont toutes équivalentes. Les deux types se 
coupent suivant une géodésique commune, donnée respectivement 
par $|z+1|^2=2$ et par $|z|^2=1$ dans les paramétrages  ci-dessus
($A_z$ et $B_z$);
ils ne sont pas équivalents (critère d'inclusion, voir 
table~\ref{tabl:types_geom_5}).
\end{rema}

\begin{coro}
\label{coro:ga5isod}
La constante d'Hermite isoduale $\ga_5^\mathrm{isod}$ est atteinte
sur le type orthogonal $V_{3,2}$ :
$$\ga_5^\mathrm{isod} = \mu_{3,2}\geq \frac{7}{5}.$$
\end{coro}

\begin{rema}
La forme $\La_5^b\in V_{3,2}$ donnée en~\eqref{equa:L5_abc} est de 
minimum $7/5$. On peut vérifier que c'est un maximum local de 
$\mu$ sur le type $V_{3,2}$ (car parfaite et eutactique relativement
à $V_{3,2}$, voir par exemple \cite[proposition~2.1-(1)]{Bavard1997a}).
\end{rema}

\subsection{Types géométriques principaux de rang $6$}
En rang 6, il existe 88 types algébriques
(table~\ref{tabl:types_6}), dont 16 vérifient la 
relation $[F]=[-F]$. Il reste donc 52 types algébriques 
à considérer, parmi lesquels nous trouvons 12~relations :
\renewcommand{\arraystretch}{1.1}
$$
\begin{array}{ll}
(I_2 \oplus -L_4)=(I_2 \oplus L_4), 
& (I_{2,1} \oplus G_3)=(-I_3 \oplus G_3),\\
(I_3 \oplus H_3)=(I_3 \oplus G_3),&
(J_2\oplus I_1 \oplus H_3)=(J_2\oplus I_1 \oplus G_3),\\
(H_3\oplus H_3) = (G_3 \oplus G_3), & 
 (I_1 \oplus H_5) = (I_1 \oplus K_5) \simeq
 (I_1 \oplus F_5) = (I_1 \oplus G_5) ,\\
(-I_1 \oplus H_5) = (-I_1 \oplus F_5), &
(L_6) = (K_6) \simeq (M_6) = (N_6),
\end{array}
$$
qui conduisent à une liste de 40 types géométriques principaux
(table~\ref{tabl:types_geom_6}). Cela fait apparaître 11 types maximaux : 
les types symétriques indéfinis ou symplectique (5 types),
$F\oplus J_2$ pour $F=I_{3,1},I_{2,2},U_4$ et 
$F\oplus F_4$ pour $F=I_2,I_{1,1},U_2$.
\par
Les 40 types  se distinguent généralement par leur dimension, 
leur structure métrique, la valeur de $\mu_F$ et le cardinal du 
groupe $\Ga_F$ (voir table~\ref{tabl:types_geom_6}). 
Pour les quelques exceptions restantes on utilise 
le critère d'inclusion (qui de toutes façons permet de trier
tous les types) : $(I_{3,3})$ et $(U_6)$ 
(proposition~\ref{prop:orth_symp_max}, (1)), $(I_{3,1}\oplus J_2)$ et 
$(I_2\oplus F_4)$ (proposition~\ref{prop:aut_ordre_4}) et 
$(I_{2,2}\oplus J_2)$, $(U_4\oplus J_2)$, $(I_{1,1}\oplus F_4)$,
$(U_2\oplus F_4)$ (proposition~\ref{prop:aut_ordre_4}).  

\par
Il existe un réseau isodual  remarquable de rang $6$ provenant
de la théorie des surfaces de Riemann. On sait que les surfaces
de Riemann compactes de genre~$g$ possédant un  automorphisme 
d'ordre $4g$ sont hyperelliptiques, à l'exception d'une surface
de genre $3$ dite {\em surface (ou courbe) de Wiman exceptionnelle},
ou encore {\em courbe de Picard}. La jacobienne de cette courbe est 
un réseau symplectique étudié dans 
\cite{Bavard1997a,Quin-Zhan1998a}. On peut l'expliciter comme suit :
\begin{equation}\label{equa:cte_delta}
{\small
W_6= 
 \begin{pmatrix}
\delta & 0 & \delta-2 & 1-\delta/2 & -1 & -\delta/2 \\ 
0 & 2 & 0 & -1 & 1 & 1 \\ 
\delta-2 & 0 & \delta & -\delta/2 & 1 & 1-\delta/2 \\ 
1-\delta/2 & -1 & -\delta/2 & \delta & -1 & \delta-2 \\ 
-1 & 1 & 1 & -1 & 2 & 1 \\ 
-\delta/2 & 1 & 1-\delta/2 & \delta-2 & 1 & \delta
\end{pmatrix}
\in \mfS_3,
~~~  \mathrm{o\grave{u}}~~\delta=1+\frac{1}{\sqrt{3}}.
}
\end{equation}

Le minimum de $W_6$ vaut $\de\simeq 1,5773$, qui semble 
correspondre à la plus grande densité connue pour un réseau
isodual de rang $6$.
On trouve des formes équivalentes à $W_6$ dans plusieurs types
géométriques (ceux pour lesquels $\mu_F\geq \de$, voir 
table~\ref{tabl:types_geom_6}).
\par
La table~\ref{tabl:types_geom_6} ne mentionne que l'inclusion des
types dans les types maximaux (dernière colonne). L'ordre entre
les types principaux est donné par la table~\ref{tabl:incl_tg6} ;
la valeur  $0$ signifie  que le type  est inclus dans un type
non principal. Les constantes $\alpha,\gamma,\delta$ sont définies 
aux équations  
\eqref{equa:cte_alpha},\eqref{equa:cte_gamma},\eqref{equa:cte_delta}
et les constantes $\epsilon,\zeta,\eta$ par
$$
{\textstyle
\eps = 2 \sqrt{\frac{5+4\sqrt{3}}{23}} \simeq 1,4403, ~
\zeta = \frac{1+\sqrt{13}}{3} \simeq 1,5351,~
\eta = \sqrt{1+\frac{2}{\sqrt{3}}} \simeq 1,4678.
}
$$
Pour compléter l'énoncé du théorème~\ref{theo:intro_param1_7}, noter
également que 
\begin{equation}
\label{equa:V33_pair}
{\small
V_{3,3}^{\rm II} = P\cdot V_{3,3} \esp \mathrm{avec} \esp
P=\frac{1}{\sqrt{2}} 
\begin{pmatrix} 
I_3 & -I_3 \\
I_3 & I_3 
\end{pmatrix}.
}
\end{equation}

\begin{table}[t]
\centering
\renewcommand\arraystretch{1.2}
\renewcommand\tabcolsep{2pt}
{\small \begin{tabular}{llrccccccl}
\hline
$d$  & $F F^\vee$& $n^o$ & $F$ & $V_F=(F)$& dim. & type sym.
& $\Ga_F$ & $\mu_F$ &$\subset_\simeq$\\
\hline
1&$I_6$ &1&$I_6$&$\{I_6\}$ &  0 & $\R^0$ &$\G(6)$ &1 & 2-6,8-10,14-16\\
&&2&$I_{5,1}$ & $V_{5,1}$ & 5 & $\HH^5_8$ & $\C_2$& 3/2& max.\\

&&3&$I_{4,2}$ & $V_{4,2}$ & 8 & $V_{4,2}$  & $\C_2$& $\geq \de$& max.\\
&&4&$I_{3,3}$ & $V_{3,3}$ & 9 & $V_{3,3}$  & $\C_2$& $\geq \de$& max.\\
&&5&$U_6$ & $V_{3,3}^{II}$ & 9 & $V_{3,3}$  & $\C_2$& $\geq 3/2 $& max.\\
\hline
2&$I_6^-$&6& $J_6$& $\mfS_3$ & 12 & $\mfH_3$ & $\C_2$ 
& $\geq \de$ & max.\\
& $I_{4,2}$&7&$I_4J_2$ & $\{I_4\}  P_2$ & 2  & $\HH^2_2$ 
& $\G(4) \C_2$& 1 & 6,8-10\\
&&8&$I_{3,1}J_2$ & $V_{3,1}  P_2$ & 5  & $\HH^3_8  \HH^2_2$ 
& $\C_2^2$& $2/\sqrt{3}$& max.\\
&&9&$I_{2,2}J_2$ & $V_{2,2}  P_2$ & 6  & 
$(\HH^2_4)^2  \HH^2_2$& $\C_2^2$& $2/\sqrt{3}$& max.\\
&&10&$U_4J_2$ & $V_{2,2}^{II}  P_2$ & 6  &
$(\HH^2_4)^2  \HH^2_2$ & $\C_2^2$& $2/\sqrt{3}$& max.\\
& $I_{2,4}$&11&$I_2J_4$ & $\{I_2\}  \mfS_2$ & 6  & $\mfH_2$ 
& $\D_8  \C_2$& 1 & 6\\
&&12&$I_{1,1}J_4$ & $V_{1,1}  \mfS_2$ & 7  & 
$\R  \mfH_2$& $\C_2^3$& $2/\sqrt{3}$& 6\\
&&13&$U_2J_4$ & $V_{1,1}^{II}  \mfS_2$ & 7  & 
$\R  \mfH_2$& $\C_2^3$& $1$& 6\\
&$I_2R_2R_2$&14&$I_2  F_4$ & $(I_2  F_4)$& 5 & 
$\HH^3_8  \HH^2_2$ &$\C_2^2$& $\geq \de$ & max.\\
& &15&$I_{1,1}  F_4$ & $(I_{1,1}  F_4)$& 6 & 
$(\HH^2_4)^2  \HH^2_2$ &$\C_2^2$& $\geq \zeta$ & max.\\
& &16&$U_2  F_4$ & $(U_2  F_4)$& 6 & 
$(\HH^2_4)^2  \HH^2_2$ &$\C_2^2$& $\geq \de$ & max.\\
& &17&$F_6$ & $(F_6)$& 2 & 
$\HH^2_2$ &384& $\de$ & 6,14,16\\
&$I_2^-R_2R_2$&18&$J_2  F_4$ & $(J_2  F_4)$& 7 & 
$\R  \mfH_2$ & $\C_2^3$& $\geq \de$ & 6 \\
\hline
4& $I_2J_2J_2$ &19& $I_2 K_4$ & $\{I_2\}  (K4)$ &
2 & $\HH^2_4$ & $\D_{8}  \D_{8}$& 1 & 6\\
 & &20& $I_{1,1} K_4$ & $ V_{1,1} (K4)$ &
3 & $\R  \HH^2_4$ & $\C_2^2  \D_{8}$& $2/\sqrt{3}$ & 6\\
 & &21& $U_2 K_4$ & $ V_{1,1}^{II} (K4)$ &
3 & $\R  \HH^2_4$ & $\C_2^2  \D_{8}$& 1 & 6\\
& &22& $I_2 L_4$ & $\{I_2 D_4\}$ &
0 & $\R^0$ & $\D_{8}  \Ga_{L_4}$& 1 & 3,4,6,9,14-16\\
 & &23& $I_{1,1} L_4$ & $ V_{1,1} \{D_4\}$ &
1 & $\R $ & $\C_2^2  \Ga_{L_4}$& $2/\sqrt{3}$ & 4,6,9,15\\
 & &24& $U_2 L_4$ & $ V_{1,1}^{II} \{D_4\}$ &
1 & $\R $ & $\C_2^2  \Ga_{L_4}$& 1 & 4,6,9,16\\
& $I_2^-J_2J_2$ &25& $J_2 K_4$ & $P_2  (K4)$ &
4 & $\HH^2_2  \HH^2_4$ & $\C_2  \D_{8}$& $2/\sqrt{3}$ & 6\\
& &26& $J_2 L_4$ & $P_2  \{D_4\}$ &
2 &  $\HH^2_2$ & $\C_2  \Ga_{L_4}$& $2/\sqrt{3}$ & 6,9\\
&$I_3X_3$ &27& $I_3G_3$ & $\{I_3 \La_3\}$ &0 & $\R^0$ &
$\G(3) \C_2 \D_8$ & 1 & 2-4,8,9\\
&&28& $I_3^-G_3$ & $(I_3^-G_3)$ &3 & $\HH^3_8$ &
$\C_2  \D_8$ & 7/5 & 3\\
&&29& $I_{1,2}G_3$ & $(I_{1,2}G_3)$ &4 & $(\HH^2_4)^2$ &
$\C_2  \D_8$& $\sqrt{2}$ & 4\\
&$I_{2,1}^-X_3$ &30& $J_2I_1G_3$ & $P_2  \{I_1 \La_3\}$ &2 &
$\HH^2_2$ & $\C_2^3 \D_8$ & 1 & 8,9\\
& &31& $J_2I_1^-G_3$ & $P_2  (I_1^-G_3)$ &3 &
$\HH^2_2  \R$ & $\C_2^2 \D_8$ & $2/\sqrt{3}$ & 9\\
&$X_3X_3$ &32& $G_3 G_3$ & $\{\La_3^2\} $ &0 & $\R^0$ & 512
& $\al$& 3-6,16\\
&&33& $G_3 H_3$ & $(G_3 H_3) $ &2 & $\HH^2_4$ & 32
& $\zeta$ & 6\\
&&34& $G_3 H_3^-$ & $(G_3 H_3^-) $ &1 & $\R$ & 256
& $\sqrt{2}$ & 4,6,16\\
&&35& $G_3 G_3^-$ & $(G_3 G_3^-) $ &3 & $\R  \HH^2_4$ & $\C_2\D_8$
& $\geq \de$ & 6\\
&&36& $G_6$ & $(G_6) $ &3 & $\R  \HH^2_4$ & $\C_2^2\D_8$
& $\geq \eps$ & 6\\
&$X_3^-X_3^-$ &37& $H_6$ & $(H_6)$ &4 & $\HH^2_2  \HH^2_4$ & $\C_2\D_8$
&$\geq \eta$  & 6\\
\hline
8 & $I_1X_5$ & 38& $I_1F_5$ & $\{I_1\La_5\}$ &0&$\R^0$
& $\C_2^2\D_{16}$& 1 & 3,4\\
& & 39&$I_1F_5^-$ & $(I_1F_5^-)$& 1& $\R$&$\C_2\D_{16}$& $1/\ga$&4\\
&$Y_6$&40&$K_6$&$\{\La_6\}$&0&$\R^0$&128&$1/\ga$&4\\
\hline
\end{tabular}}
\caption{Types géométriques principaux en dimension $6$ 
(40 types, 11 max)}
\label{tabl:types_geom_6}
\end{table}

\FloatBarrier

\begin{table}[b]
\centering
\renewcommand{\arraystretch}{1.2}
\renewcommand\tabcolsep{5pt}
{\scriptsize \begin{tabular}{rlrlrl}
\hline
1 & 0,1-16,18-21,25 & 15 & 15 &   29 & 4,29  \\
2 & 2 &    16 & 16            &   30 & 8,9,30,31 \\
3 & 3 &    17 & 0,6,14,16-18  &   31 & 9,31 \\
4 & 4 &    18 & 6,18          &   32 & 3-6,16,18,28,29,32-37 \\
5 & 5 &    19 & 6,11-13,19-21,25& 33 & 6,18,33,35\\
6 & 6 &    20 & 6,12,20,25    &   34 & 4,6,16,18,29,34,35,37\\
7 & 0,6-13,18,25 & 21 & 6,13,21,25 & 35 & 6,35 \\
8 & 8&     22 & 0,3,4,6,9,11-16,18-26,28,29,31 & 36 & 6,18,36,37\\ 
9 & 9&     23 &0,4,6,9,12,15,18,20,23,25,26,29,31 & 37 & 6,37\\
10& 10 &   24 & 0,4,6,9,13,16,18,21,24-26,29,31 & 38 & 3,4,38,39\\
11 & 6,11-13 & 25 & 6,25   &   39 & 4,39\\
12 & 6,12   & 26 & 0,6,9,18,25,26,31 & 40 & 4,29,39,40\\
13 & 6,13   &  27 & 0,2-4,8,9,27-31 & & \\
14 & 14 &  28 &3,28 & & \\
\hline
\end{tabular}}
\caption{Inclusions des types principaux en dimension $6$}
\label{tabl:incl_tg6}
\end{table}

Notre connaissance des constantes $\mu_F$ n'est que très
partielle; il reste en particulier à déterminer $\mu_F$ 
pour 7 des 11 types maximaux.  Nous nous 
contenterons de  quelques estimations, permettant notamment de
retrouver la non-équivalence de certains types. 
Pour les types paramétrés par l'espace hyperbolique $\HH^n$
($n$ \og petit\fg), on parvient assez facilement à trouver $\mu_F$
(voir \S\ref{subs:geom_5}). Nous déterminons ainsi $\mu_F$ pour 
tous les types $(F)$ de dimension $\leq 2$, ainsi que pour 
le type $(-I_3\oplus G3)$, paramétré par $\HH^3$ (en procédant 
comme dans \cite{Bavard2007a} pour ce dernier cas).

\begin{prop}\label{prop:maxF6}
Le maximum de $\mu$ sur le type $(F_6)$ 
vaut $\de=1+1/\sqrt{3}$, 
réalisé  de façon unique à équivalence près  dans $(F_6)$ par une 
forme équivalente à $W_6$.
\end{prop}

\pr On obtient un isomorphisme $\theta_4$ entre $\SL_2(\R)$ et la 
composante neutre du groupe $\OO_{F_6}^t$ en posant
$$
{\small
\theta_4(M)=I_2 \oplus \frac{1}{2}
\begin{pmatrix}
1+a & 1-a & b & -b \\ 
1-a & 1+a & -b & b \\ 
c & -c & 1+d & 1-d \\ 
-c & c & 1-d & 1+d
\end{pmatrix}
}
$$
pour $M=\smallmat{a}{b}{c}{d} \in \SL_2(\R)$. Ensuite, on paramètre 
le type comme d'habitude, par exemple en posant 
$A_z=\theta_4(M_z)\cdot A_i$ où  
$$
{\small
A_i=\frac{1}{2}
\begin{pmatrix}
4 & 2 & 0 & 0 & 0 & 0 \\ 
2 & 4 & 0 & 0 & 2 & 2 \\ 
0 & 0 & 3 & 1 & 1 & 1 \\ 
0 & 0 & 1 & 3 & 1 & 1 \\ 
0 & 2 & 1 & 1 & 3 & 1 \\ 
0 & 2 & 1 & 1 & 1 & 3
\end{pmatrix}
\in (F_6).
}
$$
Le sous-groupe de congruence principal $\Ga(2)$ agit {\it via}
$\theta_4$ par équivalence entière sur le type $(F_6)$. On peut 
donc supposer que le point $z\in \mfH_1$ appartient au quadrilatère 
hyperbolique de sommets $-1,0,1,\infty$ (domaine fondamental de 
$\Ga(2)$) et même que $\Re z \geq 0$ car l'action de $I_2\oplus
U_2\oplus I_2$ sur $(F_6)$ change $z$ en $-\overline{z}$. Cela
étant, on majore facilement $\mu$ en considérant 3 longueurs
$$
\begin{array}{ll}
l_0(z)=A_z[\eps_3]= 1+\frac{|z|^2}{2\Im z}, & 
l_1(z)=A_z[\eps_3-\eps_5]= 1+\frac{|z-1|^2}{2\Im z},\\
l_\infty(z)=A_z[\eps_5] =  1+\frac{1}{2\Im z}, &
\end{array}
$$
associées aux pointes $0,1$ et $\infty$ respectivement. Sur le 
triangle de sommets $0,1,\infty$, on a 
$\min(l_0(z),l_1(z),l_\infty(z))\leq \de$  
(voir~\eqref{equa:horoboules}), avec 
égalité \ssi $z=z_0= (1+i\sqrt{3})/2$. Enfin la forme $A_{z_0}$ est
équivalente à $W_6$.

\smallskip
\begin{rema}
\label{rema:F6}
À équivalence près, le type $(F_6)$ est inclus dans le type 
symplectique $\mfS_3$. On peut également vérifier que la forme
$A_i$ ci-dessus est une matrice de Gram de la jacobienne de 
la quartique de Fermat.
\end{rema}

\begin{prop}\label{prop:max-G3H3}
Le maximum de $\mu$ sur le type géométrique $(G_3\oplus H_3)$ 
vaut $\zeta=(1+\sqrt{13})/3 \simeq 1,5351$, réalisé  de 
façon unique à équivalence près  dans $(G_3\oplus H_3)$.
\end{prop}

\pr On a un morphisme injectif $\theta_5:\SL_2(\R)\to 
\OO^t_F(\R)$ ($F=G_3\oplus H_3$) en posant
pour $M=\smallmat{a}{b}{c}{d} \in \SL_2(\R)$
$$
{\small
\theta_5(M)=\frac{1}{2}
\begin{pmatrix}
2 & l(a,b,d;c)-1  & 1-l(a,c,d;b) & 0 & -l(b,c,d;a)& -l(a,b,c;d)\\
0 & a+ d & b-c  & 0 &  a-  d & - b-c\\
0 & c- b &   a+  d & 0 &  b+ c &  a- d\\
0 & l(a,b,c;d) & l(b,c,d;a) & 2 & l(a,c,d;b)-1 & 1-l(a,b,d;c)\\
0 &  a- d &  b+ c & 0 &  a+ d & c- b\\
0 & - b-c & a- d & 0 & b -c  &  a+  d
\end{pmatrix},
}
$$
avec $l(a,b,c;d)=(a+b+c-d)/2$. L'image de $\theta_5$ 
agit transitivement sur le  type $(G_3\oplus H_3)$. On peut noter que 
le groupe  $\OO^t_F(\R)$ a 8 composantes connexes (voir  
proposition~\ref{prop:dec1},c)) et que  $\theta_5$ s'étend en un 
morphisme surjectif de noyau $\{\pm 1\}$ de
$\SSS^1\times \SL_2(\R)$ sur  $\OO^t_F(\R)_0$, morphisme
correspondant au revêtement double $\SSS^1\times \SU(1,1) 
\to \U(1,1)$. Le type est paramétré par 
$B_z=\theta_5(M_z)\cdot (\La_3\oplus \La_3)$ 
(voir eqs.~\ref{equa:L3} et \ref{equa:Mz}) avec $z\in\mfH_1$.
Posons $\Gamma=\theta_5^{-1}(\GL_6(\Z)\cap \OO^t_F(\R))$, sous-groupe
discret de $\SL_2(\R)$. L'expression de $\theta_5$ montre que 
$$
{\small
\begin{pmatrix}1 & 4 \\ 0 & 1 \end{pmatrix} \in \Ga, \esp
\begin{pmatrix} 2 & 8 k+3 \\ 1 & 4 k+2 \end{pmatrix}\in \Ga \esp
\mathrm{et}\esp
\begin{pmatrix} 0 & -1 \\ 1 & 4 k \end{pmatrix}\in \Ga.
}
$$
Toute orbite de $\Ga$ contient un point
$z\in\mfH_1$ avec $\Im \ga(z)$ maximum (les orbites de $\Ga$ ne 
s'accumulent pas sur les pointes du groupe, voir la preuve de 
la proposition~\ref{prop:max-I2G3}), c'est-à-dire avec
$|cz+d|^2\geq 1$ pour tout $\ga=\smallmat{a}{b}{c}{d}\in \Ga$.
Grâce aux éléments ci-dessus, on peut prendre $z$ extérieur aux
disques de rayon 1 centrés aux points réels $2k$ ($k\in \Z$) avec
de plus $|\Re z| \leq 2$. Ce domaine peut être réduit par les actions
entières  de $U_6$ qui change le signe de $\Re z$ ($U_6\cdot
B_z=B_{\overline{z}})$ et de 
$$
{\small
S_0=\begin{pmatrix}
0 & 1 & 0 & -1 & 0 & 0 \\
0 & 1 & 1 & 0 & -1 & 0 \\
0 & -1 & 1 & 0 & 0 & -1 \\
1 & 0 & 1 & 0 & -1 & 0 \\
0 & 0 & 1 & 0 & -1 & -1 \\
0 & -1 & 0 & 0 & 1 & -1
\end{pmatrix}.
}
$$
\begin{figure}[t]
\labellist
\small\hair 2pt
\pinlabel $0$  at 97 -10
\pinlabel $1$  at 192 -10
\pinlabel $2$  at 288 -10
\pinlabel $z_0$  at 194 128
\pinlabel $z_1$  at 178 178
\pinlabel $\mcF$  at 250 175
\endlabellist
\centering
\includegraphics[scale=0.45]{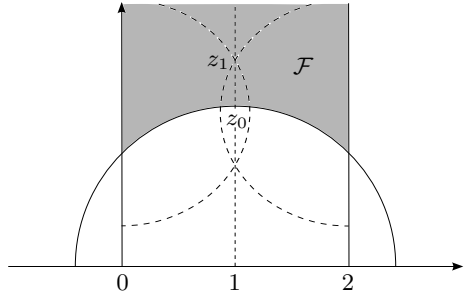}
\vspace{8pt}
\caption{Type géométrique $(G_3 \oplus H_3)$}
\label{figu:G3H3}
\end{figure}
L'action de $S_0$ sur $(G_3\oplus H_3)$ coïncide avec celle 
de $\theta_5(M_0)$ (non entière) pour 
$M_0 = 2^{-1/2}\smallmat{-1}{3}{-1}{1}\in \SL_2(\R)$, élément 
elliptique d'ordre 2 fixant $z_0=1+i\sqrt{2}$. Finalement,
on peut supposer que $z\in \mcF=\{z\in\mfH_1;
0\leq \Re z\leq 2, |z-1|^2\geq 2\}$ (figure~\ref{figu:G3H3}). 
La majoration de $\mu$ est obtenue grâce aux fonctions suivantes :
$$
\left\{
\begin{array}{l}
l_1=B_z[\eps_1]=\frac{1}{\sqrt{2}}(1+\frac{|z-i|^2}{2\Im z})
+\frac{1}{2}, \\
l_2=B_z[\eps_2-\eps_4]=
\frac{1}{\sqrt{2}}(1+\frac{|z-(2+i)|^2}{2\Im z}) +\frac{1}{2},\\
l_3=B_z[\eps_2-\eps_5] = \frac{2\sqrt{2}}{\Im z},
\end{array}
\right.
$$
où $(\eps_i)_{1\leq i \leq 6}$ désigne la base naturelle de $\Z^6$.
Soit $z_1=1+i(\sqrt{13}-1)/\sqrt{2}$. Le minimum
de la forme $B_{z_1}$ vaut $l_j(z_1)=(1+\sqrt{13})/3$ ($j=1,2,3$).
De plus, les fonctions $l_1$ et $l_2$ s'expriment au moyen des distances
respectives $d_1$ et $d_2$ aux points $i$ et $(2+i)$ par 
$l_j=(\sqrt{2}d_j +1)/2$  ($j=1,2$) (on rappelle que la distance 
$d_{z,w}$ entre  $z,w\in \mfH_1$ est donnée par 
$2 \Im z \Im w(\cosh d_{z,w} -1) =|z-w|^2$).  Par suite on a  
$\min(l_1,l_2,l_3)\leq (1+\sqrt{13})/3$ sur  le domaine $\mcF$ 
avec égalité au point $z_1$  uniquement (voir figure~\ref{figu:G3H3}).

\smallskip
\begin{rema}
Comme $(F_6)$, le type   $(G_3\oplus H_3)$ est inclus dans le type 
symplectique $\mfS_3$ (voir rem.~\ref{rema:F6}).  La forme $B_{z_0}$ 
 correspond à  la jacobienne de  la quartique de Fermat.
\end{rema}

\begin{prop}\label{prop:max-I3G3}
Le maximum de $\mu$ sur le type $(-I_3 \oplus G_3)$ 
vaut $7/5$, réalisé  de façon unique à équivalence près 
dans $(-I_3 \oplus G_3)$.
\end{prop}

\pr Nous procédons suivant la méthode de 
\cite[\S\S 1.3-1.4]{Bavard2007a}. Soit $PC_f^-$ le projectifié du 
cône négatif (dans $\R^4$) de la forme  $f=x^2+y^2+z^2-t^2$.
Le type $V_F$ ($F=-I_3 \oplus G_3$) est l'image d'un plongement
de  $PC_f^-$ dans $P_6$, équivariant pour un certain morphisme 
$\theta_6: \OO_f(\R) \to \GL_6^\pm(\R)$ (voir ci-dessous). 
 En utilisant  les  réflexions du groupe $\Pi= \theta_6^{-1}(\GL_n(\Z))$,
on se ramène à un polyèdre hyperbolique, ici une pyramide $\mcP$.
La détermination de $\mu_F$ résulte alors
 d'un simple découpage de $\mcP$ (ici en deux sous-pyramides).
\par
\begin{figure}[t]
\labellist
\small\hair 2pt
\pinlabel $p_0$  at 110 55
\pinlabel $p_1$  at 165 125
\pinlabel $s_0$  at 10 210
\pinlabel $s_1$  at 60 80
\pinlabel $s_2$  at 117 132
\pinlabel $s_3$  at 215 132
\pinlabel $s_4$  at 180 32
\pinlabel $e_5$  at 255 280
\endlabellist
\centering
\includegraphics[scale=0.40]{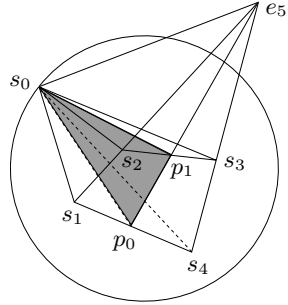}
\caption{Type géométrique $(-I_3\oplus G_3)$}
\label{figu:I3mG3}
\end{figure}

 
Notons  $[x,y,z,t]$ les coordonnées homogènes de $P\R^4$.
La forme polaire de~$f$ sera simplement notée 
$u\cdot v$ ($u,v\in \R^4$),  Pour
$p=[x,y,z,t]\in PC_f^-$, posons $\psi(p)=(x/t,y/t,z/t)'$, puis
$$
{\small
\ph(p)=Q_0\cdot \left(
\ph_{3,1}\circ \psi(p) \oplus I_2 \right)
~~~
\mathrm{o\grave{u}}
~~~
Q_0=\begin{pmatrix}
1 & 0 & 0 & 0 & 0 & 0 \\
0 & 0 & 1 & 0 & 0 & 0 \\
0 & 1 & 0 & 0 & 0 & 0 \\
0 & 0 & 0 & 2^{-1/2} & -2^{-3/4} & 2^{-3/4} \\
0 & 0 & 0 & 0 & 0 & 2^{1/4} \\
0 & 0 & 0 & 0 & 2^{1/4} & 0
\end{pmatrix},
}
$$
le plongement $\ph_{p,q}$ étant défini en \eqref{equa:Vpq}.
On obtient un paramétrage de $V_F$,  équivariant pour 
$\theta_6:\OO_f(\R) \to \GL_6^\pm(\R)$ défini par 
$\theta_6(P)=Q_0(P\oplus I_2)Q_0^{-1}$. À tout vecteur positif
$e\in\R^4$ est associée la réflexion 
$R_e=I_4-2 e I_{3,1} e'/e\cdot e$. Pour les  vecteurs suivants 
$$\begin{array}{lll}
e_1=(0,0,1,0)', & e_3=(-1,1,0,0)', & e_5= (1,0,0,1/\sqrt{2})',\\
e_2=(0,-1,1,0)', & e_4=(1,1,1,\sqrt{2})', &\\

\end{array}
$$
on constate que les réflexions $R_{e_i}$ ($i=1,\ldots 5$) appartiennent 
à $\Pi= \theta_6^{-1}(\GL_n(\Z))$. Plus précisément, les 
$\theta_6(R_{e_i})$ ($i=1,\ldots 5$) %
 sont donnés respectivement par
$I_{1,1}\oplus I_4$, $I_1\oplus U_2\oplus I_3$,
$$
{\small
\begin{pmatrix}0 & 0 &1\\0 & 1 &0\\1 & 0 & 0\end{pmatrix}
\oplus I_3,
~~~
\begin{pmatrix}
-1 & -2 & -2 & 4 & -2 & 2 \\
-2 & -1 & -2 & 4 & -2 & 2 \\
-2 & -2 & -1 & 4 & -2 & 2 \\
-2 & -2 & -2 & 5 & -2 & 2 \\
0 & 0 & 0 & 0 & 1 & 0 \\
0 & 0 & 0 & 0 & 0 & 1
\end{pmatrix}
~~~\mathrm{et}~~~
\begin{pmatrix}
-3 & 0 & 0 & 4 & -2 & 2 \\
0 & 1 & 0 & 0 & 0 & 0 \\
0 & 0 & 1 & 0 & 0 & 0 \\
-2 & 0 & 0 & 3 & -1 & 1 \\
0 & 0 & 0 & 0 & 1 & 0 \\
0 & 0 & 0 & 0 & 0 & 1
\end{pmatrix}.
}
$$
L'action du sous-groupe $\Pi'\subset \Pi$ engendré par les $R_{e_i}$
admet comme domaine fondamental le polyèdre $\mcP$ défini
par les équations $p\cdot e_i\leq 0$ ($i=1,\ldots,5$), c'est-à-dire
$0\leq z\leq y\leq x$, $x+y+z\leq \sqrt{2}t$ et $x\leq t/\sqrt{2}$.
Combinatoirement~$\mcP$ est une pyramide (cône sur un produit
de simplexes) dont le sommet est le point à l'infini
$s_0=[1,1,0,\sqrt{2}]$ ; la base de $\mcP$ est engendrée par les
points $s_1=[0,0,0,1]$, $s_2=[1,0,0,\sqrt{2}]$, 
$s_3=[1,1/2,1/2,\sqrt{2}]$ et $s_4=[1,1,1,3/\sqrt{2}]$ 
(voir figure~\ref{figu:I3mG3}). Pour $p_0=[1,1,1,4\sqrt{2}]$,
on trouve une forme de minimum $7/5$ avec~4 longueurs minimales 
$l_i(p)=1-2(p\cdot e_i)^2/p\cdot p$ ($i=1,\ldots 4$), qui varient 
comme les distances aux faces correspondantes. 

Découpons  $\mcP$ 
par $l_4=l_1$ (plan médiateur des faces $x+y+z -\sqrt{2}t=0$ et
$z=0$). Pour majorer les fonctions convexes
$l_i$ sur un polyèdre convexe, il suffit de les évaluer aux 
points extrémaux,  \cadab aux sommets (voir \cite{Bavard2007a}).
Le découpage scinde  $\mcP$ en deux pyramides de sommets respectifs
$(s_0,s_1,s_2,p_1,p_0)$ et $(s_0,p_0,p_1,s_3,s_4)$, avec
$p_1=[1,1/3,1/3,\sqrt{2}]$ (voir fig~\ref{figu:I3mG3}).
La forme correspondante à $p_1$ est de minimum $8/7 < 7/5$. Par suite,
pour $p\in \mcP\setminus\{p_0\}$, on a $l_1(p) < 7/5$ d'un côté
de $l_1=l_4$ et $l_4(p) < 7/5$ de l'autre côté, d'où le 
résultat.

\subsection{Types géométriques principaux de rang $7$}
La classification des types algébriques principaux de rang $7$ 
(tables~\ref{tabl:types_7_a} et \ref{tabl:types_7_b}) fait apparaître 
83 paires $\{[F],[-F]\}$. Parmi les types géométriques correspondant 
 nous trouvons~39 relations :  
\renewcommand{\arraystretch}{1.1}
$$
\begin{array}{ll}
(F \oplus -L_4)=(F \oplus L_4) & F=I_3,I_{2,1},I_1\oplus J_2,\\
(-I_4 \oplus G_3) = (I_{3,1} \oplus G_3),
 & (I_{1,3} \oplus G_3) = (I_{2,2} \oplus G_3),\\
(F \oplus H_3)= (F  \oplus G_3) &  F=I_4,U_4,J_4,F_4,K_4,\\
(J_2\oplus F \oplus H_3)=(J_2\oplus F \oplus G_3) &
F= \pm I_2, U_2,\\
(I_1\oplus H_3\oplus H_3)=(I_1\oplus G_3\oplus G_3) &
(I_1\oplus -G_3\oplus H_3)=(I_1\oplus -G_3\oplus -G_3) \\
(F_7)=(G_7), & (L_4\oplus H_3)=(-L_4\oplus G_3),\\
(F\oplus G)= (F\oplus F_5) & F=I_2,U_2,~G=H_5,G_5,K_5,\\
(F)=(I_{1,1}\oplus F_5) & F=-I_2\oplus F_5, I_{1,1}\oplus H_5, 
-I_2\oplus H_5, \\
(J_2\oplus F)=(J_2\oplus F_5) & F= H_5,G_5,K_5,\\
(I_1\oplus F)=(I_1\oplus K_6) & F=-K_6,\pm L_6,\pm M_6,\pm N_6, \\
(F)\simeq(H_7) & F= K_7,L_7,M_7, 
\end{array}
$$  %
dont beaucoup proviennent  de types scindés 
(voir proposition~\ref{prop:dec1}-b)  ou des dimensions inférieures. 
Finalement, cela conduit à  44 types géométriques principaux en 
dimension 7 (table~\ref{tabl:types_geom_7}), dont 12 sont 
maximaux (voir également la table~\ref{tabl:incl_tg7} des inclusions) : 
les 3 types symétriques
indéfinis, $F\oplus J_2$ ($F=I_{4,1},I_{3,2}$), $I_{2,1}\oplus J_4$,
$I_1\oplus J_6$, $F\oplus F_4$ ($F=I_3,I_{2,1},I_1\oplus J_2$) et 
$F\oplus G_3$ ($F=I_1\oplus G_3^-, K_4$).
\par 
L'analyse de la table~\ref{tabl:types_geom_7} montre que les types
sont généralement différenciés par leur  structure métrique, le cardinal
du groupe $\Ga_F$ et les estimations de la constante $\mu_F$. Il y 
a~4~exceptions : $(I_2\oplus F_5)$ et $(I_1\oplus K_6)$ (types 39 et 43)
 de dimension $0$ qui n'ont pas le même nombre de vecteurs 
minimaux, et 3 paires que l'on peut distinguer par le critère 
d'inclusion (proposition~\ref{prop:inclusion}) :
$(F\oplus G_3)$ avec $F\in \{I_{2,2},U_4\}$ (types 23,24), 
$(J_2 \oplus F \oplus G_3)$ avec $F\in \{I_2^- , U_2\}$ (types 26,27)  et
$(F\oplus F_5)$ avec $F\in \{I_{1,1},U_2\}$ (types 40,41).
\par 
Le calcul de $\mu_F$ est immédiat pour les types de dimension
$0$ et pour les types scindés.  Ensuite, certains types, non forcément
scindés, possèdent des termes diagonaux constants.
Ainsi $a_{7,7}=\sqrt{2}$ pour $A=(a_{ij})\in V_F$, 
$F=G\oplus G_3$ ($G=I_{3,1}, I_{2,2}, U_4$),
$I_1\oplus G_3^-\oplus G_3^-$,  $F_4\oplus G_3$ et $F_7$ 
(types 22,23,24,31,35 et 36) ; on trouve alors assez facilement un
point (généralement non isolé) de minimum $\sqrt{2}$ dans le type, 
d'où $\mu_F= \sqrt{2}$. Nous n'incluons pas les détails car ces
types ne sont pas maximaux. Un phénomène analogue a lieu pour 
$F=I_1\oplus G_3\oplus H_3$ (type 32, en fait scindé) avec $\mu_F=1$ et 
pour  $F=G\oplus F_5$ ($G=I_{1,1},U_2$, types 40,41) avec $\mu_F= \nu$ où
\begin{equation}\label{equa:cte_nu}
\nu = \frac{1}{2} \sqrt{4+2\sqrt{2}} \simeq 1,3065.
\end{equation}

\begin{table}[!tbp]
\centering
\renewcommand\arraystretch{1.2}
\renewcommand\tabcolsep{2pt}
{\small \begin{tabular}{llrccccccl}
\hline
$d$  & $F F^\vee$& $n^o$ & $F$ & $V_F=(F)$& dim. & type sym.
& $\Ga_F$ & $\mu_F$ &$\subset_\simeq$\\
\hline
1&$I_7$ &1&$I_7$&$\{I_7\}$ &  0 & $\R^0$ &$\G(7)$ &1 & 2-4,6,7,9-12,14\\
&&2&$I_{6,1}$ & $V_{6,1}$ & 6 & $\HH^6_8$ & $\C_2$& 5/3& max.\\
&&3&$I_{5,2}$ & $V_{5,2}$ & 10 & $V_{5,2}$  & $\C_2$&$\geq 5/3$& max.\\
&&4&$I_{4,3}$ & $V_{4,3}$ & 12 & $V_{4,3}$  & $\C_2$&$\geq 5/3$& max.\\
\hline
2 & $I_{5,2}$&5&$I_5J_2$ & $\{I_5\} P_2$ & 2  & $\HH^2_2$ 
& $\G(5) \C_2$& 1 & 6,7,9,10,14\\
&&6&$I_{4,1}J_2$ & $V_{4,1} P_2$ & 6  & $\HH^4_8 \HH^2_2$ 
& $\C_2^2$& $2/\sqrt{3}$& max.\\
&&7&$I_{3,2}J_2$ & $V_{3,2}  P_2$ & 8  & 
$V_{3,2}  \HH^2_2$& $\C_2^2$& $2/\sqrt{3}$& max.\\
& $I_{3,4}$&8&$I_3J_4$ & $\{I_3\} \mfS_2$ & 6  & $\mfH_2$ 
& $\G(3) \C_2$& 1 & 9,10\\
&&9&$I_{2,1}J_4$ & $V_{2,1} \mfS_2$ & 8 & 
$\HH^2_8 \mfH_2$& $\C_2^2$& $\al$& max.\\
&$I_{1,6}$&10&$I_1J_6$ & $\{I_1\} \mfS_3$ & 12  & $\mfH_3$
& $\C_2^2$& $1$& max.\\
&$I_3R_2R_2$&11&$I_3  F_4$ & $(I_3  F_4)$& 6 & 
$\HH^4_8 \HH^2_2$ &$\C_2^2$& $\geq 5/3$& max.\\
& &12&$I_{2,1}  F_4$ & $(I_{2,1}  F_4)$& 8 & 
$V_{3,2} \HH^2_2$ &$\C_2^2$& $\geq 5/3$ & max.\\
& &13&$I_1F_6$ & $\{I_1\}(F_6)$& 2 & 
$\HH^2_2$ &768& 1 & 10-12,14\\
&$I_{1,2}R_2R_2$&14&$I_1 J_2  F_4$ & $(I_1 J_2  F_4)$& 8 & 
$\HH^2_8  \mfH_2$ & $\C_2^2$& $\geq 5/3$ & max. \\
\hline
4& $I_3J_2J_2$ &15& $I_3 K_4$ & $\{I_3\} (K4)$ & 2 &
$\HH^2_4$ & $\G(3) \D_{8}$& 1 & 9,10\\
 & &16& $I_{2,1} K_4$ & $ V_{2,1} (K4)$ & 4 &
$\HH^2_8 \HH^2_4$ & $\C_2 \D_{8}$& $\al$ & 9\\
& &17& $I_3 L_4$ & $\{I_3 D_4\}$ & 0& 
$\R^0$ & $\G(3)  \Ga_{L_4}$& 1 & 3,4,7,9-12,14\\
 & &18& $I_{2,1} L_4$ & $ V_{2,1}\{D_4\}$ & 2 &
$\HH^2_8$ & $\C_2 \Ga_{L_4}$& $\al$ & 4,9,12\\
& $I_{1,2}J_2J_2$ &19& $I_1J_2 K_4$ & $\{I_1\}P_2(K4)$ & 4 &
$\HH^2_2 \HH^2_4$ & $\C_2^2 \D_{8}$& $1$ & 10\\
& &20& $I_1J_2 L_4$ & $\{I_1\} P_2 \{D_4\}$ & 2 & 
$\HH^2_2$ & $\C_2^2 \Ga_{L_4}$& $1$ & 7,10,14\\
&$I_4X_3$ &21& $I_4G_3$ & $\{I_4 \La_3\}$ &0 & $\R^0$ &
$\G(4)\C_2\D_8$ & 1 & 2-4,6,7,9,12,37\\
&&22& $I_{3,1}G_3$ & $(I_{3,1}G_3)$ &4 & $\HH^4_8$ &
$\C_2\D_8$ & $\sqrt{2}$ & 3\\
&&23& $I_{2,2}G_3$ & $(I_{2,2}G_3)$ &6 & $V_{3,2}$ &
$\C_2\D_8$& $\sqrt{2}$ & 4\\
&&24& $U_4G_3$ & $(U_4G_3)$ &6 & $V_{3,2}$ &
$\C_2\D_8$ & $\sqrt{2}$ & 4\\
&$I_{2,2}^-X_3$ &25& $J_2I_2G_3$ & $P_2\{I_2 \La_3\}$ &2 & $\HH^2_2$ &
$\C_2^2\D_8^2$ & 1 & 6,7,9\\
& &26& $J_2I_2^-G_3$ & $P_2(I_2^-G_3)$ &4 & $\HH^2_2 \HH^2_8$ &
$\C_2^2 \D_8$ & $2/\sqrt{3}$ & 7\\
& &27& $J_2U_2G_3$ & $P_2(U_2G_3)$ &4 & $\HH^2_2 \HH^2_8$ &
$\C_2^2 \D_8$ &$2/\sqrt{3}$  &7\\
&$I_4^-X_3$ &28& $J_4G_3$ & $\mfS_2 \{\La_3\}$ &6 &
$\mfH_2 $ & $\C_2^2\D_8$ & $\al$ & 9\\
&$I_1X_3X_3$ &29& $I_1 G_3 G_3$ & $\{I_1 \La_3^2\} $ &0 & $\R^0$ & 1024
& $1$ & 3,4,10,12,14,30\\
&&30& $I_1 G_3 G_3^-$ & $(I_1 G_3 G_3^-) $ &4 & $\HH^2_8 \HH^2_4$ & 
$\C_2 \C_4$ & $\geq 5/3$ & max.\\
&&31& $I_1 G_3^- G_3^-$ & $(I_1 G_3^- G_3^-) $ &2 &$\HH^2_8 $ & 128
& $\sqrt{2}$ & 4,30\\
&&32& $I_1 G_3 H_3$ & $(I_1 G_3 H_3) $ &2 &$\HH^2_4 $ & 64
& $1$ & 10,14,30\\
&&33& $I_1 G_6$ & $(I_1 G_6) $ &4 &$\HH^2_8 \HH^2_4 $ & $\C_2\D_8$
& $\geq \psi$ & 14\\
&$I_1X_3^-X_3^-$ &34& $I_1H_6$ & $\{I_1\}(H_6)$ &4 & $\HH^2_2 \HH^2_4$ & 
$\C_2^2\D_8$ & 1 & 10\\
&$R_2R_2X_3$ &35& $F_4G_3$ & $(F_4G_3)$ &4 & $\HH^2_2 \HH^2_8$ & 
$\C_2^2\D_8$ & $\sqrt{2}$ & 12\\
 &  & 36 & $F_7$ & $(F_7)$ & 2  & $\HH^2_2$  & 128 & $\sqrt{2}$  &  11,12 \\
&$J_2J_2X_3$ &37& $K_4G_3$ & $(K_4G_3)$ &4 & $W_{2,1}$ & 
$\C_2\C_4$  & $\geq \ph$ & max.\\
 &  & 38 & $L_4^-G_3$ & $\{D_4\La_3\}$ & 0  & $\R^0$  &
  18432 & $\al$  &  4,9,12,30,37\\
\hline
8 & $I_2X_5$ & 39& $I_2F_5$ & $\{I_2\La_5\}$ &0&$\R^0$ & $\D_8\C_2\D_{16}$
& 1 & 3,4,7\\
 &  & 40& $I_{1,1}F_5$& $(I_{1,1}F_5)$ &2&$\HH^2_8$ & $\C_2\D_{16}$
& $\nu$ & 4\\
 &  & 41& $U_2F_5$& $(U_2F_5)$ &2&$\HH^2_8$ & $\C_2\D_{16}$
& $\nu$  & 4\\
 &$I_2^-X_5$  & 42& $J_2F_5$& $\mfS_2 \{\La_5\}$ &2&$\HH^2_2$ &
$\C_2^2\D_{16}$ & $2/\sqrt{3}$ & 7\\
 & $I_1Y_6$ &43&$I_1K_6$&$\{I_1\La_6\}$&0&$\R^0$&256&$1$& 4\\
 &$T_7$&44&$H_7$&$\{\La'_7\}$&0&$\R^0$&$\C_2\D_{16}$&$\om$& 4\\
\hline
\end{tabular}}
\caption{Types géométriques principaux en dimension $7$ (44 types, 12 max.)}
\label{tabl:types_geom_7}
\end{table}

Il reste à déterminer $\mu_F$ pour  7  (sur 12) des types maximaux 
et d'un type non maximal $(I_1\oplus G_6)$ inclus dans le type maximal
$(I_1\oplus J_2\oplus F_4)$. Pour chacun de ces types, nous indiquons
dans la table~\ref{tabl:types_geom_7} un maximum local de densité. 
Voici pour terminer quelques précisions sur ces cas.
\par
Rappelons d'abord que la densité  maximale du type symétrique 
\og lorentzien\fg\  est déterminée dans \cite{Bavard2007a} 
(théorème 1 p. 44) : elle correspond à $\mu_{6,1}=5/3$, réalisée
de façon unique à équivalence près par la  forme 
$$
{\small
\La_7= \frac{1}{3}
\begin{pmatrix}
5 & 2 & 2 & 2 & 2 & 2 & 6\\
2 & 5 & 2 & 2 & 2 & 2 & 6\\
2 & 2 & 5 & 2 & 2 & 2 & 6\\
2 & 2 & 2 & 5 & 2 & 2 & 6\\
2 & 2 & 2 & 2 & 5 & 2 & 6\\
2 & 2 & 2 & 2 & 2 & 5 & 6\\
6 & 6 & 6 & 6 & 6 & 6 & 15
\end{pmatrix}.
}
$$
Cette forme admet des points équivalents dans les types maximaux $V_{5,2}$,
$V_{4,3}$, $F\oplus F_4$ ($F=I_3,I_{2,1},I_1\oplus J_2$) et 
$I_1\oplus G_3^- \oplus G_3$ ; dans chaque cas on vérifie, 
grâce notamment à  la proposition~\ref{prop:esp_tangent},
qu'il s'agit d'un maximum local. 
\par 
Le dernier type maximal non évoqué ci-dessus, $(K_4 \oplus G_3)$, est 
métriquement intéressant car il s'agit d'un plan hyperbolique  complexe. 
En réduisant son type algébrique sur les réels, on trouve que 
$$
(K_4 \oplus G_3)= P_4 \cdot (I_1\oplus W_{2,1}) ,
$$
où $P_4$ est donné par \eqref{equa:ann_P4} et $W_{2,1}$ est paramétré 
explicitement par un couple $(z,w)$
de   nombres complexes tels que $|z|^2+ |w|^2 < 1$ {\it via} les 
équations~\eqref{equa:kappa}, \eqref{equa:mcWpq} et \eqref{equa:Wpq},
à savoir $W_{2,1} =\{\kappa(H_{z,w}) ; |z|^2+ |w|^2 < 1\}$ avec
$$
{\small
H_{z,w}  = \frac{1}{ 1-|z|^2-|w|^2}
\left( 
\begin{array}{ccc}
1+|z|^2-|w|^2 & 2 z \overline{w} & 2 z \\
2\overline{z}w & 1 - |z|^2+|w|^2 & 2w \\
2\overline{z} & 2 \overline{w} & 1+|z|^2+|w|^2
\end{array}
\right).
}
$$
Posons 
\begin{equation}\label{equa:cte_ph}
\left\{
\begin{array}{l}
\ph   =  1+\frac{1}{2}(\sqrt{6}-\sqrt{2}) \simeq 1,5176,\\
97 z_0  =  -36-113 i+ (10+26 i) \ph + (31+197 i )\frac{\ph^2}{2} + 
(1-75 i ) \frac{\ph^3}{2}\\
97 w_0  =  16 + 61 i + (-26 + 10 i) \ph + (47-130 i ) \ph^2 + (-11+49 i) \ph^3 .
\end{array}
\right.
\end{equation}
Alors la forme $\kappa(H_{z_0,w_0})\in (K_4 \oplus G_3)$ est de minimum $\ph$ et 
possède 10 couples de vecteurs  minimaux. On vérifie par la méthode 
habituelle (proposition~\ref{prop:esp_tangent}) qu'elle atteint 
un maximum local (isolé) de densité sur le type $(K_4 \oplus G_3)$.

\begin{table}[!b]
\centering
\renewcommand{\arraystretch}{1.2}
\renewcommand\tabcolsep{5pt}
{\scriptsize \begin{tabular}{rlrlrl}
\hline
1 & 0-12,14-16,19 &  %
16 & 9,16 &          %
31 & 4,23,30,31\\    %
2 & 2 &              %
17 & 0,3,4,7-12,14-20,22,23,26 &  %
32 & 10,14,30,32\\   %
3 & 3 &              %
18 & 0,4,9,12,16,18,23 &    %
33 & 14,33\\         %
4 & 4 &              %
19 & 10,19 &         %
34 & 10,34\\         %
5 & 0,5-10,14,19 &   %
20 & 0,7,10,14,19,20,26 &  %
35 &  12,35\\        %
6 & 6 &              %
21&0,2-4,6,7,9,12,16,21-28,35,37 & %
36 & 11,12,35,36\\   %
7 & 7 &              %
22 & 3,22 &          %
37 & 37\\            %
8 & 0,8-10&          %
23 & 4,23 &          %
38 & 0,4,9,12,16,18,\\ %
9 & 9&               %
24 & 4,24 &          %
 & 23,28,30,31,35,37,38 \\ %
10 & 10 &            %
25 & 6,7,9,25-28 &   %
39 & 3,4,7,39-42\\         %
11 & 11 &           %
26 & 7,26 &         %
40 & 4,40\\        %
12 & 12 &           %
27 & 7,27 &         %
41 & 4,41  \\        %
13 & 0,10-14&       %
28 &9,28 &          %
42 & 7,42\\  %
14 & 14 &           %
29 & 3,4,10,12,14,22,23,29-34 & %
43 & 4,23,40,43 \\        %
15 & 0,8-10,15,16,19 & %
30 & 30 &           %
44 & 4,44   \\
\hline
\end{tabular}}
\caption{Inclusions des types principaux en dimension $7$}
\label{tabl:incl_tg7}
\end{table}

Enfin, il reste à commenter le cas du type non maximal 
$(I_1\oplus G_6)$ (type 33,
inclus à équivalence près dans $(I_1\oplus J_2\oplus F_4)$),
que l'on peut paramétrer par 
$(I_1\oplus G_6)=\{A_{X,z};(X,z)\in \mcV_{2,1} \times\mfH_1\}$ avec 
$$
A_{X,z} = \frac{1}{2} 
\begin{pmatrix}
2 & 0 & 0 & 0 & 0 & 0 & 0\\ 
0 & 1 & 1 & 0 & 0 & 0 & -2^{3/4}\\ 
0 & 0 & 0 & 0 & 2^{3/4} & 0 & 
-2^{3/4}\\ 
0 & 0 & 0 & 0 & 2^{3/4} & 0 & 2^{3/4}\\ 
0 & 1 & -1 & -2^{3/4} & 0 & 0 & 0\\ 
0 & 0 & 0 & -2^{3/4} & 
0 & -2^{3/4} & 0\\ 
0 & 0 & 0 & 2^{3/4} & 0 & -2^{3/4} & 0
\end{pmatrix}
\cdot \left( \ph_{2,1}(X)\oplus \sig_1(z) \oplus \sig_1(z) \right)
$$
(voir les équations \eqref{equa:Sig_g}, \eqref{equa:Om_pq} et 
\eqref{equa:Vpq}). Soient $\psi, \psi'$ et $\psi''$  définis par 
\begin{equation}\label{equa:cte_psi}
\left\{
\begin{array}{l}
63 \psi^4-64 \psi^3-116 \psi^2+128 \psi-36 = 0\esp 
\mathrm{et}\esp \psi~\mathrm{r\acute{e}el} >0, \\
144 \psi' = -198+404 \psi+65 \psi^2-126 \psi^3,\\
\psi''= 2\psi - 2\psi' -1,
\end{array}
\right.
\end{equation}
et soit 
$$A_0 =
\begin{pmatrix}
2\psi -1 & \psi' & 0 & 0 & -\psi' & 0 & 0  \\
\psi' & \psi & \psi/2 & - \psi/2 &  (1 - \psi)/2  & \psi'' & \psi'' \\
0 & \psi/2 & \psi & 0 & -\psi'' & 0 & 2\psi''  \\
0 & - \psi/2 & 0 & \psi & -\psi'' & -2\psi'' & 0 \\
-\psi' & (1 - \psi)/2  & -\psi'' & -\psi'' & \psi & \psi/2 & -\psi/2 \\
0 & \psi'' & 0 & -2\psi'' &  \psi/2 & \psi & 0\\
0 & \psi'' & 2\psi'' & 0 & - \psi/2  & 0 &\psi
\end{pmatrix}.
$$
Alors $A_0\in (I_1\oplus G_6)$ est une forme de minimum $\psi\simeq
1,5101$ et  possédant 16 couples de vecteurs minimaux. On vérifie qu'il
s'agit d'un maximum local  de densité (isolé) sur le type 
$(I_1\oplus G_6)$.
\par 
Précisons également que la forme $\La'_7$ (type 44) est donnée par 
$$
\La'_7  = {\small
\frac{1}{4}
\begin{pmatrix}
1+\sqrt{2}+c & \sqrt{2}+a & -\sqrt{2}-a  & -b & -c & -c & -b \\
\sqrt{2}+a & 2 \sqrt{2}+2 a & b & -2 b &  -2 a & 0 & 2 a \\
-\sqrt{2}-a & b & 2 \sqrt{2}+2 a & -2 a & 0 & 2 a & 2 b\\
-b & -2 b & -2 a & 4 a & 2 b & 0 & -2 b \\
-c & -2 a & 0 & 2 b & 4 a & 2 b & 0\\
-c & 0 & 2 a & 0 & 2 b & 4 a & 2 b\\
-b & 2 a & 2 b & -2 b & 0 & 2 b & 4 a
\end{pmatrix} 
}
$$
avec $a=\sqrt{2+\sqrt{2}}$, $b=\sqrt{4-2 \sqrt{2}}$ et 
$c=\sqrt{20+2 \sqrt{2}}$. Elle possède 8 couples de vecteurs
minimaux et son minimum vaut
\begin{equation}\label{equa:cte_omega}
\omega = \frac{1}{4} \left(1+ \sqrt{2} + \sqrt{4+ 2 \sqrt{2}}\right)
\simeq 1,2568.
\end{equation}

Les constantes $\alpha,\psi,\ph,\nu$ et $\omega$ de la table
\ref{tabl:types_geom_7} sont définies aux équations
\eqref{equa:cte_alpha},\eqref{equa:cte_psi},\eqref{equa:cte_ph},\eqref{equa:cte_nu} 
et \eqref{equa:cte_omega} respectivement. 
Enfin, les inclusions entre les types géométriques principaux 
en dimension~7 sont données dans la table~\ref{tabl:incl_tg7}.

\FloatBarrier

%

%

\addtocounter{section}{1}
\addcontentsline{toc}{section}{Annexe : données numériques}
\section*{Annexe : données numériques}
\setcounter{subsection}{0}
\renewcommand{\thesubsection}{A\arabic{subsection}}
Signalons d'abord quelques conventions générales. Dans les
tables, la somme directe $A \oplus B = \smallmat{A}{0}{0}{B}$ 
sera notée multiplicativement
pour abréger : $AB$ pour $A \oplus B$, $A^2$ pour $A \oplus A$ et 
$A^-$ pour $-A$. La matrice identité d'ordre $n$ est
notée $I_n$. On pose $I_{p,q}=I_p \oplus (-I_q)$ et 
$J_{2n}=\smallmat{0}{I_n}{-I_n}{0}$. Enfin pour 
$p$ premier on définit
$$
{\small 
R_p=\left(
\begin{array}{cc}
0 &  1\\
I_{p-1} &  0\\
\end{array}
\right), ~~~
V_{p-1}=
\left(
\begin{array}{cc}
0 & -1\\
I_{p-2}& \vdots\\
 & -1
\end{array}
\right), ~~~
W_{p-1}=
\left(
\begin{array}{cc}
0 &  -1\\
&  1\\
 I_{p-2}  &\vdots\\
& -1\\
&  1
\end{array}
\right).
}
$$

\subsection{Éléments d'ordre fini de $\GL_n(\Z)$}
\label{subs:of_GLnZ_ann}

\begin{table}[h]
\centering
\renewcommand{\arraystretch}{1.2}
\renewcommand{\tabcolsep}{7pt}
{\small \begin{tabular}{ccllr}
\hline
$n$  & $d$ & $\det =1$ & $\det =-1$&$N_{n,d}$ \\
\hline
1 & $2$  & & $I_1^{-*}$& 1\\
\hline
2 & 2 & $I_2^-$ & $I_{1,1},~R_2^*$& 3\\
  & 3 & $V_2^*$ & & 1 \\
  & 4 & $J_2^*$ && 1\\
  & 6 & $W_2^*$ && 1\\
\hline
3 & 2  & $I_{1,2},~I_1^-R_2$ & $I_{2,1},~I_1R_2,~ I_3^-$&5\\
  & 3  & $I_1V_2,~R_3^*$&&2\\
  & 4  & $I_1J_2,~X_3^*$ & $I_1^-J_2,~X_3^{-*}$&4\\
  & 6  & $I_1W_2$ & $I_1^-W_2,~I_1^-V_2,~R_3^{-*}$&4\\
\hline
4 & 2  & $I_4^-,~I_{2,2},~I_{1,1}R_2,~R_2^2$& 
$I_{3,1},~I_{1,3},~I_2R_2,~I_2^-R_2$& 8\\
& 3  & $I_2V_2,~I_1R_3,~V_2^2$&& 3\\
& 4  & $I_2J_2,~I_2^- J_2,~J_4,~I_1X_3,~I^-_1X^-_3$&
$I_{1,1}J_2,~R_2J_2,~I_1X_3^-,$& \\
& & & $I_1^-X_3,~Z_4^*,~Z'_4\null^{*},~Z''_4\null^{*}$ & 12\\
& 5  & $V_4^*$& & 1\\
& 6  & $I_2 W_2,~I_2^- W_2,~I_2 V_2,~I_1^-R_3^-$,& 
$I_{1,1}V_2,~R_2V_2,~I_1^-R_3,~T'_4\null^{*},$&\\ 
&  & $V_2W_2,~W_2^2,~T_4^*$ &
 $I_{1,1}W_2,~R_2W_2,~I_1R_3^-,~T''_4\null^{*}$&  15\\
& 8  & $X_4^*$& &1\\
& 10  & $W_4^*$& & 1\\
& 12  & $V_2J_2,~W_2J_2,~Y_4^*$& & 3\\
\hline
\end{tabular}}
\caption{\'Eléments d'ordre fini $d \geq 2$
de $\GL_n(\Z)$ ($1 \leq n \leq 4$)}
\label{tabl:of1-4}
\end{table}

Dans les tables~\ref{tabl:of1-4} et \ref{tabl:of5-7}
l'entier $d$ désigne l'ordre des éléments
($d=1$ n'est pas rappelé) ;
les in\-dé\-com\-po\-sa\-bles sont signalés par un astérique (*).
Noter qu'en dimension $n$ impaire, les éléments de déterminant 1
correspondent bijectivement à ceux de déterminant -1.
On pose 
\renewcommand{\arraycolsep}{0.3 em}
\renewcommand{\arraystretch}{1}
{\small
$$X_3=\left(
\begin{array}{ccc}
1 & 0 & 1\\
0 & 0 & 1\\
0 & -1&0
\end{array}
\right),~
X_4=\left(
\begin{array}{cc}
0 & -1 \\
I_3 & 0
\end{array}
\right),~
Y_4=\left(
\begin{array}{cccc}
0 & 0 &0 &-1 \\
1 &  0&0 &0 \\
 0& 1 & 0 & 1\\
 0&  0& 1&0 
\end{array}
\right),
$$
$$~Z_4=\left(
\begin{array}{cc }
 I_{1,1}  & E'_{2,2}\\
0 & J_2
\end{array}
\right),~
Z'_4=\left(
\begin{array}{cc }
 R_2 & E'_{2,2}\\
0 & J_2
\end{array}
\right),~
Z''_4=\left(
\begin{array}{cc }
 R_2  & E_{2,2}\\
0 & J_2
\end{array}
\right),~
$$
$$
T_4=\left(
\begin{array}{cc }
V_2  & E_{2,2}\\
0 & W_2
\end{array}
\right),~
T'_4=\left(
\begin{array}{cc }
 R_2 & E_{2,2}\\
0 & V_2
\end{array}
\right),~
T''_4=\left(
\begin{array}{cc }
 R_2  & E_{2,2}\\
0 & W_2
\end{array}
\right),
$$}

\noindent
avec $E_{r,s} = e^{r,s}_{r,s}$ et $E'_{r,s} =
e^{r,s}_{r,s}+e^{r-1,s}_{r,s}$,
où $e^{i,j}_{r,s}=(\delta^i_k\delta^j_l)_{k,l} \in M_{r,s}(\Z)$
(matrice élémentaire).

Pour les dimensions supérieures, posons

\renewcommand{\arraycolsep}{0.25 em}
\renewcommand{\arraystretch}{1}
{\small
$$X_5=\left(
\begin{array}{cc}
1 & E_{1,4} \\
0 & X_4
\end{array}
\right),~
X_6=\left(
\begin{array}{cc}
J_2 & E_{2,4} \\
0 & X_4
\end{array}
\right),~
Y_6=\left(
\begin{array}{cc}
J_2 & E'_{2,4} \\
0 & X_4
\end{array}
\right),~
Z_6=\left(
\begin{array}{cc }
 I_{1,1}  & E'_{2,4}\\
0 & X_4
\end{array}
\right),~
$$
$$
Z'_6=\left(
\begin{array}{cc }
 R_2 & E_{2,4}\\
0 & X_4
\end{array}
\right),~
Z''_6=\left(
\begin{array}{cc }
 R_2  & E'_{2,4}\\
0 & X_4
\end{array}
\right),~
X_7=\left(
\begin{array}{cc }
J_2  & E_{2,5}\\
0 & X_5
\end{array}
\right),~
Y_7=\left(
\begin{array}{cc }
J_2  & E'_{2,5}\\
0 & X_5
\end{array}
\right),~
$$
$$
Z_7=\left(
\begin{array}{cc }
 X_3  & E_{3,4}\\
0 & X_4
\end{array}
\right),~
Z'_7=\left(
\begin{array}{cc }
 X_3  & e^{1,4}_{3,4}\\
0 & X_4
\end{array}
\right),~
Z''_7=\left(
\begin{array}{cc }
 X_3  & E'_{3,4}\\
0 & X_4
\end{array}
\right).
$$}

\begin{table}[h]
\centering
\renewcommand{\arraystretch}{1.2}
\renewcommand{\tabcolsep}{5pt}
{\small \begin{tabular}{cllll}
\hline
$d$  & $n$ & $\det =1$ & $\det =-1$&$N_{n,d}$ \\
\hline
2 & 5,~6,~7 & 
\multicolumn{2}{c}{$I_{p,q}R_2^r,~~ q+r >0,~~ 
\det=(-1)^{q+r}$} & 11,~15,~19\\
\hline
4 & 5,~6,~7 & 
\multicolumn{2}{c}{$I_{p,q}R_2^rJ_2^s X_3^t X_3^{-u}
Z_4^v Z'_4\null^{w} Z''_4\null^{x}$}  &
24,~48,~84\\
&& \multicolumn{2}{c}{$s+t+u+v+w+x >0, ~~
\det = (-1)^{q+r+u+v+w+x}$}&\\
\hline
8 & 5 & $I_1X_4,~X_5$&$I_1^-X_4^-,~X_5^-$ &4 \\
  & 6 & $I_2X_4,~I_2^-X_4,~I_1X_5,~I_1^-X_5^-,~J_2X_4,$& 
$I_{1,1}X_4,~I_1X_5^-,~I_1^-X_5,$& \\
  &  & $X_6^*,~Y_6^*$& 
$R_2X_4,~Z_6^*,~Z'_6\null^*,~Z''_6\null^*$& 14\\
  & 7 & $I_{p,q}X_4,~I_{p,q}X_5$ ($q$ pair), $I_1^-R_2X_4$,& 
$ -X$& 40\\
  &  &$I_1J_2X_4,~I_{1,1}X_5^-,~R_2X_5^-,~J_2X_5,$ & 
$(X \in \SL_7(\Z), ~{\rm ordre}~ 8$) & \\
  &  &$I_1X_6,~I_1Y_6,~I_1^-Z_6,~I_1^-Z'_6,~I_1^-Z''_6,$  & & \\
  &  &$X_3X_4,~X_7^*,~Y_7^*,~ Z_7^*,~Z'_7\null^*,~Z''_7\null^*$ & & \\
\hline
\end{tabular}}
\caption{\'Eléments d'ordre $2^k\geq 2$
de $\GL_n(\Z)$ ($5 \leq n \leq 7$)}
\label{tabl:of5-7}
\end{table}

\subsection{Types algébriques indécomposables}
\label{subs:alg_indec}

\noindent{\boldmath \bf Rang $\leq 2$ : }
$\pm I_1$, $U_2$, $J_2$, $\pm F_2$, avec \begin{equation}
\label{equa:indec_rg2}
\renewcommand{\arraycolsep}{0.2 em}
\renewcommand{\arraystretch}{1}
{\small
I_1=\left( \begin{array}{r} 1 \end{array} \right),~
U_2= \left( \begin{array}{rr} 0 & 1 \\ 1 & 0 \end{array} \right),~
J_2= \left( \begin{array}{rr} 0 & 1 \\ -1 & 0 \end{array} \right),~
F_2= \left( \begin{array}{rr} 1 & -1 \\ 0 & 1 \end{array} \right).}
\end{equation}

\medskip
\noindent{\bf Rang 3 : }
$\pm F_3$, $\pm G_3$ et $\pm H_3$, avec
\begin{equation}
\label{equa:indec_rg3}
\renewcommand{\arraycolsep}{0.2 em}
\renewcommand{\arraystretch}{1}
{\small
F_3 =\left( \begin{array}{rrr}
0 & 1  & 0 \\ 0 & 0 & 1 \\ 1 & 0 & 0 \end{array} \right),~
G_3 = \left( \begin{array}{rrr}
0 & -1  & ~0\\ 
0 & -1 & -1 \\ 
1 & ~1 & -1 \end{array} \right),~
H_3 = \left( \begin{array}{rrr}
1 & 1  & ~0\\ 0 & 1 & 1 \\ -1 & -1 & 1 \end{array} \right).}
\end{equation}

\medskip
\noindent{\bf Rang 4 : }
$F_4,~\pm G_4,~H_4,~K_4,~\pm L_4,~\pm M_4,~\pm N_4,~\pm O_4$, avec
\renewcommand{\arraystretch}{1}
\begin{equation}
\label{equa:indec_rg4}
{\small
\renewcommand{\arraycolsep}{0.2 em}
\begin{array}{lll}
F_4 = \left( \begin{array}{rrrr}
0 & 0 & 1 & 0 \\
0 & 0 & 0 & 1 \\
0 & 1 & 0 & 0 \\
1 & 0 & 0 & 0 \end{array}\right),&
G_4=  \left( \begin{array}{rrrr}
2 & 1 & 1 & 1 \\
1 & 1 & 0 & 1 \\
1 & 1 & 1 & 0 \\
1 & 0 & 1 & 1 
\end{array} \right),&
H_4=\left( \begin{array}{rrrr}
0 & ~0 & ~0 & -1 \\
0 & 0 & 1 & 0 \\
1 & 0 & 0 & 0 \\
1 & 1 & 0 & 0 
\end{array} \right),\\ \\
K_4=\left( \begin{array}{rrrr}
0 & 0 & 1 & 0 \\
0 & 0 & 0 & 1 \\
0 & 1 & 0 & 0 \\
-1 & 0 & 0 & 0 
\end{array} \right),&
L_4= \left( \begin{array}{rrrr}
1 & ~1 & 0 & ~1 \\
-1 & 1 & -1 & 0 \\
1 & 0 & 1 & 1 \\
0 & 1 & -1 & 1 
\end{array} \right),&
M_4= \left( \begin{array}{rrrr}
1 & -1 & 1 & -1 \\
0 & 1 & -1 & 1 \\
0 & 0 & 1 & -1 \\
0 & 0 & 0 & 1 
\end{array} \right),\\ \\
N_4= \left( \begin{array}{rrrr}
1 & ~0 & ~0 & -1 \\
0 & 0 & 1 & 0 \\
1 & 0 & 0 & 0 \\
-1 & 1 & 0 & 1 
\end{array} \right),&
O_4=  \left( \begin{array}{rrrr}
0 & 1 & -1 & 0 \\
0 & 0 & 1 & -1 \\
1 & -1 & 0 & 1 \\
-1 & 1 & 0 & 0 
\end{array} \right).&
\end{array}
}
\end{equation}

\medskip
\noindent{\bf Rang 5 (principaux) : }
$\pm F_5,~\pm G_5,~\pm H_5,~\pm K_5$, avec 
\begin{equation}
\label{equa:indec_rg5}
{\small
\renewcommand{\arraycolsep}{0.2 em}
\begin{array}{ll}
F_5 = \left( \begin{array}{rrrrr}
0& 0& 1& 0& 0 \\
0& 0& -1& 1& 0\\ 
0& 0& 0& -1& 1\\
1& -1& 0& 0& -1\\
0& 1& -1& 0& 0
\end{array} \right),&
G_5 = \left( \begin{array}{rrrrr}
1& 0& -1& 0& 0 \\
0& 0& 1& -1& 0\\ 
0& 0& 0& 1& -1\\
-1& 1& 0& 0& 1\\
0& -1& 1& 0& 0
\end{array} \right),\\ \\
H_5 = \left( \begin{array}{rrrrr}
-1& 1& 1& 1& 0 \\
0& -1& 0& 0& 1\\ 
1& -1& -1& 0& 0\\
1& 0& -1& -1& 0\\
1& 0& 0& -1& -1
\end{array} \right),&
K_5 = \left( \begin{array}{rrrrr}
2& -1& -1& -1& 0 \\
0& 1& 0& 0 &-1 \\ 
-1 & 1 & 1& 0& 0\\
-1& 0 & 1&  1& 0\\
-1& 0 & 0&  1& 1
\end{array} \right).
\end{array}
}
\end{equation}

\medskip
\noindent{\bf Rang 6 (principaux) : }
$\pm F_6,~G_6,~H_6,~\pm K_6,~\pm L_6,~\pm M_6,~\pm N_6$, avec 
\begin{equation}
\label{equa:indec_rg6}
{\small
\renewcommand{\arraycolsep}{0.2 em}
\begin{array}{ll}
F_6 = \left( \begin{array}{rrrrrr}
2 & 1 & 0 & 0 & 0 & 0 \\
1 & 2 & 0 & 0 & 1 & 1 \\
0 & 0 & 1 & 1 & 1 & 0 \\
0 & 0 & 1 & 1 & 0 & 1 \\
0 & 1 & 0 & 1 & 1 & 1 \\
0 & 1 & 1 & 0 & 1 & 1
\end{array} \right),&
\begin{array}{l}
G6 = \left( \begin{array}{cc}
0& G_3 \\
G_3& 0
\end{array} \right),\\
\\ 
H6 = \left( \begin{array}{cc}
0& H_3 \\
-H_3& 0
\end{array} \right),
\end{array} \\ \\
K_6 = \left( \begin{array}{rrrrrr}
-1 & 0 & 1 & 0 & -1 & -1 \\
0 & 0 & 0 & -1 & 0 & 0 \\
1 & 0 & -1 & 0 & 0 & 1 \\
1 & 0 & -1 & -1 & 0 & 0 \\
0 & -1 & 0 & -1 & -1 & 0 \\
-1 & 0 & 0 & 0 & -1 & -1
\end{array} \right),&
L_6 = \left( \begin{array}{rrrrrr}
2 & ~1 & -1 & ~0 & ~1 & 1 \\
-1 & 1 & 0 & 1 & 0 & 0 \\
-1 & 0 & 1 & 0 & 0 & -1 \\
-1 & 0 & 1 & 1 & 0 & 0 \\
0 & 1 & 0 & 1 & 1 & 0 \\
1 & 0 & 0 & 0 & 1 & 1
\end{array} \right),\\ \\
M_6 = \left( \begin{array}{rrrrrr}
1 & 0 & 0 & 0 & 0 & 1 \\
0 & 0 & 1 & -1 & 1 & 0 \\
-1 & 0 & 0 & 1 & -1 & 0 \\
0 & 1 & 0 & 0 & 1 & -1 \\
0 & -1 & 1 & 0 & 0 & 1 \\
0 & 1 & -1 & 1 & 0 & 0
\end{array} \right),&
N_6 = \left( \begin{array}{rrrrrr}
0 & 1 & 0 & 0 & 0 & -1 \\
-1 & 1 & -1 & 1 & -1 & 0 \\
1 & 0 & 0 & -1 & 1 & 0 \\
0 & -1 & 0 & 0 & -1 & 1 \\
0 & 1 & -1 & 0 & 0 & -1 \\
0 & -1 & 1 & -1 & 0 & 0
\end{array} \right).
\end{array}
}
\end{equation}

\medskip
\noindent{\bf Rang 7 (principaux) : }
$\pm F_7,~\pm G_7,~\pm H_7,~\pm K_7,~\pm L_7,~\pm M_7$, avec
\begin{equation}
\label{equa:indec_rg7}
{\small
\renewcommand{\arraycolsep}{0.2 em}
\begin{array}{ll}
F_7 = \left( \begin{array}{rrrrrrr}
~1 & ~1 & ~0 & ~1 &  ~0 &  ~0 &  ~~0\\
1 & 1 & 1 & 0 & 0 & 0 & 0\\
1 & 0 & 1 & 1 & 1 & 0 & 0\\
0 & 1 & 1 & 1 & 1 & 0 & 0\\
0 & 0 & 1 & 1 & 2 & 1 & 0\\
0 & 0 & 0 & 0 & 0 & 1 & 1\\
0 & 0 & 0 & 0 & -1 & -1 & 1
\end{array} \right),&
G_7 = \left( \begin{array}{rrrrrrr}
~1 &  ~1 &  ~0 &  ~1 &  ~0 &  ~0 &  ~0\\
1 & 1 & 1 & 0 & 0 & 0 & 0\\
1 & 0 & 1 & 1 & 1 & 0 & 0\\
0 & 1 & 1 & 1 & 1 & 0 & 0\\
0 & 0 & 1 & 1 & 1 & -1 & 0\\
0 & 0 & 0 & 0 & 0 & -1 & -1\\
0 & 0 & 0 & 0 & 1 & 1 & -1
\end{array} \right),\\
\\
H_7 = \left( \begin{array}{rrrrrrr}
-1 & -1 & 0 & 1 & 1 & 1 & 0 \\
0 & -1 & -1 & 1 & 0 & 0 & -1 \\
1 & 0 & -1 & 0 & 0 & -1 & -1 \\
0 & 1 & 1 & -1 & -1 & 1 & 1 \\
1 & 1 & 0 & -1 & -1 & -1 & 1 \\
1 & 0 & 0 & -1 & -1 & -1 & -1 \\
1 & 0 & -1 & 1 & -1 & -1 & -1
\end{array} \right),&
K_7 = \left( \begin{array}{rrrrrrr}
1 & 0 & 0 & -1 & -1 & -1 & 0 \\
0 & 0 & 0 & -1 & 0 & 0 & 1 \\
0 & 1 & 0 & 0 & 0 & 1 & 1 \\
0 & -1 & -1 & 1 & 1 & -1 & -1 \\
-1 & -1 & 0 & 1 & 1 & 1 & -1 \\
-1 & 0 & 0 & 1 & 1 & 1 & 1 \\
-1 & 0 & 1 & -1 & 1 & 1 & 1
\end{array} \right),\\
\\
L_7 = \left( \begin{array}{rrrrrrr}
0 & 0 & 0 & 0 & 1 & 0 & 0 \\
0 & 0 & 0 & 0 & 1 & -1 & 0 \\
0 & 0 & 0 & 1 & -1 & 0 & 0 \\
0 & 0 & 0 & -1 & 1 & -1 & 1 \\
0 & 0 & 1 & -1 & -1 & 1 & -1 \\
1 & 1 & -1 & 1 & -1 & -1 & 1 \\
0 & -1 & 0 & -1 & 1 & -1 & -1
\end{array} \right),&
M_7 = \left( \begin{array}{rrrrrrr}
1 & 1 & 0 & 0 & -1 & 0 & 0 \\
0 & 1 & 1 & 0 & -1 & 1 & 0 \\
-1 & -1 & 1 & -1 & 1 & 0 & 0 \\
0 & 0 & 0 & 1 & -1 & 1 & -1 \\
0 & 0 & -1 & 1 & 1 & -1 & 1 \\
-1 & -1 & 1 & -1 & 1 & 1 & -1 \\
0 & 1 & 0 & 1 & -1 & 1 & 1
\end{array} \right).
\end{array}
}
\end{equation}

\subsection{Types géométriques maximaux indécomposables}
\label{subs:ann_tg_non_scin}
Afin de compléter l'énoncé du théorème~\ref{theo:intro_param1_7} de 
l'introduction, nous explicitons ici les types géométriques maximaux
non scindés en rang $\leq 7$, autres que les types symplectiques ou 
orthogonaux. On rappelle que $\mfS_g$, $V_{p,q}$, $V_{2,2}^\mathrm{II}$ 
et $W_{p,q}$ sont  définis par \eqref{equa:Sig_g}, \eqref{equa:Vpq},
\eqref{equa:V22_pair} et \eqref{equa:Wpq} (voir aussi \eqref{equa:W11}
pour $W_{1,1}$).

\medskip
\noindent{\bf Rang 5}
$$ 
{\small 
V_{I_1F_4} = P\cdot  (V_{2,1}\oplus \mfS_1)  
 ~~~ \mathrm{avec} ~~~
P = \frac{1}{2}
\begin{pmatrix}
2 & 2 & -2 & 0 & 0 \\
2 & 0 & -2 & -2 & 0 \\
2 & 0 & -2 & 2 & 0 \\
0 & -1 & 1 & 0 &  -1 \\
0 & -1 & 1 & 0 & 1
\end{pmatrix}
}
$$
(voir également \eqref{equa:tg_I1F4}).

\medskip
\noindent{\bf Rang 6}
$$
V_{I_2 F_4} = P\cdot  (V_{3,1}\oplus \mfS_1), ~~
V_{I_{1,1}F_4} = P\cdot  (V_{2,2}\oplus \mfS_1), ~~
V_{U_2F_4} = Q \cdot (V_{2,2}^\mathrm{II}\oplus \mfS_1),
$$
avec
$$
{\small
P = \frac{1}{2}
\begin{pmatrix}
0 & 2 & 0 & 0 & 0 & 0 \\
0 & 0 & 2 & 0 & 0 & 0 \\
1 & 0 & 0 & -1 & -1 & 0 \\
1 & 0 & 0 & -1 & 1 & 0 \\
1 & 0 & 0 & 1 & 0 & -2 \\
1 & 0 & 0 & 1 & 0 & 2
\end{pmatrix}
~~~ \mathrm{et} ~~~
Q = \frac{1}{2}
\begin{pmatrix}
2 & 0 & 0 & 0 & 0 & 0 \\
0 & 0 & 2 & 0 & 0 & 0 \\
0 & 0 & 0 & 2 & 0 & 1 \\
0 & 0 & 0 & 2 & 0 & -1 \\
0 & 1 & 0 & 0 & -2 & 0 \\
0 & 1 & 0 & 0 & 2 & 0
\end{pmatrix}.
}
$$

\medskip
\noindent{\bf Rang 7}
$$
\begin{array}{ll}
V_{I_3 F_4} = P_1\cdot (V_{4,1}\oplus \mfS_1),  & 
V_{I_{2,1} F_4} = P_1\cdot (V_{3,2}\oplus \mfS_1), \\
V_{I_1 J_2 F_4} = P_2\cdot (V_{2,1}\oplus \mfS_2), &
V_{I_1 G_3 G_3^-} = P_3\cdot (V_{2,1}\oplus W_{1,1}), \\
V_{K_4 G_3}= P_4\cdot (\{I_1\}\oplus W_{2,1}), & 
\end{array}
$$
avec
$$
{\small
P_1 = \frac{1}{2}
\begin{pmatrix}
2 & 0 & 0 & 0 & 0 & 0 & 0 \\
0 & 2 & 0 & 0 & 0 & 0 & 0 \\
0 & 0 & 0 & 2 & 0 & 0 & 0 \\
0 & 0 & -1 & 0 & -1 & -\sqrt{2} & 0 \\
0 & 0 & -1 & 0 & -1 & \sqrt{2} & 0 \\
0 & 0 & -1 & 0 & 1 & 0 & -\sqrt{2} \\
0 & 0 & -1 & 0 & 1 & 0 & \sqrt{2}
\end{pmatrix},
~~~
P_2 = \frac{1}{2}
\begin{pmatrix}
2 & 0 & 0 & 0 & 0 & 0 & 0 \\
0 & 0 & 0 & 0 & 2 & 0 & 0 \\
0 & 0 & 0 & 0 & 0 & 0 & 2 \\
0 & 1 & 1 & -\sqrt{2} & 0 & 0 & 0 \\
0 & 1 & 1 & \sqrt{2} & 0 & 0 & 0 \\
0 & 1 & -1 & 0 & 0 & -\sqrt{2} & 0 \\
0 & 1 & -1 & 0 & 0 & \sqrt{2} & 0
\end{pmatrix},
}
$$
$$
{\small
P_3 = \frac{1}{2^{7/4}}
\begin{pmatrix}
2^{7/4} & 0 & 0 & 0 & 0 & 0 & 0 \\
0 & 2^{5/4} & 0 & 1 & -1 & -1 & -3 \\
0 & 0 & 0 & 0 & -2 & -4 & -2 \\
0 & 0 & 0 & -2 & 0 & -2 & 4 \\
0 & 0 & 2^{5/4} & -3 & -1 & -1 & 1 \\
0 & 0 & 0 & -4 & 2 & 0 & 2 \\
0 & 0 & 0 & 2 & 4 & 2 & 0
\end{pmatrix}
~~~\mathrm{et}
}
$$
\begin{equation}
\label{equa:ann_P4}
{\small
P_4 = \frac{1}{2}
\begin{pmatrix}
0 & 0 & 0 & \sqrt{2} & 0 & -\sqrt{2} & 0 \\
0 & 0 & 0 & 0 & -\sqrt{2} & 0 & \sqrt{2} \\
0 & 0 & 0 & -2 & 0 & 0 & 2 \\
0 & 0 & 0 & 0 & 2 & 2 & 0 \\
\sqrt{2} & -2^{1/4} & 2^{1/4} & 0 & 0 & 0 & 0 \\
0 & 0 & 2^{5/4} & 0 & 0 & 0 & 0 \\
0 & 2^{5/4} & 0 & 0 & 0 & 0 & 0
\end{pmatrix}.
}
\end{equation}

\vspace{-11pt}

\phantomsection{}
\addcontentsline{toc}{section}{Bibliographie}

\bibliographystyle{alphabav}
\bibliography{bavard}

\noindent {\small 
Institut de Mathématiques de Bordeaux\\
U.M.R. 5251 C.N.R.S.\\
Université Bordeaux\\
351, cours de la Libération\\
F-33405 TALENCE Cedex}

\end{document}